\documentclass[11pt,a4paper,maplems]{article}
\usepackage{color,amsmath,amssymb, lscape,graphicx,ifthen,pdfpages}
\DeclareSymbolFont{boldsymbols}{OMS}{cmsy}{b}{n} %% for bf  \mathcal
\DeclareSymbolFontAlphabet{\mathbfcal}{boldsymbols} %% for bf  \mathcal
\usepackage[colorlinks=true,
linkcolor=webgreen,
filecolor=webbrown,
citecolor=webgreen]{hyperref}
\definecolor{webgreen}{rgb}{0,.5,0}\definecolor{webbrown}{rgb}{.6,0,0}
\definecolor{Bittersweet}{cmyk}{0,0.75,1,0.24}

\newcommand{\seqnum}[1]{\href{http://oeis.org/#1}{\underline{#1}}}
\def\dstyle#1{$\displaystyle #1 $}

\def\pn{\par\noindent}
\def\ps{\par\smallskip}
\def\psn{\par\smallskip\noindent}
\def\pbn{\par\bigskip\noindent}
\def\Beq{\begin{equation}}
\def\Eeq{\end{equation}}
\def\Beqarray{\begin{eqnarray}}
\def\Eeqarray{\end{eqnarray}}
\def\sspgeq{\,\geq} 
\def\sspleq{\, \leq \,}
\def\sspkl{\, < \,}
\def\sspgr{\, > \,}
\def\sspeq{\, =\,}
\def\speq{\ =\ }
\def\sspdef{\, :=\,}
\def\spdef{\ :=\ }
\def\sspfed{\, =:\,}
\def\sspin{\, \in \,}

\def\sspp{\, +\ }
\def\sspm{\, -\ }
\def\sppm{\ \pm\ }
\def\ssppm{\,\pm\,}
\def\sspto{\,\to\,}

\def\sspcdot{\,\cdot\,}
\def\sspneq{\, \neq \,}
\def\sspequiv{\,\equiv\,}
\def\sspnotequiv{\,\not\equiv  \,}
\def\Aequ{\,\Leftrightarrow\, }
\def\Aequiv{\,\Leftrightarrow\, }
\def\sspconj{\,\buildrel ! \over = \,}
\def\binomial#1#2{{#1} \choose {#2}}
\def\Modd{$\rlap{\it M}\ \ \, odd $}
 
\def\Moddn#1{(\rlap{\it M}\ \ \, odd\, $#1$ ) } 
\def\ogf{{\it o.g.f.\ }}

\def\ie{{\it i.e.},\, }
\def\eg{{\it e.g.},\, }
\def\Eg{{\it E.g.},\, }
\def\viz{{\it viz}\, }

\def\lhs{{\it l.h.s.\, }}
\def\rhs{{\it r.h.s.\, }}
\def\set#1#2{\left\{\,#1\, \Biggm | \, #2\, \right\}}
\def\sspunion{\,\cup\,} 
\def\sspiso{\,\cong\,}
\def\Caseszwei#1#2#3#4{\left\{ \begin{array}{ll}#1&\mbox{#2}\\ &\\#3&\mbox{#4}\end{array}\right.}
\def\Cases3#1#2#3#4#5#6{\left\{ \begin{array}{ll}#1&\mbox{#2}\\#3&\mbox{#4}\\#5&\mbox{#6}\end{array}\right.}
\def\sspeqv{$\,\equiv\,$}
\def\floor#1{\left\lfloor{#1}\right\rfloor}
\def\ceil#1{\left\lceil{#1}\right\rceil}
\def\spfed{\ =:\ }% space def reversed order
\def\sspfed{\, =:\, }% space def reversed order
\def\sspleftrightarrow{\,\leftrightarrow\, }
\def\sspmapsto{\,\mapsto\,}
\def\sspfollows{\,\Rightarrow\,}
\def\sspdoteq{\,\doteq\, }
\def\Simn{\,{\lower2pt\hbox{$\buildrel {\lower3pt\hbox{$n$}} \over \sim$}}\,}
\def\simn1#1{\,{\lower2pt\hbox{$\buildrel {\lower3pt\hbox{$#1$}} \over \sim$}}\,}
\def\sspnotin{\,\notin\,}
\def\Chi{\raise2.5pt\hbox{$\chi$}} 
\def\ssptimes{\,\times\,}
\normalbaselines
\begin{document}
\bibliographystyle{unsrt}
\rightline{Karlsruhe} \par\smallskip\noindent
\rightline{October 03  2012}
\rightline{Revised: March 03 2017}
\vbox {\vspace{6mm}}
\begin{center}
{\Large {\bf The field $\bf \mathbb Q(2\, cos\left(\frac{\boldsymbol{\pi}}{n}\right))$, its Galois group, and length ratios in the regular $\bf n$-gon}}\\ [9mm]
Wolfdieter L a n g \footnote{ 
\href{mailto:wolfdieter.lang@partner.kit.edu}{\tt wolfdieter.lang@partner.kit.edu},\quad 
\url{http://https://www.itp.kit.edu/~wl}
                                          } \\[3mm]
\end{center}
\vspace{2mm}
\begin{abstract}
\par\smallskip\noindent
%%\vskip 1cm  \par \noindent
The normal field extension $\mathbb Q(\rho(n))$,  with the algebraic number  \dstyle{\rho(n)\spdef  2\,\cos\left ({\frac{\pi}{n}}\right)}, for $n\in \mathbb N$, is related to ratios of the lengths between diagonals and the side of a regular $n-$gon.  This has been considered in a paper by P. Steinbach. These ratios, numbered $k\sspeq 1, ..., n-1$, are given by {\sl Chebyshev} polynomials $S(k-1,x=\rho(n))$. The product formula for these ratios was found by {\sl Steinbach} and is re-derived here from a known formula for the product of {\sl Chebyshev} $S-$polynomials. It is shown that it follows also from the $S-$polynomial recurrence and certain rules following from  the trigonometric nature of the argument $x\sspeq \rho(n)$.  The minimal integer polynomial $C(n,x)$ for $\rho(n)$ is presented, and its simple zeros are expressed in the power-basis of $\mathbb Q(\rho(n))$. Also the positive zeros of the {\sl Chebyshev} polynomial $S(k-1,\rho(n))$ are rewritten in this basis. The number of positive and negative zeros of $C(n,x)$ is determined. The coefficient $C(n,0)$ is computed for special classes of $n$ values. Theorems on $C(n,x)$ in terms of  monic integer {\sl Chebyshev} polynomials of the first kind (called here $\hat t$) are given. These polynomials can be factorized in terms of the minimal $C$-polynomials. A conjecture on the discriminant of these polynomials is made.  The {\sl Galois} group is either $Z_{\delta(n)}$, the cyclic group of the order given by the degree $\delta(n)$ of $C(n,x)$, or it is a direct product of certain cyclic groups. In order to determine the cycle structure  a novel modular multiplication, called \Modd\, $n$ is introduced. On the  reduced odd residue system \Modd\, $n$ this furnishes a group which is isomorphic to this {\sl Galois} group.
\end{abstract}
\section{Introduction and Summary}
Length ratios between diagonals and the side of a regular $n-$gon (called diagonal/side ratios, abbreviated DSRs) have been considered by Steinbach \cite{Steinbach}. He gave a product formula for these ratios (called by him diagonal product formula (DPF)), and for the pentagon and heptagon details were given. In the pentagon case the quadratic number field $\mathbb Q(\sqrt{5})$ with the basis  $<1,\varphi>$ for integers of this field turns up. Here $\varphi$, the golden section, is identified with $\rho(5)$ which is the length ratio between any of the two diagonals and the pentagon side. For the heptagon there are two different diagonal length and the two ratios between them and the side length have been called $\rho:=\rho(7)$ and $\sigma$. The DPF allowed to reduce all products and quotients of  $\rho$ and $\sigma$ to $\mathbb Q$-vectors with the basis $<1,\rho,\sigma>$ . For example, $\rho^2 \sspeq 1\sspp \sigma$, and one can use instead the power basis  $<1,\rho,\rho^2>$ of the algebraic number field $\mathbb Q(\rho(7))$ (the use of the same acute bracket notation for different bases should not lead to a confusion). The minimal polynomial of  $\rho(7)$, in the present work called C(7,x), is $x^3\sspm x^2\sspm2\,x \sspp 1$, and this was also given in \cite{Steinbach}. Therefore,  $[\mathbb Q(\rho(7)):\mathbb Q]\sspeq 3$, which is the degree of this field extension. For odd $n$, $n=2\,k+1$, Steinbach gave a polynomial with one of its roots $\rho(2\,k+1)$. In the present work the minimal polynomial for $\rho(n)$, called C(n,x), for every $n\in \mathbb N$, is given in terms of the known one for \dstyle{\cos\left ({\frac{2\,\pi}{N}}\right )}, called $\Psi(N,x)$, with $N=2\,n$. The connection between $\Psi(N,x)$ and divisor product representations of $N$, and  to {\sl Chebyshev} $T-$polynomials has been worked out earlier by the author. See the links under OEIS \cite{Sloane} \seqnum{A181875} and  \seqnum{A007955} (in the sequel OEIS Anumbers will appear without repeating this reference).  $C(n,x)$ turns out to be an integer polynomials of degree $\delta(n)\sspeq$\seqnum{A055034}$(n)$. All its zeros are known, and they are simple. $\mathbb Q(\rho(n))$ is  the splitting field for the minimal polynomial $C(n,x)$ of $\rho(n)$. It is a normal extension of the rational field. These $C$-polynomial zeros are written in  the power basis of  $\mathbb Q(\rho(n))$ with the help of certain scaled {\sl Cheyshev} $T-$polynomials (called here $\hat t$ with their integer coefficient array \seqnum{A128672}). This array can also be used to rewrite the positive zeros of the {\sl Chebyshev} $S(n-1,x)$ polynomial, related to the DSRs, in this power basis. The use of $S-$polynomial technology allows also for a re-derivation of the product formulae (DPF) for the length ratios between the regular $n-$gon diagonals and side.\psn
In section $2$  $S(k-1,\rho(n))$ is shown to yield the DSRs of the regular $n$-gon. The DPF is also re-derived there, and the independent products are extracted. This is called reduced algebra over $\mathbb{Q}$.  In section $3$ the minimal polynomial $C(n,x)$ is presented, and its zeros are related to the power basis of the algebraic number field  $\mathbb Q(\rho(n))$. In section $4$ the $\mathbb Q$-automorphisms of this field, the {\sl Galois} group $G(\mathbb Q(\rho(n))/\mathbb Q)$, is treated. It is the cyclic group $Z_{\delta(n)}$, except for infinitely many $n$-values where it is the direct product of cyclic groups, hence non-cyclic. For $n$ from $1..100$ there are $30$ non-cyclic groups. A novel modular multiplication on the odd numbers, called \Modd\, $n$, is introduced which serves to determine the cycle graph structure of this {\sl Galois} group. The multiplicative group \Modd\, $n$ is isomorphic to this {\sl Galois} group. On the reduced odd residue system \Modd\, $n$ this multiplicative group is isomorphic to the {\sl Galois} group. 
%%%%%%%%%%%%%%%%%%%%%%%%%%%%%%%%%%%%%%%%%
\section{Regular polygon diagonals/side ratios (DSRs)}
%%%%%%%%%%%%%%%%%%%%%%%%%%%%%%%%%%%%%%%%%
Motivated by a paper of Steinbach \cite{Steinbach} we consider the diagonals  and the side in the regular $n-$gon (inscribed in the unit circle). The vertices on the unit circle are called $V^{(n)}_k$, $k=0, ...., n-1,\ n\sspeq 2,3,...$; $V^{(n)}_0$ has coordinates $(1,0)$;  the side length (with the radius' length 1 unit) is \dstyle{s(n)\sspeq 2\,\sin\left ({\frac{\pi}{n}}\right)}. The length of $\overline {V^{(n)}_0\,V^{(n)}_k}$ is \dstyle{d^{(n)}_k\sspeq 2\,\sin \left ({\frac{\pi\,k}{n}}\right)}, for $k=1, ..., n-1$, hence $d_1\sspeq s(n)$. We only need one side and the diagonals in the upper half plane, {\it i.e.} it suffices to consider $k\in\{1, ... ,\floor{n/2}\}$. For even $n$, the largest diagonal (\dstyle{k\sspeq \frac{n}{2}}), of  length $2$, lies on the negative real axis. See {\it Fig. 1}, and {\it Fig. 2}, for the case $n\sspeq 10$ (decagon),  and for $n\sspeq 11$ (hendecagon), respectively. The length ratios of interest, the DSRs,  are \dstyle{R^{(n)}_k\sspeq {\frac{d^{(n)}_k}{d^{(n)}_1 }}} for the given $k-$values. Here {\sl Chebyshev} $S-$polynomials in their trigonometric form enter the stage (see \seqnum{A049310} for their coefficients):
\Beq
R^{(n)}_k\sspeq S(k-1,\rho(n)),
\Eeq
where we use the second ratio (for the smallest diagonal) \dstyle{\rho(n)\sspdef R^{(n)}_2\sspeq 2\, \cos\left ({\frac{\pi}{n}}\right)}.
Remember that $S(n,x)\sspdef U(n,x/2)$ with {\sl Cheyshev} $U-$polynomial (second kind).  It is well known that the zeros of the polynomials $S(n-1,x)$ are \dstyle{ x^{(n-1)}_{k,\pm}\sspeq \pm 2\, \cos \left ({\frac{k\,\pi}{n}}\right)}, for $k\sspeq 1,\, ...,\, \ceil{{\frac{n-1}{2}}}$. Note that in the even $n$ case the zero $x=0$ appears twice from $\pm 0$. The $S-$polynomials are orthogonal polynomials defined on the real interval $[-2,+2]$.  Hence their zeros are guaranteed to be simple. They belong to the {\sl Jacobi} class. Due to (anti)symmetry it is sufficient to consider only the positive zeros (we disregard $x=0$ in the even $n$ case). Now also {\sl Chebyshev} $T-$polynomials (first kind; see the coefficient table \seqnum{A053120}) enter, which are written in terms of $S-$polynomials. (The $T-$polynomials appear as trace polynomials in a $2\times2$ transfer matrix when one solves the three term recurrences of these classical orthogonal polynomials.)
\Beqarray
x^{(n-1)}_{k}\sspequiv  x^{(n-1)}_{k,+}&\sspeq& 2\,T(k,\rho(n)/2)\sspfed {\hat t}_k(\rho(n))\sspequiv {\hat t}(k,\rho(n))\\
&\sspeq& S(k,\rho(n))\sspm S(k-2,\rho(n))\, , \ k=1, ...,\floor{{\frac{n-1}{2}}}\, .  
\Eeqarray
Note that these positive zeros decrease with $k$, and that ${\hat t}_k$ is an integer polynomial of degree $k$ because the $S-$polynomials are integer. The coefficient array for these monic  $\hat t$-polynomials is \seqnum{A127672}. In \cite{ASt} these polynomials are called {\sl Chebyshev} $C$-polynomials, but we use the letter $C$ for the minimal polynomials of $\rho(n)$. The zeros can be replaced by DSRs from the above given eq. $(1)$. \psn
\Beq
\text{\bf Positive zeros of}\ {\bf S(n-1,x)} \ \text{\bf from DSRs:} \quad\quad  x^{(n-1)}_{k}\sspeq R^{(n)}_{k+1}\sspm R^{(n)}_{k-1}\  . \hskip 4cm 
\Eeq
%%%%%%%%%%%%%%%%%%%%%%
%%% 2 Figures 
\parbox{16cm}{
{\includegraphics[height=8cm,width=.5\linewidth]{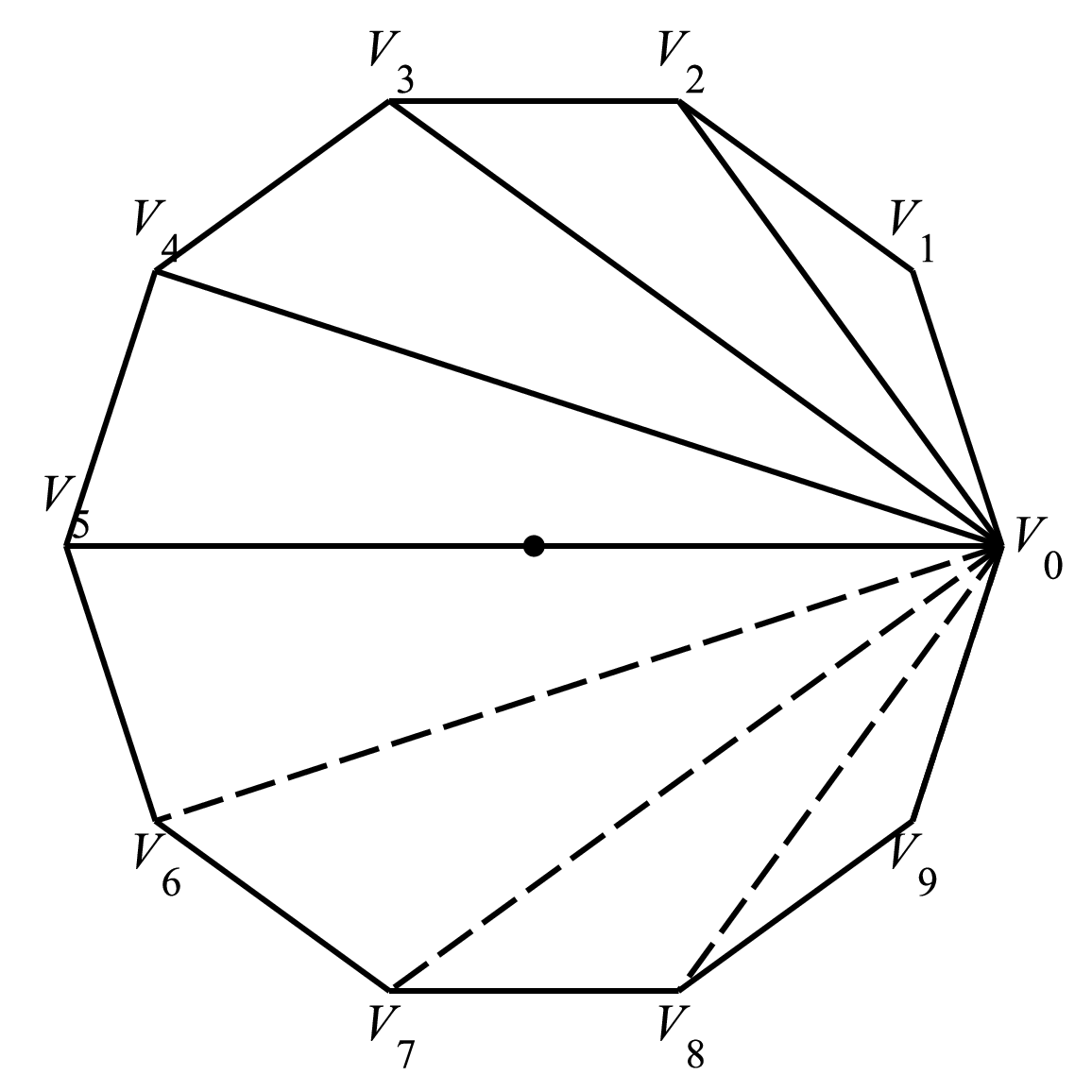}}
{\includegraphics[height=8cm,width=.5\linewidth]{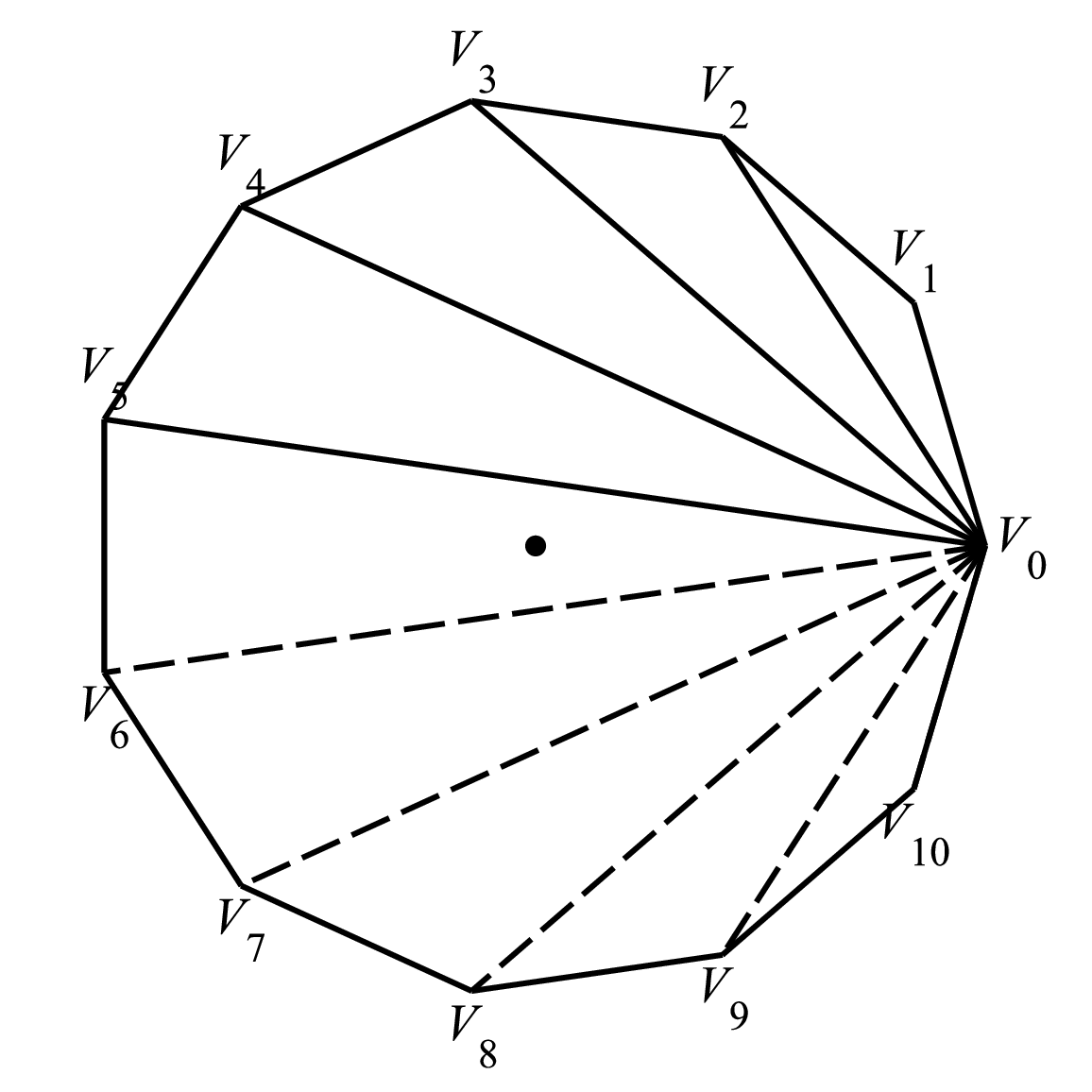}}
}
\psn
Figure 1: Decagon $n\sspeq 10$\hskip 5cm Figure 2: Hendecagon (n=11)
\psn
%%%%%%%%%%%%%%%%%%%%%%
As an aside we give the factorization of $S(n-1,x)$ in terms of the polynomial \pn
\dstyle{P(k,x)\sspdef \prod_{l=1}^k\,\left (x\sspm 2\, \cos \left ({\frac{l\,\pi}{2\,k+1}}\right)\right)} with positive zeros.
\Beq
S(n-1,x)\sspeq x^{\Theta(n-1)}\, P\left(\floor{{\frac{n-1}{2}}},x\right)\,(-1)^{\floor{{\frac{n-1}{2}}}}\, P\left(\floor{{\frac{n-1}{2}}},-x\right)\ ,
\Eeq 
where $\theta(n-1)\sspdef 0$ if $n$ is odd, and $1$ if $n$ is even. 
This factorization is also considered in \cite{Chan} for the $U(n,x)$ polynomials for even and odd $n$ separately. For example, for the heptagon $n\sspeq 7$ one finds the factorization of $S(6,x)$ with the following real polynomial $P(3,x)$. 
\Beq 
P(3,x)\sspeq x^3\sspm (2\, \sigma(7)-1)\, x^2\sspp 2\,\rho(7)\, x\sspm 1 \ ,
\Eeq 
with $\rho(7)\sspeq R^{(7)}_2$ and $\sigma(7)\sspdef R^{(7)}_3$ (see also the introduction). The three positive zeros are, written in terms of the DSRs, $\rho(7),\, \sigma(7)\sspm 1$ and $\sigma(7)\sspm \rho(7)$. 
Here we used the general translation formulae eq. $(4)$ to rewrite the positive zeros of $S(n-1,x)$  in terms of  DSRs. One should also write this in terms of powers of $\rho(7)$ with the help of eq. $(1)$ (see eq. $(7)$ below).\pn
Each \dstyle{P\left(\floor{{\frac{n-1}{2}}},x\right)} can be considered as characteristic polynomial for a ($\floor{{\frac{n-1}{2}}}\sspp 1$)-term recurrence sequence $\{f_k\}_{k=0}^{\infty}$ in the integral domain of $\mathbb Q(\rho(n))$ (represented by a $\delta(n)$-tuple of ordinary (rational) integers $(A_{1,k},A_{2,k},..., A_{\delta(n),k})$. Here $\delta(n)$ is the degree of the minimal polynomial of  $\rho(n)$ which will be discussed in sect. $3$. In fact,  one has $\floor{{\frac{n-1}{2}}}$ such $\delta(n)$-tuple sequences corresponding to the independent inputs. As a simple example take the $n=5$ pair ($\delta(5)\sspeq 2$) of sequences corresponding to  $P(2,x)\sspeq x^2\sspm (2\,\varphi\sspm 1)\,x \sspp 1$, with the golden section $\varphi:=\rho(5)$ and the simplest input. This is $S(k,\sqrt{5})\sspeq A_{1,k}\, 1\sspp  A_{2,k}\, \varphi$ with input $S(-1,\sqrt{5})\sspeq 0,\ S(0,\sqrt{5})\sspeq1$. See \seqnum{A005013}$(k+1)\, (-1)^k$, and  $2\,$\seqnum{A147600}$(k-1)$, for $A_{1,k}$ and $A_{2,k}$, respectively.
\psn
For later purposes it is also useful to give a dictionary between the positive zeros of $S(n-1,x)$, {\it i.e.} those of $P(\floor{{\frac{n-1}{2}}},x)$,  and powers of $\rho(n)$. One direction follows from the above given formula, eq. $(2)$, for $x^{(n-1)}_{k}$ in terms of the monic integer polynomials ${\hat t}(k,\rho)\sspeq 2\,T(k,\rho/2)$. The dependence on $n$ has here been omitted.\psn
\Beq
\text{\bf Positive zeros of }\  {\bf S(n-1,x)}\  \text{\bf from} \ {\bf {\boldsymbol \rho}(n)}\text{\bf -powers:}\quad\quad x^{(n-1)}_{k} \sspeq \sum_{l=0}^{\floor{{\frac{k}{2}}}}\,(-1)^l\,{\binomial{k-l}{l}}\, \rho(n)^{k-2l}\ .\hskip .5cm
\Eeq
The coefficients are given in the triangle \seqnum{A127672}. \Eg $x_5 \sspeq 5\, \rho \sspm 5\, \rho^3\sspp \rho^5$ (omitting $n$).   For the {\sl heptagon} this shows that the three positive zeros of $S(6,x)$ are $\rho(7), \ \sigma(7)-1\sspeq \rho^2(7)\sspm 2$  and $\sigma(7)\sspm \rho(7)\sspeq \rho(7)^2 \sspm \rho(7)\sspm 1\, $. In this case these three zeros could be used as a $\mathbb Q$-vector space basis. Note that the $n$ dependence is {\it via} the ${\hat t}$-polynomial variable $\rho(n)$ (which satisfies $C(n,\rho(n))=0$ by definition of the minimal polynomial $C$). For example, it is true for all $n$ that the second largest positive zero of $S(n-1,x)$, namely $x^{(n-1)}_2$, is always $\rho(n)^2 \sspm 2$. Hence, from above, $R^{(n)}_3\sspm 1 \sspeq \rho(n)^2\sspm 2$; or if one calls $R^{(n)}_3\spfed \sigma(n)$ then $\rho(n)^2\sspeq \sigma(n)\sspp 1$. This has been noted for $n=7$ above but is holds in general, showing that $\sigma(n)$ is algebraically dependent on $\rho(n)$. This fact has been observed already in \cite{Steinbach} where it appeared as a special product formula (see the later DPF eq. $(13)$, $m\sspeq k\sspeq 2$, or eq. $(11)$ with $k\sspeq 2$).  \psn
The inverse formulae (DSRs in terms of $\rho$-powers) have already been given above in eq. $(1)$. Here one uses the $S$-triangle \seqnum{A049310} for the translation. \Eg $R^{(n)}_4\sspeq -2\, \rho(n)$ $\sspp \rho^3(n)$ $\sspeq \rho(n)\,(\rho(n)^2\sspm 2)$ $\sspeq \rho(n)\,(\sigma(n)\sspm 1)$\ . Of course, $R^{(n)}_4$ is interesting only for $n\sspgeq 9$ because only then the corresponding diagonal lies in the upper half-plane (including the negative real axis).\psn
Note that the number of positive zeros of $S(n-1,x)$, \ie $\floor{{\frac{n-1}{2}}}$, is less or equal to $\floor{{\frac{n}{2}}}$, the number of DSRs for diagonals in the upper half plane including the negative real axis. The zero $x\sspeq 1$ never appears. As can be seen in the $n=7$ case these zeros can nevertheless be used as $\mathbb Q$-vector space basis (in sect. 3 it will become clear that the degree $\delta(7)$ is also 3). In the case $n\sspeq 9$ the two numbers are both 4, but the degree $\delta(9)$ is $3$, hence only  three of the zeros and three of the DSRs are rationally independent. To wit: $x^{(8)} \sspeq -5\,x^{(8)}_3\sspp 3\,x^{(8)}_1\sspm 4\,x^{(8)}_2 $ and $R^{(9)}_4\sspeq R^{(9)}_1 \sspp R^{(9)}_2 \sspeq 1\sspp \rho(9)$. \pbn
Because of symmetry only DSRs for diagonals of the upper half plane (including the negative real axis), \ie \dstyle{k\sspin \{0,1, ..., \floor{\frac{n}{2}}\}} , are of interest. The reduction for other $k$ values, also negative ones, is accomplished by the rules
\Beq
\text{\bf o)}\ \ R^{(n)}_{n+k}\sspeq - R^{(n)}_{k} \ \  , \ \  \text{\bf i)}\ \  R^{(n)}_{-|k|}\sspeq - R^{(n)}_{|k|}\ \ , \ \ \text{\bf ii)}\ \ {\text {\bf for}} \ k\sspin \{\floor{\frac{n}{2}}+1,\,...,\,n\}\ :\ \ R^{(n)}_k\sspeq R^{(n)}_{n-k}\ .
\Eeq
These rules follow  from eq. $(1)$ with the trigonometric definition (with the specific value of $\rho(n)$) of the $S$-polynomials, and therefore  also negative values for the DSRs  show up (the interesting DSRs are, of course, positive). $\bf i)$ follows also from $S(-|n|,x)\sspeq -S(|n|-2,x)$ which derives from a backward use of the recurrence (given later in eq. $(12)$). This rule does therefore not depend on the special choice of the variable $x$, in contrast to the rules $\bf o)$ and $\bf ii)$. With negative $k$ one counts the diagonals in the clockwise direction.  {\bf ii)} is used to translate  from the lower to the upper half-plane.\psn
The product formula for {\sl Chebyshev} $S$-polynomials of different degree, but with the same argument, is well known. See \eg \cite{ASt}, p.782, 22.7.25, for the $U$-polynomials, and replace the $T$-polynomials by the $S$-polynomials  {\it via} the trace formula given in eqs $(2)$ and $(3)$.
\Beqarray
S(m-1,x)\, S(n-1,x)\sspeq {\frac{1}{4\,(({\frac{x}{2}})^2\sspm 1)}}&& \left[ S(n+m,x)\sspm S(n+m-2,x) \sspm \right. \nonumber \\
 && \left. S(n-m,x)\sspp S(n-m-2,x)\,\right]\  \ n\sspgeq m\ .
\Eeqarray
 For our purpose $x\sspeq \rho(n)\sspeq R^{(n)}_2 \sspneq \pm 2$ but its special form is for this formula not of interest .  It is in fact symmetric under exchange of $n$ with $m$ if the rules for negative indices on the $S$-polynomials, stated above, are employed. Here we restrict to $n\sspgeq m$ in order to have non-negative indices.
Thus the product formula for the DSRs, now with $x\sspeq \rho(n)$, but without using its specific value, is
\Beq
(4\sspm (R^{(n)}_2)^2)\, R^{(n)}_m\,R^{(n)}_k\sspeq R^{(n)}_{k-m+1}\sspm R^{(n)}_{k-m-1}\sspm R^{(n)}_{k+m+1} \sspp R^{(n)}_{k+m-1}\ ,\ k\sspgeq m\sspgr 0\ .
\Eeq
For $m\sspeq 0$ this becomes trivial. This is not yet the DPF given by {\sl Steinbach} in \cite{Steinbach} which linearizes the product of $R^{(n)}_m\,R^{(n)}_k$. In order to eliminate the pre-factor $(4\sspm \rho(n)^2)$ one can use the three term recurrence relation of the orthogonal $S$-polynomials, written for the DSRs by eq. $(1)$ as a special product formula
\Beq
{\bf [Rec,k]^{(n)}}\ :\hskip 2cm R^{(n)}_2\, R^{(n)}_k \sspeq R^{(n)}_ {k+1} \sspp R^{(n)}_{k-1}\ ,\ \ k\sspgeq 0\ .
\Eeq
From now on we also use the specific form of $\rho(n)$. Thus also rules $\bf o)$ and $\bf ii)$  of eq. $(8)$ will be applicable. For $k\sspeq 0$ one uses (see {\bf ii)}) \dstyle{R^{(n)}_{-1}\sspeq -\,R^{(n)}_{1} \sspeq -1}, and it becomes trivial. For $k=2$  it shows that $\rho(n)^2\sspeq \sigma(n)\sspp 1$ for all $n$, if one uses \dstyle{\sigma(n)\sspdef R^{(n)}_3 }.  This recurrence can now be used twice  as \dstyle{R^{(n)}_2\, (R^{(n)}_2\, R^{(n)}_k)} in eq. $(10)$ to produce the following recurrence for two step differences of \dstyle{R^{(n)}_m\, R^{(n)}_k \sspfed p^{(n,m)}_k\sspequiv p_k}
\Beq
(p_{k+2}\sspm p_k) \sspm (p_k\sspm p_{k-2})\sspeq c_k\sspm c_{k-2}\ ,
\Eeq
where we used the  abbreviation \dstyle{c_k\sspequiv c^{(n,m)}_k\sspdef  R^{(n)}_{k+m+1}\sspm R^{(n)}_{k-m+1}}. This shows that \dstyle{p_{k+2}\sspm p_k \sspm c_k} is $k$ independent. From the inputs $p_2\sspm p_0\sspm c_0\sspeq 0$, due to the recurrence relation  eq. $(11)$, and $p_1\sspm p_{-1}\sspm c_{-1}\sspeq 0$ with the help of the rule {\bf i)} from eq. $(8)$, this leads to  the recurrence \dstyle{p_{k+2}\sspm p_{k} \sspeq c_k}, with the inputs $p_0\sspeq 0$ and $p_1\sspeq -p_1\sspeq R^{(n)}_m$. The solution for $p_k$ is found for even and odd $k$ separately, where again the rules for negative indices are employed. Both solutions can then be combined as
\Beq
\text {\bf DPF}\ \  {\bf [m,k]^{(n)}}:\hskip 2cm  R^{(n)}_m\, R^{(n)}_k\sspeq \sum_{j=0}^{k-1}\,R^{(n)}_{m+k-(2j+1)} \ ,\ 1\sspleq m\sspleq k\,  . 
\Eeq
This is finally the DPF found by {\sl Steinbach} in \cite{Steinbach} and gives a linearization of the DSR products. We need only to consider $m \sspgeq 2$ ($m\sspeq 1$ becomes trivial), and due to the symmetry of this formula under the transformation $k\to n-k$ (using the rules of eq. $(8)$) is suffices to consider $k\sspin \{1,2,\,...\, \floor{\frac{n}{2}}\}$. Hence one has only to consider the 
\dstyle{ {\binomial{\floor{\frac{n}{2}}}{2}} } products for $2\sspleq m\sspleq k\sspleq \floor{\frac{n}{2}}$.
\psn
The DPF formula looks un-symmetric with respect to $m\sspleftrightarrow k$, but it is, in fact, symmetric because \dstyle{ \sum_{j=0}^{m-1}\,R^{(n)}_{k+m-(2j+1)} } for $k\sspkl m$ reduces, due to cancellations after using the rule ${\bf i)}$ from eq $(8)$, to the expected sum with only $k$ terms. One can see this for even and odd $m\sspgr k$ separately, remembering that $R^{(n)}_0\sspeq 0$ in the former case.\psn 
The idea is  to work out, for a given $n\sspgeq 4$, all DPFs of interest, using the rules from eq. $(8)$, especially ${\bf i)}$ and ${\bf ii)}$ in order to write all  \dstyle{ {\binomial{\floor{\frac{n}{2}}}{2}} } products as linear combinations of the DSRs.\pn
{\bf Example 1}: $n\sspeq 7$\ {\bf heptagon}, treated in detail in reference \cite{Steinbach}. Here we recapitulate and link to OEIS \cite{Sloane} sequences (the analoga of {\sl Fibonacci} numbers in the pentagon case).  Superscripts and arguments $7$ are suppressed. $[2,2]$: $\rho^2\speq \sigma\sspp 1$.  $[2,3]$: $\rho\,\sigma\sspeq \sigma \sspp \rho \sspp 0$. $[3,3]$: $\sigma^2\sspeq \rho\sspp \sigma\sspp 1$. Here, $R_5\sspeq R_2$ was used. This shows that the $\mathbb Q$-vector space basis is at most $<1,\rho,\sigma>$. It will be shown to be indeed the heptagon basis in sect. $4$, where the power basis will be used instead. \pn
With these DPFs one can compute all powers of interest in the heptagon basis $<1,\rho,\sigma>$ (this has already been done explicitely for $\sigma$ in \cite{Steinbach}, p. 28). \psn
$\rho^k\sspeq$\, \seqnum{A052547}$(k-2) \,1 \sspp$\, \seqnum{A052547}$(k-1)\, \rho\sspp$\, \seqnum{A006053}$(k)\, \sigma$, $k\sspgeq 0$,  \pn
$\rho^{-k}\sspeq$\, \seqnum{A077998}$(k) \,1 \sspp$\, \seqnum{A077998}$(k-1)\, \rho\sspm$\, \seqnum{A006054}$(k+1)\, \sigma$, $k\sspgeq 0$ ,\pn
$\sigma^k\sspeq$\, \seqnum{A106803}$(k-1) \,1 \sspp$\, \seqnum{A006054}$(k-1)\, \rho\sspp$\, \seqnum{A106803}$(k)\, \sigma$, $k\sspgeq 0$,  \pn
$\sigma^{-k}\sspeq (\sigma\sspm \rho)^{k}\sspeq$\, \seqnum{A052547}$(k) \,1 \sspm$\, \seqnum{A006053}$(k+1)\, \rho\sspm$\, \seqnum{A052547}$(k-1)\, \sigma$, $k\sspgeq 0$, \pn
$(\rho\,\sigma)^k\sspeq  (\rho\sspp \sigma)^k\sspeq$\, \seqnum{A120757}$(k) \,1 \sspp$\, |\seqnum{A006054}$(k-1)|\, \rho\sspp\, 4$\, \seqnum{A181879}$(k)\, \sigma$, $k\sspgeq 0$,  \pn
$(\rho\,\sigma)^{-k}\sspeq  (\rho\sspp \sigma)^{-k}\sspeq$\ \seqnum{A085810}$(k) \,1 \sspp\, (-1)^k$\, \seqnum{A181880}$(k-2)\, \rho\sspp\, (-1)^{k+1}$\, \seqnum{A116423}$(k+1)\, \sigma$, $k\sspgeq 0$.\pn
One can also compute \dstyle{{\frac{1}{a\,1\sspp b\,\rho\sspp c\,\rho}} \sspeq A\,1\sspp B\, \rho\sspp \sigma} and find with  $N(a,b,c)\speq a^3 \sspm b^3\sspm c^3\sspm 2\,a\, b^2\sspm a\,b\,c \sspp a^2\,b \sspp b^2\,c\sspp 2\,a^2\,c\sspm a\,c^2\sspp 2\,b\,c^2$
\Beqarray
A&\sspeq& \frac{1}{N(a,b,c)}\, (a^2\sspm b^2\sspp a\,b\sspp 2\, a\,c \sspm b\,c)\ \ , \ \ B\sspeq \frac{1}{N(a,b,c)}\, (b^2\sspm c^2\sspp a\,b)\ \ , \\
C&\sspeq& \frac{1}{N(a,b,c)}\, (c^2\sspm b^2\sspp a\,c\sspm b\,c).
\Eeqarray
{\bf Example 2:}\ $n\sspeq 9$ ({\bf Enneagon}), where we use $\rho\sspeq  R^{(9)}_2,\ \sigma\sspeq  R^{(9)}_3, \tau\sspeq  R^{(9)}_4$.\ $[2,2]:\ \rho^2\sspeq 1\sspp \sigma$,\ $[2,3]:\ \rho\sigma\sspeq \rho\sspp \tau$, \ $[2,4]:\ \rho\tau\sspeq \sigma\sspp \tau$,\ $[3,3]:\ \sigma^2\sspeq 1\sspp \tau\sspp \sigma$,\  $[3,4]:\ \sigma\,\tau\sspeq \rho\sspp \sigma\sspp \tau$, \ $[4,4]:\ \tau^2\sspeq 1\sspp \rho\sspp \sigma\sspp \tau$. This DSR-algebra (over $\mathbb Q$) shows that not all ratios (including $R^{(n)}_1\sspeq 1$) are linear independent: $\tau\,(\rho\sspm 1)\sspeq \sigma\sspeq \rho^2\sspm 1 \sspeq (\rho\sspp 1)\, (\rho\sspm 1)$, \ie $\tau\sspeq \rho\sspp 1$. Therefore, $\rho^3\sspeq \rho\,(\sigma \sspp 1)\sspeq 2\, \rho\sspp \tau\sspp 2\sspeq 3\,\rho\sspp 1\, $.  In sect. $3$ we will see that the algebraic number $\rho(9)$ has degree $3$, and its minimal polynomial is indeed  $C(9,x)\sspeq x^3\sspm 3\,x\sspm 1$.  The enneagon basis is thus $<1,\rho,\sigma>$ which can be related to the power basis $<1,\rho,\rho^2>$ (remember that we use the same notation for different bases). The reduced DSR-algebra (the algebra modulo $C(9,\rho)\sspeq 0$ ) is $\rho^2\sspeq 1\sspp \sigma$, $\rho\,\sigma\sspeq 1\sspp 2\,\rho$, and $\sigma^2\sspeq 2\sspp \rho\sspp \sigma$. 
%% We refrain from giving formulae for the powers of $\rho$ and the corresponding sequences (see {\it Table 4}).
\psn
{\bf Remark 1:} $[2,2]$  in eq. $(13)$ becomes  \dstyle{ (R^{(n)}_2)^2\sspeq R^{(n)}_3\sspp 1}, showing that for all interesting values $n\sspgeq 4$ one has  $\sigma(n)\sspdef  R^{(n)}_{3} \sspeq \rho(n)^2 \sspm 1\ $. This is known already from eq $(1)$ for $k=3$, and it will  be used in sect. $3$ as \dstyle{R^{(n)}_{3}} rewritten in terms of the power basis of $\mathbb Q(\rho(n))$. 
\psn
{\bf Remark 2:} For $n=4$ and $n=5$ the second diagonal (defining \dstyle{R^{(n)}_{3}})  is not of interest because it is in the lower half plane. According to eq. $(8)$, {\bf ii)} \dstyle{R^{(5)}_{2} \sspeq R^{(5)}_{3}}, hence $\sigma(5)\sspeq \rho(5)$ and the general relation $[2,2]$ from eq. $(13)$ between $\sigma(n)$ and $\rho(n)$ leads to $\rho(5)\sspeq \rho(5)^2\sspm 1$, the golden section formula $\rho(5)\sspeq \varphi\sspdef (1\sspp \sqrt{5})/2$. In sect. $2$ it will be seen that the minimal polynomial for $\rho(5)$ is indeed $C(5,x)\sspeq x^2 \sspm x\sspm 1$.\psn 
Instead of deriving the DPF from the product formula for $S$-polynomials and the rules of eq. $(8)$ one can prove it directly by induction from the $S-$recurrence including negative indices in accordance with these rules.\psn
{\bf Proposition 1:}\psn
 The DPFs  $[m,k]^{(n)}$, eq. $(13)$\, but now for all  $ 2\sspleq m ,\ 2\sspleq k$, follow from the recurrence, eq. $(11)$,  and the rules {\bf i)} and {\bf ii)} from eq. $(8)$.
\psn
{\bf Proof:} This is shown by double induction. First one shows this for given $m \sspgeq 2$ for all $k \sspgeq 2$ by induction over $k$. Then by induction over $m$ for all $k\sspgeq 2$. \pn
The $k$-induction uses as starter $[2,2]^{(n)}$ which is the recurrence. Assume that   $[2,k^{\prime}]^{(n)}$ is true for all $k^{\prime}\sspeq 2,3,...,k-1 $. Now (we omit the superscripts) $R_2\,R_{k}\sspeq R_2\,(R_2\, R_{k-1}\sspm R_{k-2})$ from the recurrence. Then use for both terms the induction hypothesis. Note that one obtains for the first term $k-1$, and for the second one $k-2$ terms. However, due to rule ${\bf ii)}$ only two terms survive in each case. (The same type of cancellation was at work when we remarked above on the symmetry of the DPF formula in $m$ and $k$.) Therefore, one obtains $R_2\, (R_k\sspm R_{k-2}\sspp 0)\sspm (R_{k-1}\sspm R_{k-3}\sspp 0)$ which becomes, after use of the recurrence applied twice,  $R_{k+1}\sspp R_{k-1}$. This is indeed the desired result for $[2,k]$ if one uses again rule {\bf ii)} to get a truncation of  the sum after two terms.\pn
The $m$-induction uses as starter $[2,k]$ for all $k\sspgeq 2$, which has  just been established.
Then assume that $[m^{\prime},k]$ is true for all $m^{\prime}\sspeq 2,3,...,m-1$.
$R_m\,R_k\sspeq R_2\,R_{m-1}\,R_k\sspm R_{m-2}\,R_k$ from the recurrence. Assume the induction hypothesis for each term, obtaining, after use of the recurrence for the first sum, \dstyle{ \sum_{j=0}^{k-1}\,(R_{m+k-(2j+1)} \sspp R_{m+k-2-(2j+1)})\sspm  \sum_{j=0}^{k-1}\,(R_{m-2+k-(2j+1)}} which is indeed the assertion after cancellation of the last two sums.\hskip 14cm $\square$\psn
As mentioned in the {\it  example 2} the DSR-algebra turns sometimes out to be reducible because some DSRs can be expressed as rational linear combinations of other ones. Later it will become clear that this happens precisely whenever $\floor{\frac{n}{2}}\sspm \delta(n)\sspgr 0$, and this is the number of linear dependent  DSRs. See {\it Table 1} for details for $n\sspeq 3,...,12$.\psn
We shall see in {\it sects. 3} and {\it 4} that it is simpler to use the power basis of $\mathbb Q(\rho(n))$ and the minimal $C(n,x)$-polynomials of the algebraic number $\rho(n)$ instead of the DSR-algebra. Then one obtains automatically the reduced algebra.
%%%%%%%%%%%%%%%%%%%%%%%%%%%%%%%%%%%%%%%% 
\section{Minimal polynomial of $\bf {\boldsymbol\rho}(n)$, its zeros, absolute term and factorization} 
%%%%%%%%%%%%%%%%%%%%%%%%%%%%%%%%%%%%%%%% 
\hskip 1cm The minimal polynomial of an algebraic number $\alpha$ of degree $d_\alpha$ is the monic, minimal degree  rational polynomial which has as root, or as one of its roots, $\alpha$. This degree $d_\alpha$ is $1$ iff $\alpha$ is rational, and the minimal polynomial in this case is $p(x)=x\sspm \alpha$. For the notion `minimal polynomial of an algebraic number' see, {\it e.g.}, \cite{Niven}, p. 28 or \cite{Ribenboim} p. 13.
\psn
For the algebraic number  \dstyle{\rho(n)\sspdef 2\,\cos\left ({\frac{\pi}{n}}\right)}, for $n\sspin \mathbb N$, the degree is $\delta(1)\sspeq 1$, and \dstyle{\delta(n)\sspeq {\frac{\varphi(2\,n)}{ 2}} } for $n\sspgeq 2$ , with {\sl Euler}'s totient function $\varphi(n)\sspeq \seqnum{A000010}(n)$. This is the sequence \seqnum{A055034}$(n)$. The sequence of allowed $\delta$ values is given in \seqnum{A207333}. The array with the indices of the polynomials for given allowed $\delta$ values is shown in \seqnum{A207334}. Of course, $\delta$ is not multiplicative, \eg $2\sspeq \delta(6)\sspneq \delta(3)\,\delta(3)\sspeq 1\cdot 1\sspeq 1 $. These polynomials can be obtained from the ones of  \dstyle{\cos\left ({\frac{2\,\pi}{n}}\right)} which are found, {\it e.g.}, under  \seqnum{A181875}/ \seqnum{A181876}, and they have been called there  $\Psi(n,x)$. See also \cite{Lehmer}, and \cite{Niven}, Theorem $3.9$, p. 37, for the degree $d(n)$ of these polynomials.  From the trivial identity \dstyle{\cos\left ({\frac{\pi}{n}}\right) \sspeq \cos\left ({\frac{2\,\pi}{2\,n}}\right) } one finds the minimal polynomial of  \dstyle{2\,\cos\left ({\frac{\pi}{n}}\right)}, called here $C(n,x)$, from 
\Beq
C(n,x)\sspeq 2^{\delta(n)}\, \Psi\left(2\,n,{\frac{x}{2}}\right)\ . 
\Eeq
Therefore the above formula for $\delta(n)$ derives from $d(2\,n)$. These polynomials  are given for $n=1,\, ...,\,30$ in {\it Table 2}, and the first $15$ rows of their coefficient array are shown in {\it Table 3}. Concerning the parity of the degree $\delta$ one has the following lemma.\psn
{\bf Lemma 1:\  Parity of the degree} $\boldsymbol{\delta}$\psn
$\delta(n)$ is odd iff $n\sspeq 1,\, 2,$ and $p^e$, with an odd prime p of the form $4\,k\sspp 3$, $k\sspin \mathbb N_0$, and $e\sspgeq 1$.\psn
{\bf Proof:} The case $n=1$ is clear by definition. For $n\sspeq 2$ one uses the prime number factorization $n\sspeq 2^\beta\, \prod_{j=1}^N\, p_j^{e_j}$, with $\beta\sspgeq 0$ and $e_j\sspeq 1$. It follows from the definition that  \dstyle{\delta(n)\sspeq 2^{\beta-1}\, \prod_{j=1}^N\, (p_j-1)\, p_j^{e_j-1}}. For $\beta\sspgeq 2$ this is always even. For $\beta\sspeq 0$ one needs $N\sspeq 1$, otherwise this will be even. For $N\sspeq 1$ one needs  \dstyle{\frac{p-1}{2}} to be odd, \ie $p\sspeq 4\,k\sspp3$ with $k\sspin \mathbb N_0$. If $\beta\sspeq 1$ one needs $N\sspeq 0$ in order to have an odd value. Thus $n\sspeq 2$. \hskip 13cm $\square$ \psn
In order to relate to sect. 2 we are interested here in $n\sspgeq 4$. It turns out that these polynomials are in fact integer (not only rational) polynomials. The proof  will use the following lemma based on our divisor product representation paper \cite{Lang2}.\psn
{\bf Lemma 2:\  Minimal polynomial} $\bf C(n,x)$ {\bf written as a rational function.}
\Beq
C(n,x)\sspeq {\frac{p(d_p(n),x)}{q(d_q(n),x)}}\ ,
\Eeq
with monic {\it integer} polynomials $p$ and $q$ with a certain degree $d_p(n)$ and $d_q(n)\sspeq \delta_p(n)\sspm \delta(n)$, respectively.\pbn
{\bf Proof}: $\Psi\left(2\,n,{\frac{x}{2}}\right)$ is obtained from the unique divisor product representation $dpr(2\,n)$ defined in \cite{Lang2} by replacing each $a(k)$ in the numerator, as well as in the denominator,  by $t(k,{\frac{x}{2}})$ which is given as a difference of monic {\it integer} polynomials $\hat t$ which have been defined already in eq. $(2)$, multiplied with a certain prefactor. 
\Beq
t(k,{\frac{x}{2}})\sspeq{\Caseszwei{{\frac{1}{2^{{\frac{k}{2}}+1}}}\,\left( \hat t(\frac{k}{2}\sspp 1,x)\sspm \hat t(\frac{k}{2}\sspm 1,x)\right)}{\text{if}\ $k$ \text{is even\, },}{{\frac{1}{2^{{\frac{k-1}{2}}+1}}}\,\left( \hat t(\frac{k+1}{2},x)\sspm \hat t(\frac{k-1}{2},x)\right) }{\text{if}\ $k$ \text{is odd\,}. }}
\Eeq  
The monic $\Psi(2\,n,x)$ polynomials have degree $d(2\,n)$ (see \eg \cite{Lang1}) which implies that  all the prefactors in this replacement of the $a(k)$s in $dpr(2\,n)$ have to become $1/2^{d(2\,n)}$. (This could be formulated as a separate lemma). \Eg $n\sspeq 34$ with $dpr(34)\sspeq {\frac{a(34)\,a(1)}{a(17)\,a(2)}}$ has from the numerator the factor $1/2^{18+1}$ and from the denominator $2^{9+2}$, fitting with $1/2^{d(34)}$ because $d(34)\sspeq 8$. Therefore, one may in the calculation of $C(n,x)$, found above from $\Psi(2\,n,{\frac{x}{2}})$, forget about these prefactors in the replacements altogether. Thus the numerator, {\it resp.} denominator, is a product of  monic integer polynomials $\hat t$ leading to the monic integer polynomials $p$, {\it resp.} $q$, of a certain degree $d_p(n)$, {\it resp.} $d_q(n)$. (One could give more details on these degrees but this is not important here. Trivially, $d_p(n)\sspm d_q(n)\sspeq d(2\,n)$. For the given example $n=34$: $(18+1)\sspm (9+2)\sspeq 8$ from the degrees of the $\hat t$ polynomials.)  Because C(n,x) is a minimal {\it polynomial} this rational function allows polynomial division without remainder.  \hskip 15cm $\square$
\psn
{\bf Proposition 2:}  $\bf C\in {\mathbb Z}[x]$\psn
 $C(n,x)$, the minimal polynomial of \dstyle{\rho(n)\sspeq 2\,\cos\left ({\frac{\pi}{n}}\right)}, is an {\it integer} monic polynomial.
\psn
{\bf Proof:} From  {\it lemma 2} we have $q(d_q(n),x)\, C(n,x)\sspeq p(d_p(n),x)$ which leads by induction to the result that the integer coefficients of the monic polynomials $p$ and $q$ imply integer coefficients for the monic polynomial $C$. Call these monic  polynomials $q(N,x)\sspdef \sum_{l=0}^N\, q_l\,x^l$, $C(M,x)\sspdef \sum_{k=0}^M\, c_k\,x^k$ and   $p(N+M,x)\sspdef \sum_{l=0}^{N+M}\, p_l\,x^l$ with $q_N=1$, $c_M=1$ and  $p_{N+M}=1$. Collecting terms for $x^{N+M-j}$, for  $j=(0),1,\,2,\, ...,\,M$, one obtains the formula for the $C$-coefficient $c_{M-j}$ in terms of lower indexed  ones: \dstyle{c_{M-j}\sspeq p_{N+M-j} \sspm \sum_{k=1}^{j}\,c_{M-j-k}\,q_{N-k}}. With this the inductive proof on $j$, using the integer coefficients of $q, p$ and the integer higher $C$-coefficients due to the induction hypothesis becomes obvious. The starting point is the trivial $j=0$ case. \hskip 15cm $\square$
\psn
Next we give all the zeros of $C(n,x)$. This follows directly from the knowledge of all zeros of $\Psi(2\,n,{\frac{x}{2}})$ (see, \eg \cite{Lang1}, from p. 473 of \cite{WatkinsZeitlin}) \psn
{\bf Proposition 3:} {\bf All zeros of $\bf C$} \psn
 \Beq
 C(n,x)\sspeq \prod_{{k=1}\atop {gcd(k,2\,n)=1}}^{n-1}\, \left(x\sspm 2\,\cos\left( {\frac{\pi\, k}{n}}\right)\right) \ ,\  n\in \mathbb N\, .
\Eeq   
In the product the index $k$  starts with $1$ because $gcd(0,2\,n)\sspeq 2n$, and it stops at $n-1$ for all $n\sspgeq 2$ because $gcd(n,2\,n)=n$. One has to omit the even $k$ values $\sspgr 0$ which leads to (for $n\sspeq 1$ one takes only $l \sspeq 0$ in the following product)
\Beq
 C(n,x)\sspeq \prod_{{l=0}\atop {gcd(2\,l+1,n)=1}}^{\floor{{\frac{n-2}{2}}}}\, \left(x\sspm 2\,\cos\left( {\frac{\pi\, (2\,l+1)}{n}}\right)\right) \ .
\Eeq  
{\bf Proof:} This is a simple consequence of the mentioned known results for $\Psi(2\,n,{\frac{x}{2}})$. It implies, {\it en passant}, a formula for the known degree $\delta(n)$ (see \seqnum{A055034}) in terms of the number of factors of this product. \hskip 6 cm $\square$ \psn
The following proposition lists the nonnegative and negative zeros of $C$ for prime $n$. A vanishing zero (which will be counted later as positive) appears only for $n=2$.
\psn
{\bf Proposition 4:} {\bf Non-negative and negative zeros of  $\bf C(p,x)$}\psn
{\bf i)} If $n=2$ then $C(2,0)\sspeq 0$ , and this is the only case with a vanishing zero.\psn
{\bf ii)} If $n$ is an odd prime $1\, (mod\, 4)$ then the number of positive and negative zeros coincides, and this number is $\frac{n-1}{4}$.  These zeros are  $ (-1)^{k+1}\,2\, \cos\left( {\frac{\pi\, k}{n}}\right)$, with $k$ values $1,2,...,{\frac{n-1}{2}}\sspeq \delta(n)$. \pn
{\bf iii)} If $n$ is an odd prime $3\, (mod\, 4)$ then the number of positive zeros exceeds the number of negative ones by one, and the number of negative ones is $\frac{n-3}{4}$. These zeros are $ (-1)^{k+1} \,2\, \cos\left( {\frac{\pi\, k}{n}}\right)$ with $k$ values $1,2,...,{\frac{n-3}{2}}$, and an extra positive zero appears for $k={\frac{n-1}{2}}\sspeq \delta(n)$ .\psn
{\bf Proof:} \psn
{\bf i)} Vanishing zeros require $n\sspequiv 2\, (mod\,4)$ from $\frac{k}{n}\sspeq\frac{1}{2} $  and $k$ odd in eq. $(20)$. A vanishing zero appears for $n=2$. For other such $n$ values, $n\sspeq 4\,K\sspp 2$, $K\sspgeq 1$, the odd  $k\sspeq \frac{n}{2}\sspgr 1$ divides $n$, hence it does not appear in the  product of eq. $(20)$.\psn  
For odd primes $n\sspeq p\sspeq 2\,q\sspp 1$ the product in eq. $(20)$ is unrestricted, and there are \dstyle{ (q-1) \sspp 1\sspeq {\frac{p-1}{2}}}  factors, in accordance with the degree \dstyle{\delta(p)\sspeq \frac{\varphi(2\,p)}{2}}.  The $cos$ function will produce negative zeros for \dstyle{{\frac{2\,l\sspp 1}{p}}\sspgr {\frac{1}{2}} }, \ie for $2\,l\sspgeq q$.\psn
{\bf ii)} $n\sspeq p\sspeq 4\,K\sspp 1$ means that $l\sspeq 0,\, 1,\, ...,\, K-1$ lead to positive zeros, and $l\sspeq K,\,K+1,\,...,\,2K-1$ to negative ones. In both cases there are \dstyle{K\sspeq {\frac{n-1}{4}}} such zeros. For the $l$ values leading to  negative zeros one can use  the formula \dstyle{\cos\left({\frac{\pi}{p}}\,(2\,l+1)\right)\sspeq -\cos\left({\frac{\pi}{p}}\,(p\sspm (2\,l+1))\right)}. In this way even multiples of \dstyle{{\frac{\pi}{p}}} appear, and one produces, counted backwards, beginning with the largest $l$ value, the new values $2,4,..., 2\,K$. This  proves {\bf ii)}. As a test consider $n\sspeq 13$, $K\sspeq 3$: the $2\, l\sspp 1$ values for the positive zeros are $1,\,3,\,5$, and the ones for the negative ones are $7,\,9,\,11$. The latter become (we use underlining to indicate that they come with a minus sign in the final formula) $\underline{6},\,\underline{4},\,\underline{2}$. Rearranged these values are $1,\,\underline{2},\,3,\,\underline{4},\,5,\,\underline{6}$, as given in {\bf ii)}.\psn
{\bf iii)} For $n\sspeq p\sspeq 4\,K\sspp 3$ a similar analysis leads to $l\sspeq 0,\,1,\,...,\,K$, and $l\sspeq K+1,\, ...,\, 2\,K$ for positive, and negative zeros, respectively. This shows that the number of negative zeros is \dstyle{K\sspeq {\frac{p\sspm 3}{4}}}, and the number of positive ones  exceeds the one for negative zeros by one. Again, the odd values  leading to negative zeros are transformed to even ones $2,\, 4,\,...,\,2\,K$ with a minus sign in front of $2\, \cos$. Take, \eg $n=19$, $K\sspeq 4$, with the new $k$ values of {\bf iii)}\, $1,\, $\underline{2},\, 3,\, $\underline{4},\, 5,\,$\underline{6},\,7,\,$\underline{8},\,9 $.   \hskip 7 cm $\square$\psn
For general values $n$ one can also find the number of positive and negative zeros of $C(n,x)$. In the non-prime or non-power-of-$2$ case the $gcd$ restriction in the product excludes certain $l$ values. Therefore one has to eliminate from the unrestricted $2\,l\sspp 1$ values in the product  of eq. $(20)$ all odd multiples for each odd prime dividing $n$. Multiples of $2$ are irrelevant for even $n$ , therefore only the odd primes dividing n are of interest. This shows that the prime factors of the {\it squarefree kernel} of $n$, denoted here by $sqfk(n)$, \ie the largest squarefree number dividing $n$ (see \seqnum{A007947}) will be of interest. This squarefree kernel is also known as radical of $n$, denoted by $rad(n)$. We will denote the set of primes of this kernel by $sqfkset(n)$, \Eg $sqfkset(2^3\cdot 3^2\cdot11)\sspeq \{2,3,11\}$, $sqfk(2^3\cdot 3^2\cdot11)\sspeq 66$ (strip off all exponents in the prime number factorization of $n$). Because one may encounter multiple counting (\eg $15$ is hit by the 
multiples  of $3$ as well as $5$ for any $n$ which has in its squarefree kernel set the primes $3$ and $5$, and which is larger than $16$) one can employ {\it PIE}, the principle of inclusion-exclusion, \eg \cite{Charalambides}, p. 134,  to get a correct counting.\psn
{\bf Proposition 5:} {\bf Number of positive and negative zeros of} $\bf C(n,x)$\psn
{\bf i)} If $n=2^m,\, m\sspgeq 1$, one has for $m=1$ a vanishing zero (see {\it proposition 4 i)}). For $m\sspgeq 2$ the number of positive zeros, which is  \dstyle{{\frac{n}{4}}},  coincides with the one for negative ones. This is in accordance with the degree \dstyle{\delta(2^m)\sspeq 2^{m-1}\sspeq {\frac{n}{2}}}. These zeros lie symmetric: \dstyle{\pm\,2\,\cos\left({\frac{\pi}{n}}\,(2\,l+1)\right), l\sspeq 0, 1,...,{\frac{n}{4}}-1}. \psn {\bf ii)} If $n=1$ then there is only the negative zero $-2$. If $n\sspgeq 3$ is not a power of $2$ then the number of positive zeros, called $\delta_+(n)$, is
\Beq
\delta_+(n)\speq K(n) \sspp \sum_{r=1}^{M(n)}\, (-1)^r\,\sum_{<i_1,i_2,...,i_r>}\,\floor{{\frac{1}{2}}\,\left( {\frac{L(n)}{p_{i_1}\cdots p_{i_r}}}\sspm 1 \right) } \ ,
\Eeq
where $M(n)$ (sometimes called $\omega(n)$) is in the odd $n$ case the number of elements (cardinality) of the set $sqfkset(n)$, denoted by $|sqfkset(n)|$, (see \seqnum{A001221}). $M(n)$ is for even $n$ the number of odd primes in $sqfkset(n)$, \ie $|sqfkset(n)|\sspm 1$. This is because multiples of $2$ are irrelevant here. The sum $<i_1,i_2,...,i_r>$ extends over the \dstyle{{\binomial{M}{r}}} combinations $1\sspleq i_1\sspkl i_2\sspkl\, ...\, \sspkl i_r\sspleq M$. Only the  odd primes from the set $sqfkset(n)$  enter. The values $K(n)$ and $L(n)$ depend on the parity of $n$ and they are given by  
\Beqarray
&&\alpha)\ \ \,  n\,  {\text{even}}:\hskip 1.8cm \ K(n) \sspeq \floor{{\frac{n-2}{4}}}\, ,\ \ L(n)\sspeq 2\,K(n)\sspp 1 \ , \\
&&\beta 1)\  n\,\ \ {\text{odd}},\ 1 (mod\, 4):\ K(n) \sspeq {\frac{n-5}{4}} \ \ ,\ \quad L(n)\sspeq 2\,K(n)\sspp 1 \sspeq {\frac{n-3}{2}} , \\
&&\beta 2)\  n\,\ {\text{ odd}},\ 3 (mod\, 4):\ K(n)\sspeq {\frac{n-3}{4}}\ \ ,\ \quad L(n)\sspeq 2\,K(n)\sspp 1 \sspeq {\frac{n-1}{2}}.
\Eeqarray    
Note that also negative values for the floor function may appear. In this way the pure prime case from {\it proposition $4$} is also included. \psn
 {\bf iii)} If $n\sspgeq 3$ is not a power of $2$ then the number of negative zeros, called $\delta_-(n)$, is
\Beq
\delta_-(n)\speq N(n) \sspp \sum_{r=1}^{M(n)}\,(-1)^r \,\sum_{<i_1,i_2,...,i_r>}\,\left\{ \floor{{\frac{1}{2}}\,\left( {\frac{P(n)}{p_{i_1}\cdots p_{i_r}}}\sspm 1 \right)} \sspm  \floor{{\frac{1}{2}}\,\left( {\frac{L(n)}{p_{i_1}\cdots p_{i_r}}}\sspm 1 \right)}    \right\}\ ,
\Eeq
where $M(n)$ and $L(n)$ are as above, and 
\Beqarray
&&\alpha)\ \ \,  n\ \  {\text{even}}:\hskip 1.8cm \ N(n) \sspeq {\frac{n-2}{2}}  \sspm \floor{{\frac{n-2}{4}}}\, ,\ \ P(n)\sspeq n\sspm 1 \ , \\
&&\beta 1)\  n\,\ \ {\text{odd}},\ 1 (mod\, 4):\ N(n) \sspeq {\frac{n-1}{4}} \ \ ,\ \quad P(n)\sspeq n\sspm 2, \\
&&\beta 2)\  n\,\ {\text{ odd}},\ 3 (mod\, 4):\ N(n)\sspeq K(n)\sspeq {\frac{n-3}{4}}\ \ ,\ \quad P(n)\sspeq n\sspm 2.
\Eeqarray    
{\bf Proof:}\psn
{\bf i)} All \dstyle{l\sspeq 0,\,..,\,\floor{\frac{n-2}{2}}\sspeq 2^{m-1}\sspm 1} values contribute to the product of eq. $(20)$, compatible with the degree \dstyle{\delta(2^m)\sspeq  {\frac{\varphi(2^m)}{2}} \sspeq 2^{m-1}\sspeq {\frac{n}{2}}}. The negative zeros appear for $2\,l+1\sspeq 2\,2^{m-2}+1,...,n-1$, hence there are  ${\frac{n}{4}}$ of them. With the $\cos$ formula given in the proof of {\it proposition 4} {\bf ii)} with $p\sspto n$, they can be rewritten in the notation, also used in the above context, when read backwards, as $\underline{1},\underline{3},...,\underline{2^{m-1}-1}$. Therefore one obtains the same $\cos$ arguments like for the ${\frac{n}{4}}$  positive zeros, only the overall sign is different.\psn
{\bf ii)} For $n\sspgeq 3$, not a power of $2$, one uses PIE to count the positive zeros. For this consider the odd multiples of some odd number $a$, say $(2\,k\sspp 1)\,a$,  (only such multiples come up as $2\,l\sspp 1$ values in the product $(20)$) up  to  $k_{max}$ such that all positive zeros are reached.  If the largest $2\,l+1$ value in the product which leads to  a positive zero is  $\bar n$,  then \dstyle{k_{max} \sspeq \floor{{\frac{1}{2}}\, \left({\frac{\bar n}{a}}\sspm 1\right)}}, and there are $k_{max} \sspp 1$ of these odd multiples of $a$.  In the following application of PIE $a$ will be taken as any odd prime or product of odd primes from the set $sqfkset(n)$.\psn
$\alpha$) Case even $n$, not a power of $2$: The counting of the  positive zeros starts at level $r\sspeq 0$ with all positive ones in the unrestricted product eq. $(20)$. There are \dstyle{\floor{\frac{n-2}{4}}\sspp 1} of them (the $l$ values are \dstyle{0,\,1,\,...,\, \floor{{\frac{n-2}{4}}})}. In the next step, $r\sspeq 1$, all odd multiples of every odd prime $p_{i_1}$ in $sqfkset(n)$ not exceeding \dstyle{2\,\floor{{\frac{n-2}{4}}} \sspp 1} are discarded (there are $M(n)$ such primes, where $M(n)\sspdef |sqfkset(n)|-1$). This leads to a subtraction of the form \dstyle{-\sum_{i_1=1}^{M(n)}\,\floor{{\frac{1}{2}}\,\left({\frac{L(n)}{p_{i_1}}}\sspm 1 \right)}\sspp 1 } with $L(n)$ given in the {\it proposition} {\bf ii)} $\alpha)$. Now double subtractions may have appeared and one adds, at step $r\sspeq 2$, all odd multiples of the product of two odd primes from $sqfkset(n)$, call them $p_{i_1}\,p_{i_2}$ with $i_1\sspkl i_2$ (because no problem of over-subtraction appeared for $i_1\sspeq i_2$). This leads to the term  \dstyle{+\sum_{1\sspleq i_1<i_2\sspleq M(n)}\,\floor{{\frac{1}{2}}\,\left({\frac{L(n)}{p_{i_1}\,p_{i_2}}}\sspm 1 \right)}\sspp 1}. Now one has to correct for triple products in the step $r\sspeq 3$ , etc., up to level $r\sspeq M(n)$. This leads to the formula given in part {\bf ii)} $\alpha$) if one observes that the $+1$ term in each $r$-sum becomes a \dstyle{\binomial{M(n)}{r}} term because of the number of terms of this sum.  Then, because the alternating sum over row $M(n)$ in {\sl Pascal}'s triangle (see \seqnum{A007318}) vanishes ($(1-1)^{M(n)}\sspeq 0$), all these terms add up to $-1$ (the missing $r=0$ term). Therefore the formula starts with  \dstyle{\floor{\frac{n-2}{4}}}, the K(n) value given in {\bf ii)} $\alpha$), and only the floor term remains for each $r$-sum .
\psn
$\beta 1$) The counting of the positive zeros in the odd $n$ case distinguishes the two cases $1\, (mod\, 4)$ and $3\,(mod\,4)$. Here we discuss the former one. The procedure is the same as the one given for case $\alpha)$. One has just to determine the boundary values for the $l_{max}$ value leading still to a positive zero. This is $l_{max}\sspeq K-1$, if $n\sspeq 4\,K+1$, \ie \dstyle{{\frac{n-5}{4}}}. Hence there are $l_{max}\sspp 1$ such $l$ values to start with at level $r\sspeq 0$. Again the extra $+1$ can later be taken as $r\sspeq 0$ term \dstyle{\binomial{M(n)}{0}} for the vanishing alternating sum over row No. $M(n)$  in {\sl Pascal}'s triangle. Here $M(n)\sspdef |sqfkset(n)|$. Therefore the PIE formula  starts with the $K(n)$ given in $\beta 1)$ of the {\it proposition}. The $L(n)$ of the PIE sum is now \dstyle{2\,{\frac{n-5}{4}}\sspp 1\sspeq {\frac{n-3}{2}}}, the largest $2\,l\sspp 1$ value leading to a positive zero. 
\psn
$\beta 2$) In the  case $M(n)\sspdef |sqfkset(n)|$,  $n\sspeq 4\,K\sspp 3$ the maximal $2\,l\sspp 1$ value leading to a positive zero, the $L(n)$ in the formula,  is \dstyle{2\,K\sspp 1\sspeq \frac{n-1}{2}}, coming from the largest $l$ value which is $K$. This $K$ is the $K(n)$ in the formula claimed for this case $\beta 2$).\psn
{\bf iii)} For the number of negative zeros the counting is done by finding all the odd multiples of odd primes, or products of them, which cover the $2\,l\sspp 1$ range for the negative zeros. This is the difference of the number of such multiples for the whole range and the range for the positive zeros. \psn
$\alpha)$ If $n$ is even, not a power of $2$, the whole range is determined by $l_{max}\sspp 1\sspeq n\sspm1$. From the above {\bf ii)} $\alpha)$ case one knows the corresponding upper bound for the positive zeros, therefore the PIE formula has \dstyle{N(n) \sspeq ({\frac{n-2}{2}}\sspp 1)\sspm (\floor{{\frac{n-2}{4}}}\sspp 1)}, the number of factors in the unrestricted product which lead to negative zeros, which is the value claimed in the {\it proposition} {\bf iii)} $\alpha)$. Accordingly, $P(n)\sspeq n-1$ and $L(n)\sspeq 2\, \floor{{\frac{n-2}{4}}}\sspp 1$ as given in eq. $(26)$. \psn
$\beta 1)$ In this $n \sspequiv 1\,(mod\,4)$ case the maximal number $2\,l\sspp1$ in the product  is $n-2$, determining $P(n)$.  The corresponding number $L(n)$ for the positive zeros is taken from above as \dstyle{\frac{n-3}{2}}. There are \dstyle{N(n)\sspeq \left({\frac{n-3}{2}}\sspp1\right)\sspm \left({\frac{n-5}{4}} \sspp 1\right)\sspeq {\frac{n-1}{4}}} unrestricted factors in the product with negative zeros.\psn
$\beta 2)$ In this $n \sspequiv 3\,(mod\,4)$ case one has also $P(n)\sspeq n\sspm2$ like  for any odd $n$. $L(n)$ is taken from the corresponding {\bf ii)} case as  \dstyle{{\frac{n-1}{2}}}
, and $N(n)$ is obtained from \dstyle{ \left({\frac{n-3}{2}} \sspp 1\right)-\left({\frac{n-3}{4}} \sspp 1\right)\sspeq {\frac{n-3}{4}}}. \hskip 6cm $\square$ 
\psn
We give the first entries of the $\delta_+(n)$ and $\delta_-(n)$ sequences. Remember that the vanishing zero for $n=2$ is here counted as positive. The A-numbers given for $\delta_-(2\,k)$ and $\delta_-(2\,k)$ are conjectured.\psn
$n$ even:\psn
$\delta_+(2\,k),\  k\sspeq 1,2,...$;\  \seqnum{A055034}\ : $[1, 1, 1, 2, 2, 2, 3, 4, 3, 4, 5, 4, 6, 6, 4, 8, 8, 6, 9, 8, 6, 10,...]$,\psn
$\delta_-(2\,k),\  k\sspeq 1,2,...$;\  \seqnum{A055034} with first entry $0$\ :  $[0,1,1,2,2,2,3,4,3,4,5,4,6,6,4,8,8,6,9,8,6,10,...]$,\psn
$n$ odd, $1 (mod\, 4)$: \psn
$\delta_+(4\,k\sspp 1),\  k\sspeq 0,1,...$;\ \ $[1,5,9,13,17,21,25,29,33,37,41,45,49,53,57,61,65,69,73,...]$,\psn
$\delta_-(4\,k\sspp 1),\  k\sspeq 0,1,...$;\ \ $[1,1,2,3,4,4,5,7,6,9,10,6,11,13,10,15,12,12,18,16,14,...]$,\psn
$n$ odd, $3 (mod\, 4)$: \psn 
$\delta_+(4\,k\sspp 3),\  k\sspeq 0,1,...$;\ \  $[1,2,3,2,5,6,5,8,6,6,11,12,8,10,15,10,17,18,10,20,...]$,\psn
$\delta_-(4\,k\sspp 3),\  k\sspeq 0,1,...$;\ \  $[0,1,2,2,4,5,4,7,6,6,10,11,8,10,14,8,16,17,10,19,20,...]$.\psn
The start of the sequences for odd $n$ is therefore\psn
$\delta_+(2\,k\sspp 1),\  k\sspeq 0,1,...$;\ \ $[0, 1, 1, 2, 1, 3, 3, 2, 4, 5, 2, 6, 5, 5, 7, 8, 4, 6, 9, 6, 10, 11,...]$,\psn 
$\delta_-(2\,k\sspp 1),\  k\sspeq 0,1,...$;\ \ $[1, 0, 1, 1, 2, 2, 3, 2, 4, 4, 4, 5, 5, 4, 7, 7, 6, 6, 9, 6, 10, 10, 6, 11,...]$.\psn
Finally, combining for all $n$:\psn
 $\delta_+(n),\  n\sspeq 1,2,...$;\ \seqnum{A193376}\ : $[0, 1, 1, 1, 1, 1, 2, 2, 1, 2, 3, 2, 3, 3, 2, 4, 4, 3, 5, 4, 2, 5, 6, 4, 5,...]$,\psn 
$\delta_-(n),\  n\sspeq 1,2,...$;\ \seqnum{A193377}\ : $[1,0,0,1,1,1,1,2,2,2,2,2,3,3,2,4,4,3,4,4,4,5,5,4,5,6,...]$.\psn
On can check that  $\delta(n) \sspm (\delta_+(n)\sspp \delta_-(n))$ vanishes in each case.
\pbn
Observe that from these examples it seems that $\delta_+(2\,k)\sspeq \delta_-(2\,k)\sspeq$\seqnum{A055034}$(k)\sspeq \delta(k)$ for $k\sspgeq 2$. This is compatible with $\delta(2\,k)\sspeq 2\,\delta(k)$, \ie \dstyle{{\frac{\varphi(4\,k)}{2}}\sspeq 2\, \varphi(2\,k)} for  these $k$ values. The latter eq.,  which holds  also for $k\sspeq 1$, can be checked by considering the two cases $k\sspeq 2^{e(1)}\, \cdot \text{(odd number)}$ and $k$ odd. Recall that $\delta(2)\sspeq \delta(1)\sspeq 1$, because $\delta(1)$ was special (it is not \dstyle{{\frac{\varphi(2)}{2}}}\,).\pbn 
{\bf Zeros of $\bf C(n,x)$ in the power basis $\bf <{\boldsymbol\rho}(n)^0=1,{\boldsymbol\rho}(n)^1, ..., {\boldsymbol\rho}(n)^{{\boldsymbol\delta}(n)-1}>$}\psn
As will be explained in the next section, the minimal polynomial $C(n,x)$ with degree $\delta(n)$ leads to the splitting field $\mathbb Q(\rho(n))$ which is a simple field extension of the rational field  $\mathbb Q$, and the degree of  $\mathbb Q(\rho(n))$ over $\mathbb Q$, the dimension of $\mathbb Q(\rho(n))$, considered as a vector space over $\mathbb Q$, is just the degree of $C$. In standard notation $[\mathbb Q(\rho(n)):\mathbb Q]\sspeq \delta(n).$ Therefore it is interesting to write the zeros of $C$ in the power basis for this $\mathbb Q$-vector space. This task is accomplished by using for the zeros 
\dstyle{\tilde{\xi}^{(n)}_l\sspeq  \cos\left({\frac{\pi}{n}}\,(2\,l\sspp 1)\right)},\ \dstyle{l\sspin \{0,1,..., \floor{\frac{n-2}{2}}\}} with $gcd(2\,l\sspp 1,n)\sspeq 1$ (see eq. $(20)$),the formula \dstyle{{\tilde \xi}^{(n)}_l\sspeq \hat t(2\,l\sspp 1,\rho(n))\sspeq 2\, T(2\,l\sspp 1,{\frac{\rho(n)}{2}})}\  (compare with eq. $(2)$). Remember that the integer coefficient array for these $\hat t$-polynomials is  shown in \seqnum{A127672}.  Of course, one has to reduce this integer polynomials using $C(n,\rho(n))\speq 0$, \ie one has to replace $\rho(n)^{\delta(n)}$ (and higher powers) by an integer polynomial of lower degree. This can be done on a computer and {\it Table 4} has been found with the help of Maple13 \cite{Maple}. For example, the zeros of $C(15,x)$ are $\tilde \xi^{(15)}_1\sspeq \rho(15),\  \tilde\xi^{(15)}_2\sspeq -2 \sspp 3\,\rho(15)\sspp \rho(15)^2 \sspm \rho(15)^3 ,\ \tilde\xi^{(15)}_3\sspeq -1\sspm 4\,\rho(15)\sspp \rho(15)^3,\  \tilde\xi^{(15)}_4\sspeq 2\sspm \rho(15)^2\ $. Here the reduction has been performed with $\rho(15)^4\sspto -\rho(15)^3 \sspp 4\,\rho(15)^2\sspp 4\, \rho(15)\sspm 1$. One can check these zeros which have approximate values  $1.956295201$, $0.209056928,\  -1.338261216$ and $-1.827090913$, respectively. We will call the zeros with increasing value $\xi_1,\xi_2,...,\xi_{\delta(n)}$.
\pbn
Our next objective is to compute $C(n,0)$ (the sign and the number multiplying $x^0$, also called the absolute term) and relate it first to the cyclotomic polynomial $cy(n,x)$ (the minimal polynomial of the complex algebraic number $e^{\frac{2\,\pi\,i}{n}}$, see \eg \cite{WatkinsZeitlin} or \cite{GKP}, p. 149, Exercise 50 a and b), evaluated at $x=-1$, which is  
\Beq
cy(n,-1)\sspeq  (-1)^{\varphi(n)}\,\prod_{{k=1}\atop {gcd(k,n)=1}}^{n-1}(1\sspp e^{2\pi\,i\frac{k}{n}}) \ ,\ n\sspeq 2,\,3,\,... 
\Eeq
We used the fact that the product in the definition of $cy(n,x)$ has $\varphi(n)$ factors (the degree as minimal polynomial). For $n=1$ one has $cy(1,-1)\sspeq -2$, which fits this formula if one defines $\varphi(1)\sspdef 1$ and takes the undefined product as $1$. One may rewrite this for 
$n\sspgeq 2$ by extracting $e^{\frac{\pi\, i}{n}}$, using for the sum of the $\varphi(n)$ terms the formula
\Beq
 {\frac{2}{n^2}}\, \sum_{{k=1}\atop {gcd(k,n)=1}}^{n-1}k\  \sspeq \  {\frac{\varphi(n)}{n}}\ ,\  n\sspeq 2,\,3,\,...  
\Eeq
This is a formula listed  by  {\sl R. Zumkeller} under \seqnum{A023022} which is found in\cite{Apostol}, p. 48, exercises 15 and 16, written such that both sides depend only on the distinct primes in the prime number factorization of $n$.  For the \rhs there is the well known formula \cite{Apostol}, p. 27,
\dstyle{{\frac{\varphi(n)}{n}} \sspeq \prod_{j=1}^{M(n)}\,\left(1\sspm \frac{1}{p_j}\right) \sspeq \sum_{d|n}\,\mu(d)\,{\frac{n}{d}}}, if $n$ has the distinct prime factors $p_j,\ j=1,2,...,M(n)\sspeq \seqnum{A001221(n)}\sspeq |sqfkset(n)|$ (see above for the notion  `squarefree kernel'). The {\sl M\"obius}  function $\mu$ entered here (see \seqnum{A008683}).  The PIE proof of eq. $(30)$, given in the {\it appendix A}, uses this form of  $\varphi(n)$. If one uses also $e^{i\,\pi}\sspeq -1$ the rewritten eq. $(29)$ becomes
\Beq
cy(n,-1)\sspeq (-1)^{{\frac{\varphi(n)}{2}}} \, \prod_{{k=1}\atop {gcd(k,n)=1}}^{n-1}\, 2\,\cos\left (\pi\,{\frac{k}{n}}\right)  \ ,\ n\sspeq 2,\,3,\,.... 
\Eeq
This can also be related to a special {\sl Sylvester} sequence. In general the  complex cyclotomic {\sl Sylvester}-numbers $Sy(a,b;n)$ are related to the sequence with three term recurrence $f_n\sspeq a\, f_{n-1} \sspp b\,f_{n-2}$, which has characteristic polynomial $x^2\sspm a\, x\sspm b$ with zeros $\alpha\sspequiv \alpha(a,b)\sspeq (a\sspp \sqrt{a^2\sspp 4\, b})/2$ and $\beta\sspequiv \beta(a,b)\sspeq (a\sspm \sqrt{a^2\sspp 4\, b})/2$. The definition is (see \eg \cite{mathworld})
\Beq
Sy(a,b;n)\sspdef \prod_{{k=1}\atop {gcd(k,n)=1}}^{n-1}\, (\alpha \sspm \beta\, e^{2\pi\, i\frac{k}{n}}) \ ,\ n\sspeq 2,3,... 
\Eeq
 For $n\sspeq 1$ one takes $Sy(a,b;1)\sspeq \alpha\sspm \beta$. For our purpose the values are $(a,b)\sspeq (0,-1)$, \ie $(\alpha,\beta)\sspeq (i,-i)$, and 
 \Beq
Sy(0,-1;n)\sspdef (-1)^{\frac{\varphi(n)}{2}}\, \prod_{{k=1}\atop {gcd(k,n)=1}}^{n-1}\, (1\sspp  e^{2\pi\, i\frac{k}{n}}) \ ,\ n\sspeq 2,\,3,\,..., 
\Eeq
with $Sy(0,-1;1)\sspeq 2\,i$.
This can be rewritten, using again eq. $(30)$ and $e^{i\,\pi}\sspeq -1$, as
\Beq
Sy(0,-1;n)\sspeq (-1)^{\varphi(n)}\, \prod_{{k=1}\atop {gcd(k,n)=1}}^{n-1}2\,\cos\left (\pi{\frac{k}{n}}\right)  \ ,\ n\sspeq 2,\,3,\,.... 
\Eeq
Therefore
\Beq
Sy(0,-1;n)\sspdef (-1)^{\frac{\varphi(n)}{2}}\, cy(n,-1)\ ,  n\sspeq 2,\,3,\,.... 
\Eeq
One could also include the $n=1$ case, with $\varphi(1)\sspdef 1$  if one takes \dstyle{{\frac{1}{i^{\varphi(n)}}}} instead of the prefactor.\psn
After these preliminaries back to the number $C(n,0)$. \psn
{\bf Proposition 6:} $\bf C(2\,m, 0),\ m\sspin \mathbb N$\psn
\Beq
C(2\,m,0)\sspeq (-1)^{{\frac{\varphi(4\,m)}{2}}\sspp \varphi(2\,m)}\,Sy(0,-1;2\,m)\sspeq   (-1)^{{\frac{\varphi(4\,m)}{2}}\sspp {\frac{\varphi(2\,m)}{2}}} \,cy(2\,m,-1)\ , \  m\in\mathbb N\, .
\Eeq
{\bf Proof:}  Use eqs. $(19)$ and $(35)$. The number of factors in eq. $(19)$ with $n\sspeq 2\,m$ is \dstyle{{\frac{\varphi(4\,m)}{2}}}, the degree $\delta(2\,m)$. Note that the restriction in the $Sy(0,-1;2\,m)$ product is  $gcd(k,2\,m)\sspeq 1$, while in the $C(2\,m,0)$ product it is $gcd(k,4\,m)\sspeq 1$, but this just says that only odd $k$s can contribute in both cases, and those odd numbers $\sspleq 2\,m\sspm 1$ dividing $4\,m$ or $2\,m$ have to be omitted. Both restrictions exclude the same odd numbers.\hskip 17cm $\square$
\psn
{\bf Corollary 1:}\  $\bf C(2\,p,0)$, $\bf p$ an odd {\bf a prime} \psn
\Beq
{\text {For odd  prime}}\ p:\ \  C(2\,p,0)\sspeq (-1)^{{\frac{p-1}{2}}}\, p\ .
\Eeq
{\bf Proof:}\psn
This follows immediately from eqs. $(19)$ and $(31)$ if one uses the known formulae $cy(2\,p,x)\sspeq cy(p,-x)$, $cy(p,1)\sspeq p$, and evaluates the $\varphi$ functions to obtain the sign. This formula follows from \dstyle{cy(n,x)\sspeq \prod_{d|n}\, (x^d\sspm 1)^{\mu(\frac{n}{d})}} where the multiplicative {\sl M\"obius} function (with $\mu(p)\sspeq -1$) enters (see, \eg \cite{GKP}, exercise 50 b., solution p.506).
\psn 
{\bf Proposition 7:}\  $\bf C(2^m,x),\ m\sspin \mathbb N_0$\psn
\Beq
C(2^m,x) \sspeq {\hat t}(2^{m-1},x)\ ,\  m\sspin \mathbb N,
\Eeq
with the integer polynomials $\hat t$ defined in eqs. $(2)$ and $(3)$ in terms of {\sl Chebyshev} $T$- or $S-$polynomials.
\psn
{\bf Proof:} Use eq. $(16)$ and the {\it Note added} to the link \cite{Lang1}, where it was shown in eq. $(8)$ that\pn
\dstyle{2^{2^{m-1}}\, \Psi\left(2^{m+1},{\frac{x}{2}}\right)\sspeq 2\,T\left(2^{m-1},{\frac{x}{2}}\right) \sspfed {\hat t}(2^{m-1},x)}, for $m\sspeq 1,2,...\,$. 
It is clear that $C(1,x)\sspeq x+2$. This {\it proposition} will later be generalized in {\it theorem 1A}. \hskip 9cm $\square$
\psn
{\bf Corollary 2:}\  $\bf C(2^m,0),\ m\sspin \mathbb N_0  $\psn
From the {\it proposition 7} and the known fact that $T(2\,n,0)\sspeq (-1)^n$ follows that
\Beq
C(1,0)\sspeq +2, \ \ C(2,0)\sspeq 0,\ \ C(2^2,0)\sspeq -2,\ \  C(2^m,0)\sspeq +2, \  m\sspgeq 3\ . 
\Eeq
Note that a standard formula for the cyclotomic polynomials $cy(2^m,x)$ which involves divisors (see \eg \cite{GKP}, p. 149, Exercise 50 b, with $\Psi_m(x) \sspequiv cy(m,x)$) leads to undetermined expressions for $cy(2^m,-1)$ if  $m\sspgeq 1$. 
\psn
For the following {\it theorem 1A} on a formula for $C( n,x)$ for even $n$, based on the divisor product representation ($dpr$) of numbers \cite{Lang2}, we need  two {\it lemmata}.\psn
{\bf Lemma 3:}\  {\bf Pairing in $\bf dpr(2\,m)$}\psn
In the divisor product representation of an even number  $2\,m$, $m\in \mathbb N$, for each $a$-factor in the numerator (resp. denominator) one finds exactly one $a$-factor in the denominator (resp. numerator) with argument ratio either $1:2$ or $2:1$.  \psn
{\bf Proof:}\psn
Due to the $*$-multiplication property of $dpr$s (see the {\it theorem} in \cite{Lang2})
it is sufficient to prove this {\it lemma} for $dpr(2\,p_1\,\cdots\, p_N)$ with odd primes $p_1,\,...,\, p_N$. Indeed, if the prime factorization of $2\,m$ is $2^{k+1}\,p_1^{e_1}\,\cdots\, p_N^{e_N}$, with $k\sspin \mathbb N_0$ and odd primes $p_j$ and positive exponents $e_j$, for $j=1,\,...,\,N$, then $dpr(2\,m) \sspeq (2^k\, p_1^{e_1-1}\,\cdots p_N^{e_N-1})*dpr(2\,p_1\cdots p_N)$ and the pairing property of the latter $dpr$ will be preserved after multiplication of each argument. Now the proof of  the pairing property for $dpr(2\,p_1\cdots p_N)$ becomes elementary, once one splits the products in eq. $(7)$ of \cite{Lang2} (note that $N$ in this reference is now  $N+1$) into those $a$-factors which are even (\ie contain the even prime $2$) and those which are odd (\ie those composed of only odd primes). Indicate $a$-factors with odd arguments, being the product of $k$ odd primes, by $a_o(.(k).)$.   Then any $a$-factor in the numerator (resp. denominator) with some even argument $2\,n$ is found in the denominator (resp. numerator) exactly once with the argument $n$. This is because for each factor in the numerator product $\Pi\, a(2\,.(N- 2\,j).)$, with \dstyle{j\sspin \{1,...,\floor{\frac{N}{2}}\}}, say, $a(2\,p_{i_1}\cdots p_{i_{N-2\,j}})$, there is exactly one factor in the denominator product $\Pi\, a_o(.(N-2\,j).)$, namely the one which uses the same $N-2\,j$  odd primes $a(p_{i_1}\cdots p_{i_{N-2\,j}})$. Similarly, for each $a$-factor in the denominator  product $\Pi\, a(2\,.(N-(2\,j+1)).)$, with \dstyle{j\sspin \{0,...,\floor{\frac{N-1}{2}}\}}, there is exactly one factor in the numerator product $\Pi\,a_o(.(N-(2\,j+1)).)$. The first factor in the numerator, $a(2\,p_1\,\cdots\, p_N)$  $(j=0)$ is paired with the single odd argument $a$-factor of the first product in the denominator $\Pi\, a_o(.(N).)\sspeq a(p_1\,\cdots\, p_N)$ . \pn
Note, that the example for the primorial \seqnum{A002110}$(5)\sspeq 2310$, given in \cite{Lang2}, {\it table}\, $2$ shows that the pairing occurs in general not between numerator and denominator factors at the same position, given the ordering prescription indicated by $\cal O$ (falling arguments) in eq. $(7)$ in \cite{Lang2}. See the fifth position in the $N\sspeq 5$ case there.\ps
Recall  also, for later purposes, that because $dpr(2\,m)$ has the same number of factors in the numerator and in the denominator, \viz\  $2^N$ ( this balance holds true for all $dpr(n)$, $n\sspgeq 2$, due to {\it proposition 5} in \cite{Lang2}), after the split into even and odd arguments the number of pairs with even argument in the numerator matches the one with even argument in the denominator. This will later lead to a balanced formula for $C(2\,m,x)$  in terms of $\hat t$-polynomials (the same number of $\hat t$s in the numerator and in the denominator, \viz\  $2^{N-1}$). \hskip 17.5cm $\square$  
\pbn
{\bf Lemma 4:}  \pbn
\Beq
{\frac{{\hat t}({\frac{n}{2}}\sspp 1,x)\sspm {\hat t}({\frac{n}{2}}\sspm 1,x)}{{\hat t}({\frac{n}{4}}\sspp 1,x)\sspm {\hat t}({\frac{n}{4}}\sspm 1,x)}} \sspeq {\frac{S({\frac{n}{2}}\sspp 1,x)}{S({\frac{n}{4}}\sspm 1,x)}} \sspeq 2\,T\left({\frac{n}{4}} ,{\frac{x}{2}}\right) \sspeq {\hat t}\left({\frac{n}{4}},x\right) \  , \ {\text {for}}\ n\sspequiv 0\,(mod\, 4)\ .
\Eeq
This identity involves the integer polynomials $\hat t$ introduced earlier in eqs. $(2)$ and $(3)$. \psn
{\bf Proof:} \psn
Use the known identity for {\sl Chebyshev} polynomials, \cite{MOS}, p. 261, first line, specialized,  $T(m+1,x)\sspm T(m-1,x) \sspeq 2\,(x^2-1)\,U(m-1,x)$, written for $x$ replaced by \dstyle{\frac{x}{2}}, with \dstyle{S(m-1,x)\sspeq U(m-1,\frac{x}{2})} and the definition of the $\hat t$-polynomials from eq. $(2)$. For  $\ n\sspequiv 0\,(mod\, 4)$ this identity can be used in the numerator as well as in the denominator. The $n-$independent factor \dstyle{({\frac{x}{2}})^2\sspm 1} drops out if $x\sspneq \pm 2$. However, the lemma holds also for these $x-$values as can be seen after applying {l'H\^opital}'s rule, using the well known identity $T^{\prime}(x,n)\sspeq n\,U(n-1,x)$. In the second to last step the well known identity, \cite{MOS} p. 260, last line, written as \dstyle{2\, T\left(m,{\frac{x}{2}}\right)\, S(m-1,x)\sspeq S(2\,m\sspm 1,x)} has been employed. \hskip 5cm  $\square$\pbn
 {\bf Theorem 1A:\ $\bf C(2\,m,x),\ m\sspgeq 1$, in terms of $\bf \hat t$-polynomials}\psn
With the prime number factorization of \dstyle{2\,m \sspeq 2^k\,p_1^{e_1}\,\cdots p_N^{e_N}} with $k\sspin \mathbb N$,  odd primes $p_j$ and positive exponents $e_j$, $j=1,...,N$, one has
\Beq
C(2\,m,x)\sspeq {\frac{\hat t(2^{k-1}\,p_1^{e_1}\cdots p_N^{e_N},x)\, \Pi\,\hat t(2^{k-1}\,_*(N-2).,x)\, \Pi\, \hat t(2^{k-1}\,_*(N-4).,x)\cdots }{\Pi\,\hat t(2^{k-1}\,_*(N-1).,x)\, \Pi\, \hat t(2^{k-1}\,_*(N-3).,x)\cdots \hskip 3.5cm}}\,.
\Eeq
with $\hat t(2^{k-1}\,_*(0).,x)\sspeq \hat t(2^{k-1}\,p_1^{e_1-1}\,p_2^{e_2-1}\cdots p_N^{e_N-1},x)$, and the products $\Pi\, \hat t(2^{k-1}\,_*(K).,x)$ are over the \dstyle{\binomial{N}{K}} factors with $K$ primes from the set $\{p_1,p_2,...,p_N\}$ multiplied by $p_1^{e_1-1}\,p_2^{e_2-1}\cdots p_N^{e_N-1}$, \ie  $K$-products of the form $\hat t(2^{k-1}\, p_1^{e_1-1}\,p_2^{e_2-1}\cdots p_N^{e_N-1}\,   p_{i_1}\,p_{i_2}\cdots p_{i_K},x)$. We used $_*(k).$ instead of $.(k).$ for the indices of the $\hat t$ polynomials to remind one of the extra factor (without the powers of $2$ because the even prime has been extracted everywhere) due to  the $*$-multiplication.\psn
Before we give the proof an example will illustrate this {\it theorem}.\psn
{\bf Example 3:}\  $\bf k\sspeq 1$, $\bf N\sspeq 2$\psn
\Beq
C(1350,x)\sspeq C(2\cdot 3^3\,5^2,x)\sspeq {\frac{\hat t(3^3\,5^2,x)\,\hat t(3^2\,5,x)}{\hat t(3^3\,5,x)\,\hat t(3^2\,5^2,x )}}\sspeq {\frac{\hat t(675,x)\, \hat t(45,x) }{\hat t(135,x)\, \hat t(225,x)}}\ .
\Eeq
{\bf Proof:}\psn
 This is based on \cite{Lang2} applied for the minimal polynomials $\Psi\left (4\,m,\frac{x}{2}\right)$\, (see the definition for $C(2\,m,x)$ given in eq. $(16)$). Let $2^{k+1}\,p_1^{e_1}\,\cdots p_N^{e_N}$ be the prime number factorization for $4\,m$. Only the squarefree kernel $2\,p_1\cdots p_N$  is important due to the {\it theorem} in \cite{Lang2} for the divisor product representations ($dpr$s) and the $*$-multiplication. This means that one has to multiply after the computation of $\Psi(2\,p_1\cdots p_N)$, using {\it proposition 1} of \cite{Lang2}, the first argument, the index in conventional notation, of each factor $t(n_i,x)$ in the the numerator, and of each $t(m_i,x)$ in the denominator (see \cite{Lang2}, eq. $(1)$) with the number  $2^k\, p_1^{e_1-1}\,p_2^{e_2-1}\cdots p_N^{e_N-1} $. Here the pairing {\it lemma 3} is crucial which carries over to the indices of the $t$-factors in the numerator and denominator of \cite{Lang2}, eq. $(1)$. Because there is always at least one factor of $2$ in the number which multiplies every $t$ index after the $*$-multiplication, the pairing will occur always between indices which are $ 0\,(mod\, 4)$ and $0\, (mod\, 2)$. Then only the $n$ even alternative in the definition of \dstyle{t\left(n,\frac{x}{2}\right)} from eq. $(2)$ of \cite{Lang2}, rewritten here as eq. $(18)$, is relevant. Observe that the prefactors in this definition of  \dstyle{t\left(n,\frac{x}{2}\right)} are not of interest, provided we take \dstyle{2\, \left(T\left({\frac{n}{2}}\sspp 1,\frac{x}{2}\right)\sspm T\left({\frac{n}{2}}\sspm 1,\frac{x}{2}\right)\right)} which is the monic integer polynomial \dstyle{{\hat t}\left({\frac{n}{2}}\sspp 1,x\right)\sspm {\hat t}\left({\frac{n}{2}}\sspm 1,x\right)}. See eq. $(18)$, and the discussion in connection with {\it lemma 2}. $C(2\,m,x)$ is written here as a rational function of monic integer polynomials. For each replaced  \dstyle{t\left(n,\frac{x}{2}\right){\Big /}t\left({\frac{n}{2}},\frac{x}{2}\right)} (either in the numerator or denominator, depending on where the larger index appears),  one can apply {\it lemma 4}. This is how for each such $t$-quotient one obtains \dstyle{{\hat t}\left({\frac{n}{4}},x\right)}. Here the fact that the larger index of a pair is always equivalent to $0 \,(mod\, 4)$ is important. From the structure of the numerator and denominator of the original $dpr(2\,m)$ with the separation of the even and odd $a$-arguments (which carries over to the $t$-indices),  discussed in the proof of {\it lemma 3}, one now finds the numerators and denominators of the theorem. Just search the products for the even $t$-indices. \Eg  \dstyle{\Pi\, {\hat t}(2^{k-1}\,_*(N-2).,x)} in the numerator of the {\it theorem} originates from the quotient of products \dstyle{{\frac{\Pi\,t(2\,.(N-2).,{\frac{x}{2}})}{\Pi\, t_o(.(N-2). ,{\frac{x}{2}})}}}  before the $*$-multiplication has to be applied. This leads, with eq. $(18)$ and {\it lemma 4},  to \dstyle{\Pi\, {\hat t}\left({\frac{2^{k+1}}{4}}\,_*(N-2).,x\right)} after multiplying each $\hat t$ index in the product with $2^k\, p_1^{e_1-1}\,p_2^{e_2-1}\cdots p_N^{e_N-1} $, where the $*$ reminds one to  multiply each index with this number divided by $2^k$. \ps
As announced at the end of the proof of {\it lemma 3} the number of $\hat t$ polynomials in the numerator is the same as the one for the denominator, \viz\   $2^{N-1}$, which is due to the sum in the {\sl Pascal}-triangle \seqnum{A007318} row No. $N$ over even, resp. odd numbered positions. \hskip 14cm $\square$       
\psn
Note that this {\it theorem 1A}, evaluated for $x\sspeq 0$ leads in general to undetermined expressions, remembering that ${\hat t}(n,0)\sspeq 0$ if $n$ is odd, and \dstyle{2\,(-1)^{\frac{n}{2}}} if $n$ is even. However, a correct evaluation (using {\sl l'H\^opital}'s rule) has to reproduce the result known from {\it corollary 2}. 
\pbn 
 The following factorization of the monic integer $\hat t$-polynomials is related to {\it theorem 1A}.\psn
{\bf Theorem 1B: \  Factorization of $\bf \hat t$-polynomials in terms of the minimal $\bf C$-polynomials} \pbn
\Beq
{\hat t}(n,x) \sspeq \prod_{d|op(n)}\, C(2\,n/d,x)\sspeq \prod_{d|op(n)}\,C(2^{k+1}\,d,x)\ ,\   
\Eeq
with $op(n)\sspeq$\seqnum{A000265}$(n)$, the odd part of $n$, and $2^k$ is the largest power of $2$ dividing  $n$. The exponents are  $k\sspeq k(n)\sspeq$\seqnum{A007814}$(n)$, $k\sspin \mathbb N_0\,$.
\pbn 
Before the proof we give an example.\psn
{\bf Example 4:}\  $\bf n=10$\psn
\Beqarray
 &&-2\sspp 25\,x^2\sspm 50\,x^4\sspp 35\,x^6\sspm 10\,x^8\sspp x^{10} \sspeq {\hat t}(10,x) \sspeq C(20,x)\,C(4,x)  \speq \nonumber \\
&&  (x^8\sspm 8\,x^6\sspp 19\,x^4\sspm 12\,x^2\sspp 1)\, (x^2\sspm 2)\,\ .
\Eeqarray
{\bf Proof:}  This is modeled after a similar proof in \cite{WatkinsZeitlin}.\psn
It is clear that both sides are monic integer polynomials. See the definition of $\hat t$ in eqs. $(2)$, $(3)$ and {\it proposition 2}.\pn
In order to check the degree we consider the even and odd $n$ cases separately. If $n\speq 2^{k(n)}\,op(n)$ with $k(n)\sspeq$\seqnum{A007814}$(n)\sspgeq 1$  then \dstyle{\sum_{d|op(n)}\,\delta(2^{k+1}\,d) \sspeq  \sum_{d|op(n)}\,{\frac{\varphi(2^{k+2}\,d)}{2}}} from the known degree $\delta$ of the $C$-polynomials, here for index $>1$. Due to well known properties of the {\sl Euler} totient function (see \eg \cite{Apostol}, {\it Theorem} 2.5, (a) and (c), p. 28 and {\it Theorem} 2.2., p. 26) this leads to \dstyle{{\frac{2^{k(n)+2}-2^{k(n)+1}}{2}}\,op(n) \sspeq n}, the degree of $\hat t(n,x)$. In the odd case, if $n\sspeq op(n)$ then \dstyle{\sum_{d|n}\, {\frac{\varphi(4\,d)}{2}}  \sspeq \sum_{d|n}\,d \sspeq n}, again the degree of $\hat t(n,x)$. \pn
For the proof one compares the zeros of both sides. The zeros of $\hat t(n,x)$ are \dstyle{\hat x^{(n)}_l\sspeq 2\,cos\left((2\,l+1)\,\frac{\pi}{2\,n}\right)} for $l \sspeq 0,\,1,\,...,\,n-1$. This is known from the zeros of the  {\sl Chebyshev}  $T$-polynomials. The zeros of the $C(2^{k+1}\,d,x)$-polynomials are used in the form given by eq. $(20)$:  \dstyle{2\,cos\left((2\,l'+1)\,\frac{\pi}{2^{k+1}\,d}\right)} for $l' \sspeq 0,1,...,2^k\, d-1$, where $gcd(2\,l'+1,2^{k+1}\, d)\sspeq 1$, \ie $gcd(2\,l'+1, d)\sspeq 1$. Because the degrees match, it is sufficient to show that each zero of $\hat t(n,x)$ occurs on the {\it r.h.s.}. For $n\sspeq 2^k(n)\,op(n)$, with $k(n)\sspeq$\seqnum{A007814}$(n)$, consider $gcd(2\,l+1,2\,n) \sspeq  gcd(2\,l+1,2^{k+1}\,op(n)) \sspeq gcd(2\,l+1,op(n)) \sspeq g$ (some odd number). Hence $2\,l+1\sspeq (2\,l'+1)\,g$, with some $l'$ and $op(n)\sspeq d\,g$, \ie $d|op(n)$. Therefore \dstyle{{\frac{2\,l+1}{2\,n}}\sspeq {\frac{2\,l'+1}{2^{k+1}\,d}}}. Thus the $\hat t(n,x)$ zero  \dstyle{2\,cos\left((2\,l+1)\,\frac{\pi}{2\,n}\right)} appears on the {\it r.h.s.} as one of the $C(2^{k+1}\,d,x)$ zeros.\hskip 18cm $\square$.\psn
{\bf Remark 3:\ }{\bf  Derivation of Theorem 1A from Theorem 1B}\psn
{\it Theorem 1B} can be used as recurrence for the $C$-polynomials in terms of the $\hat t$-polynomials. This is similar to the case treated in \cite{WatkinsZeitlin} for the minimal polynomials $\Psi$ of \dstyle{{\frac{2\,\pi}{n}}} (see \seqnum{A181875}/\seqnum{A181876} for their coefficients). The solution of this recurrence has been given in \cite{Lang2}. Indeed, {\it theorem 1A} has been derived above from this solution. Originally we found {\it theorem 1B} starting from {\it theorem 1A} building up iteratively  a formula for \dstyle{ {\hat t}(2^{k-1}\,p_1^{e_1}\cdots p_N^{e_N},x)} in terms of the $C$-polynomials. We give this formula as a corollary.\psn
{\bf Corollary 3:\  \  $\bf \hat t$-polynomials in terms of $\bf C$-polynomials}\psn
 \Beq
 {\hat t}(2^{k-1}\,p_1^{e_1}\cdots p_N^{e_N},x)\sspeq \prod_{q_1=0}^{e_1}\,\prod_{q_2=0}^{e_2}\, \cdots \prod_{q_N=0}^{e_N}\, C(2^k\,p_1^{q_1}\,p_2^{q_2}\cdots p_n^{q_N},x)\ , k\sspeq \mathbb N,\   N\sspin \mathbb N_0\ .
\Eeq  
In order to find a  simplified expression for $C(2\,m+1,x), m\sspin \mathbb N_0$ we need the following {\it lemma} in order to rewrite the quotient \dstyle{{\frac{t(even,\frac{x}{2})}{t(odd,\frac{x}{2})}}}, given the definition eq. $(18)$.\psn
{\bf Lemma 5:}\psn
for $n\sspeq 2\,M\sspp 1, M\sspin \mathbb N_0$ one has
\Beqarray
{\frac{ t\left(2\,n,{\frac{x}{2}}\right)}{t\left(n,{\frac{x}{2}}\right)}}&\sspeq& 2^{M-n}\,{\frac{{\hat t}(n+1, {\frac{x}{2}}) \sspm {\hat t}(n-1, {\frac{x}{2}})}{ {\hat t}(M+1, {\frac{x}{2}}) \sspm {\hat t}(M, {\frac{x}{2}}) }}\sspeq {\frac{1}{2^{M+1}}}\, {\frac{(x^2\sspm 4)\,S(n-1,x)\hskip 1cm}{(x\sspm 2)\,S(2\,M,\sqrt{2\sspp x}\,)}}  \\
&\sspeq & {\frac{x\sspp 2}{2^{M+1}}}\,{\frac{S(n-1,x)\hskip 1.3cm}{S(n-1,\sqrt{2\sspp x}\,)}}\ .  
\Eeqarray 
Here new important monic integer polynomials enter the stage:\psn
{\bf Definition 1:\ \  $\bf q$-polynomials}\psn 
With $n \sspin  \mathbb N_0$ define
\Beq
q(n,x)\sspdef {\frac{S(2\,n,x)}{S(2\,n,\sqrt{2\sspp x}\,)}} \ .
\Eeq
That this defines indeed monic integer polynomials of degree $n$ is shown by the next {\it lemma}.\psn
{\bf Lemma 6:}\psn
\Beq
q(n,x)\sspeq  (-1)^n\,S(2\,n\,,\sqrt{2\sspm x}\,) \sspeq S(n,x)\sspm S(n-1,x) \ .
\Eeq
{\bf Proof:}\psn
It is known from the \ogf of the {\sl Chebyshev} $S$-polynomials that the bisection yields \dstyle{S(2\,n,y)\sspeq S(n,y^2\sspm 2)\sspp S(n-1,y^2\sspm 2)}, or \dstyle{S(2\,n\,,\sqrt{x\sspp 2}\,)\sspeq S(n,x)\sspp S(n-1,x)}. With $x$ replaced by $-x$ one has \dstyle{S(2\,n,\sqrt{2\sspm x}\,)\sspeq (-1)^n\,(S(n,x)\sspm S(n-1,x))} which explains the second equation of the {\it lemma}.  Therefore one has to prove $S(2\,n,x)\sspeq S(n,x)^2\sspm S(n-1,x)^2$. This identity can, for example,  be proved using the \ogf for the square of the $S$-polynomials (for their  coefficient table see \seqnum{A181878}, also for the paper \cite{Lang3} given there as a link). This computation was based on the {\sl Binet}-{\sl de Moivre} formula for the $S-$polynomials. This \ogf was found to be \dstyle{{\frac{1\sspp z}{1\sspm z}}\, {\frac{1}{1\sspp (2\sspm x^2)\,z\sspp z^2}}}. The \ogf for $\{S(2\,n,x)\}_{n=0}^{\infty}$ is \dstyle{(1\sspp z)/(1\sspp(2\sspm x^2)\, z\sspp z^2)} (from the bisection). Because the \ogf for $ S(n,x)^2\sspm S(n-1,x)^2$ is $(1+z)$ times the one for $S(n,x)^2$ the completion of the proof is then obvious.\hskip 8cm $\square$\psn
{\bf Remark 4 :\ \  O.g.f. and coefficient array for the $\bf q$-polynomials}\psn
From {\it lemma 6} the \ogf \dstyle{Q(z,x):=\sum_{n=0}^{\infty}\,q(n,x)\,z^n \sspeq (1\sspm z)/(1-x\,z\sspp z^2)} from the known \ogf of the $S$-polynomials. This shows that the $q(n,x)$ coefficients constitute a {\sl Riordan}-array (infinite lower triangular ordinary convolution matrix), which is in standard notation \dstyle{\left ({\frac{1\sspm x}{1\sspp x^2}}, {\frac{x}{1\sspp x^2}}\right )}, meaning that the \ogf of the column No. $m$ sequence is \dstyle{{\frac{1\sspm x}{1\sspp x^2}}\,\left ({\frac{x}{1\sspp x^2}}\right )^m }. This is the triangle \seqnum{A130777} where more information can be found. For example, the explicit form for the coefficients is \dstyle{ Q(n,m)\sspeq (-1)^{{\frac{n-m+1}{2}}}\ {\binomial{{\frac{n+m}{2}}}{m}} } if $n\sspgeq m\sspgeq 0$ and $0$ otherwise.\psn
 The pairing {\it lemma 3} will also be used in the proof of the following {\it theorem}.\psn
{\bf Theorem 2A:\  $\bf C(n,x)$,\  $\bf n \sspgeq 3$, odd, in terms of $\bf q$-polynomials}
\pbn
With the prime number factorization of \dstyle{n\sspeq p_1^{e_1}\,\cdots p_N^{e_N}} with  odd primes $p_1, ...,p_N$ and positive exponents $e_j$, $j=1,...,N$, one has
\Beq
C(n,x)\sspeq {\frac{q\left({\frac{n-1}{2}},x\right)\, \prod_{i_1<i_2}\,q\left({\frac{n/(p_{i_1}\,p_{i_2})-1}{2}},x\right)\, \prod_{i_1<i_2<i_3<i_4}\, q\left({\frac{n/(p_{i_1}\,p_{i_2}\,p_{i_3}\,p_{i_4})-1}{2}},x\right) \cdots}{\prod_{i_1}\,q\left({\frac{n/p_{i_1}-1}{2}},x\right)\, \prod_{i_1<i_2<i_3}\,q\left({\frac{n/(p_{i_1}\,p_{i_2}\,p_{i_3})-1}{2}},x\right)\, \cdots\hskip 3.5cm}}\  .
\Eeq
\psn
Here each index $i_j$ of the products runs from $1$ to $N$.\pn
Before we give the proof an examples will illustrates this {\it theorem}.\psn
{\bf Example 5:}  $\bf n\sspeq 3^2\cdot 5\sspeq 45$\psn
\Beq
C(45,x)\sspeq C(3^2\cdot 5,x)\sspeq {\frac{q(22,x)\, q(1,x)}{q(7,x)\,q(4,x)}} \ . 
\Eeq  
This checks.\psn
{\bf Proof:} This is analogous to the proof of {\it theorem 1A}. Again the $dpr(2\,n)$ 
representation, from which one derives $\Psi(2\,n,x)$, with the prime number factorization for the odd $n$, is considered. The $*$-multiplication property allows to consider $dpr(2\, p_1\cdots p_N)$ with a subsequent multiplication of all arguments in the numerator and denominator by $p_1^{e_1-1}\cdots p_N^{e_N-1}$. Because of the paired numerator/denominator structure due to {\it lemma 3} one finds, either in the numerator or in the denominator quotients of the type \dstyle{{\frac{t(2\,k,{\frac{x}{2}})}{t(k,{\frac{x}{2}})}}} which are up to a factor \dstyle{2^{{\frac{k+1}{2}}}\, (x\sspp 2)} equal to \dstyle{q\left({\frac{n-1}{2}},x\right)}. Factors of powers of $2$ are irrelevant (they have to cancel) because on both sides of {\it theorem} monic polynomials appear. The factors of $x-2$ also cancel because the number of $q$-polynomials in the numerator and denominator have also to match (see the remark at the end of the proof  of {\it lemma 3}). This number is $2^{N-1}$. The structure of the numerator of the {\it theorem} originates from products  with the even indexed \dstyle{t\left(2k,\frac{x}{2}\right) } polynomials in the numerator after $*$- multiplication. Before this multiplication on has in the numerator \dstyle{t\left(2\,p_1\cdots p_N,{\frac{x}{2}}\right)} and $t-$products over all possibilities to leave out $2,\, 4,\,...$ of the odd primes from the set $\{p_1,\, ...,\,  p_N\}$. The $*$-multiplication then leads to  \dstyle{t\left(n,{\frac{x}{2}}\right)} and $t-$products over all possibilities to divide $n$ by these $2, 4,...$ odd primes. Together with the pairing partners from the denominator this leads to the $q$-polynomials given in the numerator of the {\it theorem}. A similar argument produces the denominator $q$-polynomials. The number of these $q-$polynomials  in the numerator, \viz\ $2^{N-1}$ , matches the one for the denominator. \hskip 7cm $\square$
\pbn
{\bf Theorem 2B:\   \  Factorization of $\bf q$-polynomials  in terms of $\bf C$-polynomials}
\pbn
\Beq
q(n,x)\speq \prod_{1\sspkl d|(2\,n+1)}\ C(d,x)\ , \ n\sspin \mathbb N\ . 
\Eeq
{\bf Example 6:}  $\bf n\sspeq17 $\psn
\Beq
q(17,x)\sspeq C(5,x)\, C(7,x)\,C(35,x) \ .
\Eeq
This checks.\psn
{\bf Proof:} \psn 
This is again modeled after a similar proof in \cite{WatkinsZeitlin}. Compare this with the proof of {\it theorem 1B}.\psn
It is clear that both sides are monic integer polynomials, and the degree fits due to 
\dstyle{n\sspeq \sum_{d|(2\,n+1)}\,{\frac{\varphi(d)}{2}}  \sspm  {\frac{1}{2}}\,\varphi(1)}. See the properties of the {\sl Euler} totient function mentioned above in the proof of {\it theorem 1B}. \pn
The zeros of $q(n,x)$ are \dstyle{x^{(n)}_l\sspeq 2\,\cos\left( \pi\,{\frac{2l+1}{2n+1}}\right)}, for $l\sspeq 0,1,...,n-1$. This follows from {\it definition 1} and {\it lemma 6} 
with the zeros of \dstyle{S(2\,n,\sqrt{2\sspm x}\,)} which are known from those of {\sl Chebyshev} $S$-polynomials. To show that each of these zeros appears on the \rhs for $C(d,x)$ with some $d|(2\,n+1), d\sspneq 1$, \ie as \dstyle{2\, \cos\left( \pi\,{\frac{2l'+1}{q}}\right)} for some  \dstyle{l'\sspin\{0,\,1,\,...,\,\frac{d-3}{2}\}} and  $gcd(2l'+1,d)\sspeq 1$ (see eq. $20$), let $gcd(2l+1,2n+1)\sspeq g$, with some odd $g$ with $2\,n+1\sspeq d\,g$. Then $2\,l\sspp 1\sspeq (2\,l'\sspp 1)\,g\sspeq (2\,l'\sspp 1){\frac{2n+1}{d}}$, and for each  $l\sspin \{0,\,1,\,...,\,n-1\}$ there is one $l'\sspin \{0,\,1,\,...\,,\frac{d-3}{2}\}$. \hskip 8 cm $\square$\pbn
{\bf Remark 5:\ }{\bf  Derivation of Theorem 2A from Theorem 2B}\psn 
{\it Theorem 2B} can be used as recurrence for the $C$-polynomials in terms of the $q$-polynomials. See the {\it remark 3}. We give the solution in the following {\it corollary}.\psn
{\bf Corollary 4:\  \  $\bf q$-polynomials in terms of $\bf C$-polynomials}\psn
With the prime number factorization of $2\,n+1\sspeq p_1^{e_1}\cdots p_N^{e_N}$, with odd primes, $n\sspin \mathbb N$, one has
\Beq
q(n,x)\sspeq \prod_{q_1=0}^{e_1}\,\cdots \prod_{q_N=0}^{e_N}\, C(p_1^{q_1}\cdots p_N^{q_N},x)/C(1,x) \ .
\Eeq
The division by $C(1,x)\sspeq x\sspp 2$ was necessary because not all $q_j$-indices were originally allowed to vanish.\pbn
These rational representations of $C(n,x)$ do not lend itself to extraction of the value $C(n,0)$ because undetermined \dstyle{{\frac{0}{0}}} quotients appear. In the following we give the absolute term of $C$ for prime indices.\psn
{\bf Proposition 8:} \ $\bf C(p,0)$\psn
For $n\sspeq p$, a prime, one has
\Beq
C(p,0)\speq  {\Cases3 {0} {{\text if}\ \ p\sspeq 2\ ,} {(-1)^{{\frac{p-1}{4}}}} {{\text if}\ \ p\sspeqv 1\,(mod\,4)\ ,} {(-1)^{{\frac{p+1}{4}}}} {{\text if}\ \ p\sspeqv 3\,(mod\,4)\ . } }
\Eeq
{\bf Proof:} In  eq. $(31)$ with $n\sspeq p$ there is no $gcd$-restriction on the product. Therefore, one can use a known formula  (see {\it appendix B} for a proof)
\Beq
\prod_{k=1}^{n-1}\,2\,\cos\left(\pi\,{\frac{k}{n}}\right)\sspeq {\Caseszwei{(-1)^{{\frac{n-1}{2}}}}{{\text if}\ n\ \text{is odd},}{0}{{\text if}\ n\ \text{is even}, }}
\Eeq 
to obtain $cy(2,-1)\sspeq 0$ (which is also clear from the definition: $cy(2,x)\sspeq x\sspp 1$), and for odd primes $cy(p,-1)\sspeq (-1)^{p-1}\sspeq +1$. Because all $k$ contribute in the product, one can use for even $k$ the formula \dstyle{\cos\left(\pi\,{\frac{2\,K}{p}}\right)\sspeq -\cos\left(\pi\,{\frac{p\sspm 2\,K}{p}}\right)} which shows that 
one generates again all odd $k$ contributions, however each with a minus sign. This leads to \dstyle{cy(p,-1)\sspeq (-1)^{p-1}\ \prod_{l=0}^{{\frac{p-3}{2}}}\, \left(2\, \cos\left(\pi\,{\frac{2\,l\sspp 1}{p}}\right)\right)^2\sspeq C(p,0)^2}. Thus $C(p,0)\sspeq \pm 1$. The sign of $C(p,0)$ is $(-1)^{\delta_+(p)}$ from eq. $(19)$ for $n\sspeq p$. This number $\delta_+(p)$ of positive zeros has been found  in {\it proposition 5} to be $\frac{p-1}{4}$ if $p\sspequiv 1\,(mod\,4)$, and $\frac{p+1}{4}$ if $p\sspequiv 3\,(mod\,4)$. For $C(2,0)\sspeq 0$, the sign is, of course, irrelevant. \hskip 10cm $\square$  \psn
We close this section with a conjecture on the discriminant of the $C-$polynomials. 
The discriminant of a monic polynomial $p$  of degree $n$ can  be written as the square of the determinant of  an $n\times n$ {\sl Vandermonde} matrix $V_n(x_1^{(n)},\,...,\,x^{(n)}_n )$ with elements $(V_n)_{i,j} \sspdef (x^{(n)}_i)^j, \, \sspeq 1,\,...,\,n$ and $j\sspeq 0,\,...,\,n-1$ with the zeros $x^{(n)}_i$ of $p$ . Here the $\delta(n)$ zeros of $C_n$ are given in eq. $(20)$. Another formula for the discriminant  (see \eg \cite{Rivlin}, Theorem 5.1, p. 218) is in terms of the derivative $C^{\prime}(n,x)$ and the zeros of $C$: \dstyle{(-1)^{\frac{\delta(n)\,(\delta(n)-1)}{2}}\, \prod_{i=1}^{\delta(n)}\,C^{\prime}(n,x^{(n)}_i)}. The result is the sequence $Discr(C(n,x))\sspeq$ \seqnum{A193681}$(n)\sspeq [1,1,1,8,5,12,49,2048,81,2000,14641,2304,371293,...]$ for $n\sspgeq 1$. The following conjecture is on the sequence \dstyle{q(n)\sspeq {\frac{ n^{\delta(n)}}{Discr(C(n,x))}}} which is \seqnum{A215041}, $[1, 2, 3, 2, 5, 3, 7, 2, 9, 5, 11, 9, 13, 7, 45, 2, 17, 27, 19, 25,...]$.\psn
{\bf Conjecture:\ explicit form of the $\bf q$-sequence}\psn
{\bf o)}\ \ $q(1)=1$ (clear).  \psn
{\bf i)} \ \ If $n\sspeq 2^k$ for $ k\sspin \mathbb N$ then $q(n)\sspconj 2$.\psn  
{\bf ii)}\ \ If $n\sspeq p^k$ for odd prime $p$ and $ k\sspin \mathbb N$ then \dstyle{q(n)\sspconj  p^{(p^{k-1}+1)/2}}.\psn
{\bf iii)} \ \ if $n\sspeq 2^{k_2}\,p(i_1)^{k_{i_1}}\, \cdots\,p(i_N)^{k_{i_N}}$ with $k_2\sspin \mathbb N_0$, the $i_j-th$  odd primes $p(i_j)$, where  $2\sspleq  i_1\sspkl i_2\sspkl ...\sspkl i_N$, with $N\sspin \mathbb N$ if $k_2\sspneq 0$ and $N\sspgeq 2$ if $k_2\sspeq 0$, then\psn
\Beqarray
q(n)&\sspconj& \prod_{j=1}^N\, p(i_j)^{2^{k_2-1}\, p(i_1)^{k_{i_1}-1}\,\cdots\, p(i_N)^{k_{i_N}-1}\, P(N,j)}\ \ , {\text {with}}\ \ P(N,j)\sspeq \prod_{l=1, l\neq j}^{N}\, (p(i_l)-1)\ ,  \nonumber\\
  &\sspeq&  \prod_{odd\  p|n}\, p^{\frac{\delta(n)}{p-1}}\ . 
\Eeqarray
The last eq. follows from the degree \dstyle{\delta(n)\sspeq \frac{\varphi(2\,n)}{2}}.
This last  formula does, however, not work in the cases {\bf i)} and {\bf ii)}. One can compare this formula with the proven one for the discriminant of the cyclotomic polynomials (the minimal polynomials of \dstyle{exp(2\,\pi\,\frac{1}{n})} (or any of the primitive $n$-th root of $1$), as given in \cite{Ribenboim}, eq. $(1)$ p. 297. For this (slightly rewritten) formula see also \seqnum{A004124} and \seqnum{A193679}. \psn
{\bf Example 7:}  \psn
{\bf ii)}\ \ $p\sspeq p(4)\sspeq 7$,\  $k\sspeq 3$: $q(7^3)\sspeq q(343)\sspeq 7^{(7^2+1)/2}\sspeq 1341068619663964900807$.  \psn
{\bf iii)}\ \ $n\sspeq 2^3\cdot 3^2\cdot 7$, $q(n)\sspeq 3^{2^2\cdot 3\cdot 6}\,7^{2^2\cdot 3\cdot 2}\sspeq  4316018525852839090954658176626149564980915348463203041$. \psn
These values have been checked with the help of Maple13 \cite{Maple}.\psn
{\bf Note added:} The conversion of the $C$ polynomials in terms of Chebyshev $S$ polynomials (\seqnum{A049310}) has been initiated by {\sl Ahmet Zahid K\"u\c c\"uk} and was finished together with the author. See \seqnum{A255237} for the details.\pbn  
%%%%%%%%%%%%%%%%%%%%%%%%%%%%%%%%%%%%%%%%%
%%%%%%%%%%%%%%%%%%%%%%%%%%%%%%%%%%%%%%%%%
\section{Splitting field $\bf \boldsymbol{\mathbb Q}({\boldsymbol \rho}(n))$ for $\bf C(n,x)$, field extension  and Galois group}
The algebraic number \dstyle{\rho(n)\sspeq 2\,\cos\left( {\frac{\pi}{n}}\right)}, $n\sspin \mathbb N$, with minimal polynomial $C(n,x)$ over $\mathbb Q$ of degree $\delta(n)$, has been studied in {sect. 3}. Each of these polynomials, being minimal, is irreducible. All roots have been given in eq. $(19)$ (or eq. $(20)$).  $C(n,x)$ is also separable because all of its roots are distinct. Because \dstyle{  2\,\cos\left( {\frac{\pi\,k}{n}}\right) \sspeq {\hat t}(k,\rho(n))} (see eq. $(2)$, and the coefficient array \seqnum{A127672}), with the monic integer $\hat t$-polynomials, each zero can be written as integer linear combination in the vector space basis $<1,\rho(n),\rho(n)^2,...,\rho(n)^{\delta(n)-1}>$, called the power basis. One has to reduce in $\hat t(k,\rho(n))$ all powers $\rho(n)^p$, $p\sspgeq \delta(n)$ with the help of the equation $C(n,\rho(n))\sspeq 0$. See {\it Table 4} for the zeros of $C(n,x)$, $n\sspeq 1,...,30$ written in this power basis.\pbn
{\bf Example 8: \ $\bf n\sspeq 8$  (octogon)} \psn
The $\delta(8)\sspeq 4$ zeros of $C(8,x) \sspeq x^4\sspm 4\,x^2\sspp 2$ are, with $\rho\sspequiv \rho(8)\sspeq \sqrt{2\sspp \sqrt{2}}$,  $\pm\, \rho$ and $\pm(\sspm 3\,\rho \sspp \rho^3)\sspeq \pm \sqrt{2-\sqrt{2}}$.  In this case the degree $4$ coincides with the number of  DSRs in the upper half plane and the negative real axis.
\pbn
This shows that the extension of the rational field $\mathbb Q$, called $\mathbb Q(\rho(n))$, obtained by adjoining just one algebraic element (called a simple field extension) is the splitting field (Zerf\"allungsk\"orper in German) for the polynomial $C(n,x)$. Note that even though the polynomial $C(n,x)$ is from the ring $\mathbb Z[x]$ one needs $\mathbb Q(\rho(n))$ with rational coefficients $r_j$ for the general element $\alpha\sspeq \sum_{j=0}^{\delta(n)-1}\, r_j\,\rho(n)^j$. For example, \dstyle{{\rho(8)^{-1}\sspeq 2\rho(8) \sspm {\frac{1}{2}} \rho(8)^3}}. Some references for field extensions and {\sl Galois} Theory are \cite{Cox}, \cite{SLang}, chpts. V and VI, \cite{Artin}, and the on-line lecture notes \cite{Lazarev}.  The dimension of of $\mathbb Q(\rho(n))$ as a vector space over $\mathbb Q$ is $\delta(n)$. This is the degree of the extension, denoted usually  by $[\mathbb Q(\rho(n))\,:\,\mathbb Q]$, and it coincides with the degree of the minimal polynomial for $\rho(n)$. Of course, it is a proper extension only if $\delta(n)\sspgeq 2$, \ie for  $n\sspgeq 4$. This extension of $\mathbb Q$ is separable, \ie the minimal polynomial for the general algebraic number $\alpha$ given above is separable (has only distinct zeros). It is a normal field extension, meaning that every irreducible rational polynomial with one root in $\mathbb Q(\rho(n))$ splits completely over $\mathbb Q(\rho(n))$. See, e.g., \cite{Cox}, Theorem 5.24, p. 108. \psn
We now consider a subgroup of the group of automorphisms of $\mathbb Q(\rho(n))$, called ${\cal A}{\sl ut}(\mathbb Q(\rho(n))$, which consists of those elements $\sigma$ which leave the subfield $\mathbb Q$ pointwise invariant (fixed point field $\mathbb Q$): \pn
$\sigma\, :\, \mathbb Q(\rho(n))\sspto \mathbb Q(\rho(n))$, $\alpha\sspmapsto \sigma(\alpha)$, with $\sigma(\beta)\sspeq \beta$ for all $\beta\sspin \mathbb Q$.
The subgroup of ${\cal A}{\sl ut}(\mathbb Q(\rho(n))$ of these so-called $\mathbb Q$-automorphisms is called the {\sl Galois} group of  $\mathbb Q(\rho(n))$ over $\mathbb Q$, denoted by ${\cal G}{\sl al}(\mathbb Q(\rho(n))/\mathbb Q )$. Occasionally we abbreviate this with ${\cal G}_n$.  In order to find the elements $\sigma$ (we omit the label $n$) of this subgroup it is sufficient to know $\sigma(\rho(n))$, because of the usual rules for automorphisms: {\it  (i)}\  $\sigma(\alpha \sspp \beta)\sspeq \sigma(\alpha) \sspp \sigma(\beta)$, {\it (ii)}\ $\sigma(\alpha\, \beta)\sspeq \sigma(\alpha)\, \sigma(\beta)$, and {\it (iii)}\ $\alpha\neq 0 \sspfollows \sigma(\alpha)\neq 0$. Indeed, because of $\sigma(1)\sspeq \sigma(1^2)\sspeq \sigma(1)\,\sigma(1)$, one has $\sigma(1)\sspeq 1$, and the images of products of $\rho(n)$ are obtained from products of $\sigma(\rho(n))$. Applying $\sigma $ on the equation $C(n,\rho(n))\sspeq 0$  (minimal polynomial), leads to $C(n,\sigma(\rho(n)))\sspeq 0$, because the integer (rational) coefficients and $0$ are invariant under  the $\mathbb Q$-automorphism $\sigma$ we are looking for. Therefore we have exactly  $\delta(n)$ distinct $\mathbb Q$-automorphisms $\sigma_j$, $j\sspeq 0,\,...,\,\delta(n)-1$,  determined from the distinct roots of $C(n,\sigma(\rho(n)))$, \viz\   $\sigma_j(\rho(n))\sspeq \tilde \xi^{(n)}_{j+1}$ with the zeros of $C(n,x)$ ordered like in eq. $(19)$ with increasing $k$ values (see the {\sl Table 4}). By the same token $C(n,\sigma(\tilde \xi ^{(n)}_j))\sspeq 0$, for $j=1,...,\delta(n)$ for every $\sigma$.  Because all roots are distinct (separable $C$) this leads to an isomorphism between ${\cal G}{\sl al}(\mathbb Q(\rho(n))/\mathbb Q )$ and a subgroup of the symmetric group $S_{\delta(n)}$. The $\mathbb Q$-automorphisms $\sigma$ can therefore be identified with permutations of the roots of $C$ (see \eg \cite{Cox} ch. 6.3, pp. 132 ff. One identifies the roots  $\tilde \xi^{(n)}_{j}$ with $j$, and  $\sigma_j$  with the permutation \dstyle{\left( {\ 1\ \ \quad 2\quad ...\quad \delta(n)}\atop{\sigma_j(1)\,\sigma_j(2)\,...\,\sigma_j(\delta(n))}\right) \sspequiv [\sigma_j(1)\,\sigma_j(2)\,...\,\sigma_j(\delta(n))]}.  This subgroup of $S_{\delta(n)}$ is Abelian because only products of  powers of $\rho(n)$ appear and due to the automorphism property of $\sigma$ this carries over to the {\sl Galois} group.\pbn
{\bf Example 9: \ $\bf n\sspeq 5$  (pentagon) $\boldsymbol{\mathbb Q}$-automorphisms} \psn      
$\alpha\sspeq r_0\, 1\sspp r_1\, \rho(5)$ with $\rho(5)\sspeq \varphi$, the golden section. In this case $\mathbb Q(\varphi)$ is as quadratic number field usually called $\mathbb Q(\sqrt{5})$ with the basis $<1,\varphi> $ for integers in $\mathbb Q(\sqrt{5})$ (see \eg \cite{HardyWright}, ch. 14.3, p. 207, where $\tau\sspeq (\varphi\sspm 1)$). The two ($\delta(5)=2 $)  $\mathbb Q$-automorphisms are obtained from the solutions $\sigma(\varphi)\sspeq \varphi$ or \dstyle{-{\frac{1}{\varphi}}\sspeq 1\sspm \varphi} of $C(5,n)\sspeq x^2\sspm x\sspm 1$. Hence $\sigma_0\sspeq id:\  \sigma_0(1)\sspeq 1,\ \sigma_0(\varphi)=\varphi $ and $\sigma_1:\  \sigma_1(1)\sspeq 1,\ \sigma_1(\varphi)=1\sspm \varphi$. Because $\sigma_1^2\sspeq \sigma_0$, ${\cal G}_5$ is generated by $\sigma_1$, hence the {\sl Galois} group ${\cal G}{\sl al}(\mathbb Q(\varphi)/\mathbb Q )$ is the cyclic group $Z_2$ (also known as additive group $\mathbb Z/2\,\mathbb Z$ or $C_2$, but we reserve $C$ for the minimal polynomials). The fixed field for ${\cal G}_5$ is $\mathbb Q$; for the trivial subgroup with element $\sigma_0$ it is $\mathbb Q(\varphi)$.\pbn
{\bf Example 10: \ $\bf n\sspeq 7$  (heptagon) $\boldsymbol{\mathbb Q}$-automorphisms} \psn    
The three zeros ($\delta(7)=3$) of $C(7,n)$ are $\tilde \xi^{(7)}_1\sspeq \rho(7)$,  $\tilde \xi^{(7)}_2\sspeq -1\sspm \rho(7)\sspp \rho(7)^2$ and $\tilde \xi^{(7)}_3\sspeq 2\sspm \rho(7)^2$ (see {\it Table 4}). In the sequel we omit the argument $7$. Computing powers of $\rho$ modulo $C(7,\rho)\sspeq 0$ one finds the $\mathbb Q$-automorphisms $\sigma_0\sspeq id$, $\sigma_1(\rho\sspeq \tilde \xi_1)\sspeq \tilde \xi_2$, $\sigma_1(\tilde \xi_2)\sspeq \tilde \xi_3$,  $\sigma_1(\tilde \xi_3)\sspeq \tilde \xi_1$ and $\sigma_2(\rho\sspeq \tilde \xi_1)\sspeq \tilde \xi_3$, $\sigma_2(\tilde \xi_2)\sspeq \tilde \xi_1$,  $\sigma_2(\tilde \xi_3)\sspeq \tilde \xi_2$. Therefore, the identification with $S_3$ permutations is $\sigma_0\sspdoteq [1\,2\,3]\sspeq e$, $\sigma_1\sspdoteq [2\,3\,1]$ and $\sigma_2\sspdoteq [3\,1\,2]$. Each $\sigma_j$ permutation can be depicted in a circle diagram  with vertices labeled $1,2$ and $3$ and directed edges, also allowing for loops. The only (Abelian) group of order three is the cyclic $Z_3$ subgroup of $S_3$.\pbn
We can also characterize these {\sl Galois} groups ${\cal G}_n$ by giving their cycle structure and depict them as cycle graphs. For cycle graphs see, \eg \cite{Wikipedia} ``Cycle graph'' and ``List of small groups'' with all cycle graphs for groups (also non-Abelian ones) of order $1,\,...\, ,16$.  In order to manage powers of {\sl  Galois} group elements we first need some fundamental identities of {\sl Chebyshev} $T$-polynomials, whence $\hat t$-polynomials.\psn
We will need iterations of $\hat t$ polynomials governed by the following well known identity.\psn
{\bf Lemma 7:\  Iteration of  $\bf \hat t$-polynomials}\psn
\Beq
\hat t(n,\hat t(m,x)) \sspeq \hat t(n\,m,x), \ n,m \sspin \mathbb N_0\ .
\Eeq
{\bf Proof:} See \cite{Rivlin}, Exercise 1.1.6, p. 5, first with the trigonometric definition of {\sl Chebyshev} T-polynomials which then carries over to the general polynomials defined by their recurrence relation. This identity is then rewritten for $\hat t$-polynomials.\hskip 11.7cm $\square$\psn
{\bf Lemma 8:\  mod $\bf n$  reduction of $\bf \hat t$-polynomials in the variable $\bf \boldsymbol \rho (n)$}\psn
\Beq
\hat t(k,\rho(n))\sspeq (-1)^{\floor{\frac{k}{n}}}\,\hat t(k(mod\,n),\rho(n)),\  n\in \mathbb N, k\in \mathbb Z\ .
\Eeq    
$(-1)^{\floor{\frac{k}{n}}}\sspfed  p_n(k)$, the parity of ${\floor{\frac{k}{n}}}$, will become important in the following. Sometimes  ${\floor{\frac{k}{n}}}$ is called quotient and denoted also by $k\backslash n$.\psn
{\bf Proof:} This follows, with $k\sspeq l\,n\sspp r$, trivially from the trigonometric identity \dstyle{\cos\left( {\frac{\pi}{n}}\,(l\,n\sspp r)\right)\sspeq \cos\left(\pi\,l\sspp {\frac{\pi\, r}{n}}\right)\sspeq (-1)^l\, \cos\left({\frac{\pi}{n}}\,r\right) }, with $r\in \{0,1,...,n-1\}$. \hskip 7cm $\square$\psn
To simplify notation we will use also $\hat t_n(x)$ for $\hat t(n,x)$. From {\sl lemmata} $7$ and $8$ we have \eg $\hat t_3(t_3(\rho(7))\sspeq \hat t_9(\rho(7))\sspeq -\hat t_2(\rho(7)) \sspeq -(\rho(7)^2\sspm 2)$. This is also $+\hat t_5(\rho(7))$ from the identity  \dstyle{\cos\left( {\frac{\pi}{n}}(n\ssppm l)\right)\sspeq  -\cos\left( {\frac{\pi}{n}}\,l\right)}. This proves the following {\sl  lemma}.\psn
{\bf Lemma 9: Symmetry relation of $\bf \hat t$ polynomials}\psn
\Beq
\hat t(n-l,\rho(n))\sspeq -\hat t(l,\rho(n))\, ,\  n\sspin \mathbb N, \, l\sspin \{0,1, ... ,n\} \ .
\Eeq
{\sl Lemma 7} which used the trigonometric $\rho(n)$ definition can be rewritten 
as a congruence for the $\hat t$-polynomials with indeterminate $x$. This is because all what is needed is that  $\rho(n)$ is a zero of $C(n,x)$. \psn
{\bf Corollary 5:\ Congruence for $\bf \hat t$ polynomials modulo $C$-polynomials}\psn
\Beq 
\hat t(k,x)\sspequiv (-1)^{\floor{\frac{k}{n}}}\, \hat t(k(mod\, n),x)\,(mod\, C(n,x))\ ,\ n\sspin \mathbb N, k\sspin \mathbb Z\ .
\Eeq
An example will illustrate this before we give an another proof of this {\sl corollary} based on known $T$-polynomial identities.\psn
{\bf Example 11:\ Congruence for  $\bf n\sspeq 7$, $\bf k\sspeq 9$}\psn
$\hat t(9,x)\sspeq x^9\sspm 9\,x^7\sspp 27\,x^5\sspm 30\,x^3\sspp 9\,x$ and $C(7,x)\sspeq x^3\sspm x^2\sspm 2\,x\sspp 1$. Polynomial division shows that  $\hat t(9,x)\sspeq (x^6\sspp x^5\sspm 6\,x^4\sspm 5\,x^3\sspp 9\,x^2\sspp 5\,x-2)\,C(7,x) \sspp (-x^2\sspp 2) $, hence $\hat t(9,x)\sspequiv  - \hat t(2,x)\,(mod\, C(7,x))$.\psn
The alternative  proof of {\sl corollary 5} is based on the following  factorization of {\sl Chebyshev} $S$-polynomials in terms of the minimal polynomials $\Psi(n,x)$ of \dstyle{\cos\left({\frac{2\,\pi}{n}}\right)} over $\mathbb Q$  (for these polynomials see \seqnum{A049310} and \seqnum{A181875/A181876}).\psn
{\bf Proposition 9: \ Factorization of $\bf S$-polynomials in terms of $\bf \Psi$-polynomials}\psn
\Beq
S(n-1,x)\sspeq 2^{n-1}\,\prod_{2<d\,|\,2n}\,\Psi\left(d,{\frac{x}{2}}\right)\, , \, n\sspin \mathbb N\ .
\Eeq
As usual, the undefined (empty)  product is defined to be $1$.\psn
{\bf Proof:} Start with the eq. $(3)$ of \cite{WatkinsZeitlin}, p. 471, written for $n\sspto 2\,n$ and for the $\hat t$-polynomials instead of the $T$-polynomials. This is \dstyle{\hat t(n+1,x)\sspm \hat t(n-1,x)\sspeq 2^{n+1}\,\prod_{d|2\,n}\,\Psi\left(d,{\frac{x}{2}}\right)}. Then use the identity \cite{MOS}, p. 261, first line,  with $m\sspto n$, $n\sspto m$ and $U$-polynomials replaced by $S-$polynomials. This leads to the identity 
\Beq
\hat t(n+m,x)\sspm \hat t(n-m,x)\sspeq (x^2\sspm 4)\, S(n-1,x)\, S(m-1,x)\,,\ m\sspleq n\sspin \mathbb N\, . 
\Eeq
Here we only need the case $m\sspeq 1$.
Then dividing out $2\,\Psi(1,{\frac{x}{2})}\sspeq x\sspm 2$ and $2\,\Psi(2,{\frac{x}{2})}\sspeq x\sspp 2$ leads to the claimed identity. \hskip 13cm $\square$\psn
To end the preparation for an independent proof of {\sl corollary 5} we state:  \psn
{\bf Corollary 6: $\bf C$-polynomial divides some family of $\bf S$-polynomials}\psn
\Beq
C(n,x)\ |\, S(l\,n-1,x)\, , \, n\sspgeq 2,\  l\sspin \mathbb N\ .
\Eeq
This is clear from the definition of $C$ in eq. $(16)$ and the fact that $2\sspkl 2\,n\,|\, 2\,l\,n$, for $l\sspgeq 1$. \Eg $C(10,x)\sspeq x^4\sspm 5\,x^2\sspp 5 $ divides the family $\{S(9,x),S(19,x), S(29,x),...\}$\psn
{\bf Proof of corollary 5:}  From the {\sl corollary 6} and the identity eq. $ (63)$  with $n\sspto l\,n$ and $m\sspto r$ one sees that $\hat t(k,x)\sspeq \hat t(l\,n\sspp r,x)\sspequiv \hat t(l\,n\sspm r,x)\,(mod\, C(n,x))$. Now $\hat t(l\,n\sspm r,x)\sspeq \hat t((l-1)\,n\sspp (n-r),x)$, and one can use $l_1\sspdef l-1$ and $r_1\sspdef n-r$ (remember that $r\sspin \{0,1,...,n-1\}$, hence $r_1$ has identical range) as new $l$ and $r$ variables in this congruence, to get $\hat t(l_1\,n\sspp r_1,x)\sspequiv \hat t(l_1\,n\sspm r_1,x)\, (mod\, C(n,x))\sspeq $\pn
$\hat t((l-2)\,n\sspp r,x)\, (mod\, C(n,x))$. This can be continued until one ends up with \pn
 $\hat t(0\,n\sspp k(mod\, n),x)\, (mod\, C(n,x))$.\hskip 11.5cm $\square$\psn
Now we are in a position to compute powers of elements of the {\sl Galois} group ${\cal G}_n$. The subset of odd numbers $2\,l\sspp 1\sspkl n$ entering the product in eq. $(20)$ will be denoted by ${\cal M}(n)$. There are $\delta(n)$ (degree of $C(n,x)$) such odd numbers.\psn
{\bf Definition 2:\  The fundamental set ${\mathbfcal {\cal M}(n)}$}
\Beq 
{\cal M}(n)\sspdef \left\{2\,l\sspp 1 {\biggm \vert } l \sspin \left\{0,....,\floor{\frac{n-2}{2}}\right\} \text{and}\, gcd(2\,l + 1,n) = 1 \right\}\sspeq\{m_1(n), ...,m_{\delta(n)}(n)\}, 
\Eeq
with $m_1(n)\sspeq 1$ and we use the order  $m_i(n)\sspkl m_j(n)$ if $i\sspkl j$. For $n\sspeq 1$ one takes ${\cal M}(1)\sspeq \{1\}$.\psn
{\bf Example 12:}  \dstyle{|{\cal M}(2)|\sspeq \delta(2)\sspeq 1,\,  {\cal M}(2) \sspeq \{1\}};\ \dstyle{|{\cal M}(14)|\sspeq \delta(14)\sspeq 6\, ,\, {\cal M}(14) \sspeq \{1,3,5,9,11,13\} }.\psn
For ${\cal M}(n)$ see the row No. $n$ of the  array \seqnum{A216319}.\psn
For later purpose we define here, for odd $n$, the extended fundamental set $\widehat {\cal M}(n)$ and its first difference set $\triangle\,\widehat {\cal M}(n)$
\psn
{\bf Definition 3:\  The fundamental extended set $\widehat {\mathbfcal {\cal M}}(\bf n)$ for odd $\bf n$}\psn
For $n$ odd, $\sspgeq 1$:\ \dstyle{\widehat {\cal M} (n)\sspdef \{0,m_1(n)=1,...,m_{\delta(n)}(n)\sspeq n-2, n+2\}}\ .\psn
Thus $|\widehat {\cal M} (n)|\sspeq \delta(n)+2$. Note that $gcd(n\sppm 2,n)\sspeq 1$ for odd $n$ (proof by assuming the contrary:  $gcd(n+2,n)\sspeq d>2$ because $n$ and $n+2$ are odd. Then $d\,|\,(n+2)$ and $d\,|\, n$, hence $d\,|\, ((n+2)-n)$, $d\,|\,2$, implying $d=2$ or $d=1$, but $d>2$.) The reason for defining this extended set is that for odd $n$ the first difference set  \dstyle{\triangle\,\widehat {\cal M} (n)\sspeq \{1,\triangle m_2(n),...,\triangle m_{\delta(n)}(n), 4\}} with $\triangle m_j(n)\sspdef m_j(n)\sspm m_{j-1}(n)$, will become important later on.  $\triangle m_1(n)\sspeq 1$ and $\triangle m_{\delta(n)+1}(n) \sspeq 4$ for each odd $n$.\psn
With {\it definition 2} eq. $(20)$ implies, 
\Beqarray
\sigma_j(\rho(n))&\sspeq& \sigma_j(\tilde\xi_1^{(n)}) \sspeq \tilde \xi_{j+1}^{(n)}\sspeq 2\,\cos\left(\frac{\pi}{n}\,m_{j+1}(n) \right) \sspeq \nonumber \\
 && {\hat t}(m_{j+1}(n),\rho(n))\  \text{for}\, \ j\sspin \{0,1,...,\delta(n)-1\}\ .
\Eeqarray
{\Eg} $n=7$: $\delta(7)\sspeq 3$, ${\cal M}(7)\sspeq \{1,3,5\}$, $\sigma_0(\rho(7))\sspeq \rho(7)$, $\sigma_1(\rho(7))\sspeq {\hat t}(3,\rho(7))$, and $\sigma_2(\rho(7))\sspeq {\hat t}(5,\rho(7))$.\psn
Because $\hat t$ is a rational integer polynomial 
\Beq
\sigma_j^2(\rho(n))\sspeq \sigma_j(\sigma_j(\rho(n)))\sspeq {\hat t}(m_{j+1}(n),\sigma_j(\rho(n)))\sspeq {\hat t}(m_{j+1}(n),{\hat t}(m_{j+1},\rho(n))) \ ,
\Eeq
and with {\sl lemma 6} this becomes \dstyle{\sigma_j^2(\rho(n))\sspeq {\hat t}((m_{j+1}(n))^2,\rho(n)) }. In general we have
\Beq
\sigma_j^k(\rho(n))\sspeq {\hat t}((m_{j+1}(n))^k,\rho(n))\,,\, \text{for}\, j\sspin \{0,1,...,\delta(n)-1\}\ , n,k\sspin \mathbb N\ .
\Eeq
Instead of powers of $\sigma_j(\rho(n))$ we can therefore consider powers of $m_{j+1}(n)$. {\it Lemmata} $8$ and $9$ are now employed to prove that a $\hat t$ polynomial with a product of elements from ${\cal M}(n)$ as its first argument (or index) is  again a $\hat t$ polynomial with first argument from ${\cal M}(n)$. In this way one can build sequences of powers, starting from any element of ${\cal M}(n)$. Trivially, $1^k\sspeq 1$. Before proving this closure of  ${\cal M}(n)$ under powers, provided the rules for $\hat t$ polynomials are taken into account, we give two examples, and then  define a new equivalence relation on the integers, called $\Moddn{n}$, denoted by  $\Simn$.\psn
{\bf Example 13: Cycle structure for $\bf n=12$ (dodecagon)}\psn
\dstyle{\delta(12)\sspeq 4,\ {\cal M}(12)\sspeq \{1,5,7,11\}}.\ $5^2\sspeq 25 
\equiv 1\, (mod\ 12)$, reflecting  {\it lemma 8}. The  sign $p_n(25)$  in eq. $(59)$ is here $+$. (If later the sign will be $-$ , the $mod\ n$ result will be underlined.) The first $2$-cycle is therefore $[5,1]$. Similarly, $7^2\sspeq 49\equiv 1\, (mod\ 12)$ (sign $+$), whence the second $2$-cycle is $[7,1]$, and finally, $11^2\sspeq 121\equiv 1\, (mod\ 12)$ (sign $+$), producing the third  $2$-cycle $[11,1]$.   This result appears as the $n\sspeq 12$ entry in {\it Table 6}. In this example the {\sl Galois} group is not generated by one element. hence it is non-cyclic. In fact,  ${\cal G}{\sl al}(\mathbb Q(\rho(12))/\mathbb Q )\sspeq Z_2\times Z_2\sspeq Z_2^2$ (see first entry in {\it Table 8}). The corresponding cycle graph is shown in {\it Figure 4} as the first entry, where the  shaded (colored) vertex stands for $1$ and the open vertices should here be labeled with $5, 7$ and $11$. The cycle structure is $2_3$ (three $2$-cycles). In this example it was not necessary to employ {\it lemma 9} because the signs were always $+$.\psn
{\bf Example 14: Cycle structure for $\bf n=7$ (heptagon)}\psn 
\dstyle{\delta(7)\sspeq 3,\ {\cal M}(7)\sspeq \{1,3,5\}}. $3^2\equiv \underline{2}\, (mod\ 7)$, where now the sign $-$ in eq. $(59)$ is remembered by the underlining. $2\sspnotin {\cal M}(7)$, and now {\it lemma 9} is used to rewrite this $\underline{2}$ as $7\sspm 2\sspeq 5$. Therefore, $3^2\simn1{7} 5\sspin {\cal M}(7)$ (or  $3^2\equiv  5\, \Moddn{7}$). The symbol $\simn1{7}$ (or $Modd\ 7$) is used for the congruence in the new sense, due to {\it lemmata} $8$ and $9$.  Then $3\cdot 5\sspeq 15\equiv 1\, (mod\ n)$ (sign $+$). The first $3$-cycle is therefore $[3,5,1]$. Here $3$ generates all the elements of  ${\cal M}(7)$, and  ${\cal G}{\sl al}(\mathbb Q(\rho(7))/\mathbb Q )\sspeq Z_3$, the cyclic group of order $3\sspeq \delta(7)$. The corresponding cycle graph is a circle with three vertices, one of them, labeled $1$, is shaded (colored) and the other two open  ones are labeled by $3$ and $5$.\psn
This brings us to the definition of an equivalence relation $\Simn$ (or \Modd\,$n$) over the integers $\mathbb Z$. Remember that the floor function for negative arguments is defined as $\floor{-x} \sspeq -\floor{x}$ if $x\sspin \mathbb N_0$, and $\floor{-x} \sspeq -(\floor{x}+1)$ if $0\sspkl x\sspnotin \mathbb N$.\psn
{\bf Definition 4: Equivalence relation $\bf \Simn$ on $\bf \mathbb Z$}\psn
For $k, l \sspin \mathbb Z, n\sspin \mathbb N $:\ $k\Simn l  \Aequ a_n(k)\sspeq a_n(l) $, with the map $a_n\,:\, \mathbb Z\to I_n\sspdef\{0,1,...,n-1\} $, $k\mapsto a_n(k)$, where
\Beq
a_n(k) \sspeq {\Caseszwei{r_n(k)}{{\text if}\ $p_n(k)\sspeq +1\ ,$}{r_n(-k)}{{\text if }$ p_n(k)\sspeq -1 \ ,$}}
\Eeq   
where we used the division algorithm to write $k\sspeq q_n(k)\, n \sspp r_n(k)$, with the quotient $q_n(k)\sspin \mathbb Z$ and the residue $r_n(k)\sspin I_n$. Note that \dstyle{q_n(k)\sspeq \floor{\frac{k}{n}}}. The sign \dstyle{p_n(k)\sspeq  (-1)^{\floor{\frac{k}{n}}} \sspeq (-1)^{q_n(k)}}  corresponding to the parity of $ q_n(k)$, appeared already in {\it lemma 8}. $a_n(1)\sspeq 0$ because $-1\sspequiv 0\,(mod\,1)$. \psn 
Instead of $k\Simn l $ we also write $k\sspequiv l\, \Moddn{n}$ (this should not to be confused with $mod\, n$).  Therefore the sequence  $a_n$ could also be called \Modd\, $n$.  The first of these $2\,n$-periodic sequences $a_n$ are found in \seqnum{A000007}$(n+1), n\sspgeq 0$, (the 0-sequence), \seqnum{ A000035}, \seqnum{A193680},  \seqnum{A193682}, \seqnum{A203571},  \seqnum{A203572} and  \seqnum{A204453}, for $n=1,...,7$, respectively. \pn
The smallest non-negative residue system  $mod\ n$, \viz \, $0,1,...,n-1$ is used here. For  the residue of $k$  modulo $n$  Maple \cite{Maple}  uses $r_n(k)\speq modp(k,n)$.  We also use $k(mod\,n)$ for $r_n(k)$. The reader should verify that $\Simn$ is indeed an equivalence relation satisfying reflexivity, symmetry and transitivity. The disjoint and exhaustive equivalence classes are given by $\{[0], [1], .., [n-1]\}$, called the smallest non-negative complete representative classes (or residue classes) \Modd\,$n$ (we omit the index $n$ at the classes: $[m]\sspeq _n[m]$, written this way in order to distinguish this class from the ordinary one $[m]_n$ used in the arithmetic $mod\, n$). These classes are defined by  $[m]\sspeq \{ l \in \mathbb Z\,\vert l\, \Simn m \}\sspeq \{ l \in \mathbb Z\,\vert \, a_n(l)=m\}$. Because $r_n(-k) \sspeq 0$ if $k\sspequiv 0\,(mod\ n)$ and  $r_n(-k) \sspeq  n \sspm r_n(k)$ if $k\sspnotequiv 0\,(mod\ n)$ (later listed as {\it lemma 17}) one can characterize these residue classes also in the following way. \psn
{\bf Lemma 10:  Complete residue classes $\bf Modd\, n$}\psn
\Beqarray
 l&\in& [0] \Aequiv l\equiv 0\, (mod\, n).\  \text{Equivalently,\ } [0]\sspeq [0]_n\ ,\nonumber \\
l&\in& [1]  \Aequiv l\equiv 1\, (mod\, 2\,n) \ \text{or}\  \equiv -1\, (mod\, 2\,n) .\  \text{Equivalently,\ } [1]\sspeq [1]_{2\,n}\sspunion  [2\,n-1]_{2\,n}\ , \nonumber \\
&\vdots & \nonumber\\
l&\in& [n-1]  \Aequiv l\equiv n-1\, (mod\, 2\,n) \ \text{or}\  \equiv -(n-1)\, (mod\, 2\,n) . \nonumber\\  
&& \text{Equivalently,\ } [n-1]\sspeq [n-1]_{2\,n}\sspunion  [n+1]_{2\,n}\ .
\Eeqarray
For example, take $n\sspeq 7$, then  $[3]\sspeq \rlap{$_7$}\ [3]\sspeq [3]_{14}\sspunion [14-3]_{14}\sspeq \{...,-25,-11,3,17,...\}\sspunion$\pn
$\{...-17, -3, 11,25, ...\}\sspeq \{ ...,-25,-17,-11,-3,3,11,17,25,... \}$. If $n\sspeq 2$ one has $ \rlap{$_2$}\ [0]\sspeq [0]_2$ and  $ \rlap{$_2$}\ [1]\sspeq [1]_2$. The first differences in the class $[0]$ are, of course, $n$, and in the class $[m]$, for $m=1,2,...,n-1$, they alternate between $\triangle_1\sspeq 2\,(n\sspm m)$ and $\triangle_2\sspeq 2\, m$, \eg $n\sspeq 7, m\sspeq 3,\  \triangle_1\sspeq 8,\  \triangle_2\sspeq 6:\  3,\, 11,\,\, 17,\, 25,\, 31,\,...$. This difference alternation invalidates certain theorems known for $mod\, n$. \Eg {\it Theorem} $53$ of reference \cite{HardyWright}, p. 50, is no longer true: take $m\sspeq 5,\, n\sspeq 7,\ a\sspeq 3,\  b\sspeq 17$.  $a$ and $b$ belong both to the class $ \rlap{$_5$}\ [3]$ as well as $ \rlap{$_7$}\ [3]$ but they do obviously not both belong to the class  $_{35}[3]$ .  The following  {\it lemma 11} shows that is is sufficient to know the positive values, and append the negative of these values for each class $[m]$, for $m\sspin\{0,1,...,n-1\}$.\psn 
{\bf Lemma 11: Antisymmetry of the classes $\bf [m]$, for $\bf m\sspin \{0,...,n-1\}$}\psn
For $n\sspin \mathbb N$ and every $m\sspin\{0,...,n-1\}$ the elements of the equivalence class [m] are antisymmetric around $0$.\psn 
{\bf Proof:} This is obvious for the class $[0]\sspeq \{...-2\,n,-n,0,n,2\,n,...\}$. For $m>0$ the negative of every positive numbers of $m\,(mod\,2\,n)$ appears as a negative number of $(2\,n-m)\,(mod\,2\,n)$, and {\it vice versa}: the negative of every  positive numbers of $[2\,n-m]_{2\,n}$ appears as a negative one of $[m]_{2\,n}$. This is true because $-(m+l\,2\,n)\sspeq (2\,n-m)\sspm (l+1)\,2\,n$, for every $l\sspin \mathbb N$, and similarly, $-((2\,n-m)\sspp l\,2\,n)\sspeq m\sspm (l+1)\,2\,n$,  for every $l\sspin \mathbb N$.\hskip15cm $\square$\psn 
This leads immediately to the following {\it corollary}.\psn
{\bf Corollary 7: Non-negative elements of the classes $\bf [m]$}\psn
For $n\sspin \mathbb N$ and every $m\sspin\{0,...,n-1\}$ one has
\Beq
[m]_{\geq} \sspeq [m]_{2\,n,\geq} \sspunion [2\,n-m]_{2\,n,>}\, ,\ \text{and}\  [m]\sspeq [m]_{\geq} \sspunion -([m]_{>}) \ .
\Eeq
Here we used the notations $[m]_{\geq}$  and  $[\,\cdot\,]_{2\,n,\geq}$ or $[\,\cdot\,]_{2\,n,>}$ to denote the subset of non-negative numbers of [m] and the non-negative or positive numbers of the ordinary residue classes $mod\,2n$, respectively. Of course $0$ appears only in the class $[0]$. $ -([m]_{>})$ is obtained from the set $[m]_{>} $(excluding $0$ in the case of class $[0]$) by taking all elements negative. Note that $[-m]$ is not used here. \psn
With $g(n,m)\sspdef gcd(m,2\,n-m)$ one has, for $m\sspin \{1,2,...,n-1\}$, $[m]_{2n,>}\sspeq  g(n,m)\,[m/g(n,m)]_{2\,n/g(n,m),>}$ and $[2\,n-m]_{2n,>}\sspeq g(n,m)\,[(2\,n-m )/g(n,m)]_{2\,n/g(n,m),>}$. This is. of course, only of interest if $g(n,m)\sspneq 1$. \Eg $n=6,\ m=3, \ g(n,m)=3$: $[3]_{12,>} = 3*[1]_{4,>}$ and $[9]_{12,>}\sspeq 3*[3]_{4,>}$, where again $k*[p]_{q,>}$ is the set with all members of the set $[p]_{q,>}$ multiplied with  $k$. This is obvious.\psn  
The trivial formula for the  members of the residue classes $[m]_{\geq}$, considered as sequences of increasing numbers called  $\{c(n,m;k)\}_{k=1}^{\infty}$, is \psn
\Beq
c(n,m;k)\sspeq {\Caseszwei{(k-1)\,n}{{\text if}\ $m\sspeq 0$\ ,}
{\floor{{\frac{k}{2}}}\, 2\, n\sspp (-1)^k\,m}{{\text if}\ $m\sspin \{1,2,...,n-1\}$\ . }}
\Eeq
The nonnegative members of the complete residue classes $\Moddn{n}$ for $n\sspeq 3,4,5,6,$ and $7$ are found in \seqnum{A088520}, \seqnum{A203575}, \seqnum{A090298},  \seqnum{A092260}, and \seqnum{A113807}. Sometimes $0$ has to be added, in order to obtain the class $[0]$. Of course, these complete residue classes can be recorded as a permutation sequence of the non-negative integers.\psn
We now list several {\it lemmata} (the trivial {\it lemma 17} has already been used) in order to prepare for the proof of the multiplicative structure of these equivalence classes.\psn
{\bf Lemma 12: Parity of $\bf Modd\, n$ residue classes}\psn
For even $n$ the parity of the members of the residue class $_n[m]$, $m\sspin \{0,1,...,n-1\}$ coincides with the one of $m$. If $n$ is odd this is also true for the classes with $m\sspin \{1,\, 2,\, ... \,,n-1\}$, and for $m=0$ the parity of the elements alternates, starting with $+$ (for even).\psn
{\bf Proof}: The case of the residue class $[0]$ is clear for even or odd $n$ because of its members $0\, mod\, n$ (see {\sl lemma 10}). Similarly, for the other $m$ values, because then $mod\, 2n$ applies.\psn 
{\bf Lemma 13: Periodicity of the parity sequence $\bf p_n$ with period length $\bf 2\,n$ }\psn
\Beq
p_n(k)\sspeq p_n(k\sspp 2\,n\,l)\ \text{for}\ l\sspin \mathbb Z,\ \text{\ie}\   p_n(k)\sspeq p_n(r_{2\,n}(k))\sspeq p_n(k\,(mod\,2\,n)),\ \text {for}\ k\sspin \mathbb Z, n\sspin \mathbb N.  
\Eeq
{\bf Proof:} This is obvious from the definition of $p_n(k)$ given in {\it lemma 8}, eq. $(59)$, with {\it definition 4}, eq. $(69)$. For the second part use $k\sspeq 2\,n\,q_{2\,n}(k)\sspp r_{2\,n}(k)$. \hskip 13cm $\square$ \psn
{\bf Lemma 14: (A)symmetry of  sequence $\bf p_n$ around $\bf k\sspeq 0$}\psn
\Beq
\text {for}\ n\sspin \mathbb N, \ k\sspin \mathbb N_0\, :\  p_n(-k)\sspeq \Caseszwei{+p_n(k)}{\text {if}\ $k\sspequiv 0\,(mod\, n)$}{-p_n(k)\ ,}{\text {if}\ $k\sspnotequiv 0\,(mod\, n)$\ . }  
\Eeq
{\bf Proof:} This follows immediately  from the property of the $\floor{-x}$ function mentioned above before {\it definition 4}. \hskip 16cm $\square$ \psn
{\bf Lemma 15: Product  formula for the residue $\bf r_n$ }\psn
For $n\sspin \mathbb N$\ \text{and}\ $k,l\sspin \mathbb Z$ one has:  \dstyle{r_{n}(k\,l)\sspeq r_n(k)\,r_n(l)\, (mod\, n)}.
\psn
{\bf Proof:} Just  multiply $k\sspeq q_n(k)\,n\sspp r_n(k)$ with $l\sspeq q_n(l)\,n\sspp r_n(l)$. \hskip 14cm $\square$\psn
{\bf Lemma 16: Residue $\bf r_{2\,n}$ from $\bf r_{n}$ }\psn
 For $n\sspin \mathbb N$\ \text{and}\  $k\sspin \mathbb Z$ one has: \dstyle{ r_{2\,n}(k)\sspeq \Caseszwei{r_n(k)}{\text {iff}\ $r_{2\,n}(k)\sspin \{0,\,1,\,...,\,n-1\}$,}{r_n(k)+n}{\text {iff}\ $r_{2\,n}(k)\sspin \{n,\,n+1,\,...,\,2\,n-1\}$. } } \psn
{\bf Proof:} Obvious for $r_{2\,n}(k)\sspin \{0,...,n-1\}$ as well as  $\sspin \{n,...,2\,n-1\}$.\hskip 5.5cm $\square$\psn
\vfill
\eject
\noindent
{\bf Lemma 17: Residue for negative numbers}\psn
For $n\sspin \mathbb N$\ \text{and}\  $k\sspin \mathbb N_0$ one has: \dstyle{  r_n(-k)\sspeq\Caseszwei{0}{\text {if}\ $k\sspequiv 0\,(mod\,n)$\ ,}{n\sspm r_n(k)}{\text {if}\ $k\sspnotequiv 0\,(mod\, n)$\ . } }
\psn
{\bf Proof:} Obvious. \hskip 14.2cm $\square$\psn
{\bf Lemma 18:  Symmetry of the $\bf a_n$ (or $\bf Modd\, n$) sequence around $\bf k=0$}\psn
For $n\sspin \mathbb N$\ \text{and}\  $k\sspin \mathbb N_0$ one has: $a_n(-k) \sspeq a_n(k)$.\psn
{\bf Proof:} This is proved for the two cases $k\sspequiv 0\,(mod\,n)$ and  $k\sspnotequiv 0\,(mod\,n)$ separately. In the first case $r_n(k)\sspeq 0$, as well as $r_n(-k)\sspeq 0$ , hence $a_n(k\sspequiv 0\,(mod\, n))=0$, which is symmetric around $k=0$. In the other case, we employ {\it lemma 13}, noting that $r_{2\,n}(k)\sspneq 0, n$ (otherwise $k\sspequiv 0\,(mod\,n)$). The two cases $r_{2\,n}(k)\sspin \{1,\,2,\,...,\,n-1 \}$ and $r_{2\,n}(k)\sspin \{n+1,\,n+2,\,...,\,2n-1 \}$ have $p_n(k)\sspeq p_n(r_{2\,n}(k))$ equal $+1$ or $-1$, respectively. Then with the second alternative of {\it lemma 14} one finds $a_n(-k)\sspeq r_{n}(-k)$ if $p_n(-k)\sspeq -p_n(k)\sspeq +1$ and  $a_n(-k)\sspeq r_{n}(k)$ if $p_n(-k)\sspeq -p_n(k)\sspeq -1$  which coincides with the definition of $a_n(+k)$.     \hskip 16cm $\square$\psn
Now we turn to the arithmetic structure of the \Modd\,$n$ residue classes. It is clear from the following counter-example that addition cannot be done class-wise. Consider $n\sspeq 6$, $k\sspeq 2$ and $l\sspeq 7$. Then $a_6(2\sspp 7)\sspeq 3$ but $a_6(a_6(2)\sspp a_6(7))\sspeq a_6(2\sspp 5)\sspeq a_6(7)\sspeq 5$. In other words, $2\simn1{6} 2$ and $7\simn1{6} 5 $, but  $2\sspp 7\sspeq 9\simn1{6} 3$ but $ (2\sspp 5)\sspeq 7\simn1{6} 5$, and $3$ is not equivalent to $5\,\Moddn{6}$. Similarly, from $(a-1)\simn1{n} 0$ does in general not follow  $a\simn1{n} 1$, because for $n>2$ in the latter case $a$ can also be of the form  $-1+k\,2\,n$ if it belongs to a class $[m]$ with positive $m$, whereas in the first case it has to be of the form  $1+k'\,n$, which can only match for $n\sspeq 1$ and $2$.  However, it turns out that multiplication can be done class-wise. This is the content of the following {\it proposition}.\psn 
{\bf Proposition 10: Modd\,$\bf n$ residue classes are multiplicative}\psn
For $n\sspin \mathbb N$ and  $k,l\sspin \mathbb Z$ one has:
\Beq
a_n(k\,l) \sspeq a_n(a_n(k)\,a_n(l))\,\ \text{\ie} \  k\,l\,\Simn\, a_n(k)\,a_n(l)\, ,\, {\text {i.e.,}}\,  kl\sspequiv a_n(k)\,a_n(l)\, ({\text{\Modd}}\,n)\ .
\Eeq 
Before giving the proof, consider the example $n\sspeq6$, $k\sspeq 2$ and $l\sspeq 7$. 
Now  $a_6(2\cdot 7)\sspeq 2$ and $a_6(a_6(2)\cdot a_6(7))\sspeq a_6(2\cdot 5)\sspeq a_6(10)\sspeq 2$. Or stated equivalently, $2\cdot 7\sspeq 14\simn1{6} 2\simn1{6} 2\cdot 5\simn1{6} a_6(2)\ a_6(7)$.\psn
{\bf Proof:}\psn
{\bf i)} Due to the symmetry of $a_n$ (see {\it lemma 18}) it is clear that it is sufficient to consider only non-negative $k$ and $l$.\psn
{\bf ii)} Consider first the cases $k\sspequiv 0\,(mod\, n)$ or $ l\sspequiv 0\,(mod\, n)$. If $m\sspequiv 0\,(mod\, n)$ then $a_n(m)\sspeq 0$.  This follows for both alternatives in eq. $(69)$. Therefore, if $k\sspequiv 0\,(mod\, n)$, $a_n(k\,l) = 0$ because $k\,l \sspequiv 0\,(mod\, n)$ for every $l$, and $a_n(0\cdot a_n(l))\sspeq a_n(0)\sspeq 0$, proving the assertion. In the other case, $l\sspequiv 0\,(mod\, n)$, the proof is done analogously. \psn
{\bf iii)} Now $r_n(k)$ and $r_n(l)$ are non-vanishing, and $k$ and $l$ are positive. Four cases are distinguished according to the signs of $(p_n(k),p_n (l))$, \viz\  $(+,+)$, $(-,-)$, $(+,-)$ and $(-,+)$. \psn
{\bf (+,+)}: In this case, due to the $2\,n$-periodicity of $p_n$ (see {\it lemma 13}), $r_{2\,n}(k)$ and $r_{2\,n}(l)$ are both from $\{1,2,...,n-1\}$, hence 
\Beq
{\bf (+,+)}:\hskip 2cm\  r_{2\,n}(k)\sspeq r_n(k) \ {\text and}\  r_{2\,n}(l)\sspeq r_n(l)\,.
\Eeq 
Also, from eq. $(69)$, one obtains in this case $a_n(k)\,a_n(l)\sspeq r_n(k)\,r_n(l)$. There is an alternative for $a_n(r_n(k)\,r_n(l))$, depending on $p_n(r_n(k)\,r_n(l))$  being $+1$ or $-1$. Both cases are possible as the following examples for $n\sspeq 6$ show: $p_6(3\cdot 4)\sspeq p_6(12\,(mod\, 12))\sspeq p_6(0)\sspeq +1$ and  $p_6(3\cdot 7)\sspeq p_6(21\,(mod\, 12))\sspeq p_6(9)\sspeq -1$. First the argument of $p_n$ is rewritten with the help of eq. $(76)$. $p_n(r_n(k)\,r_n(l))\sspeq p_n(r_{2\,n}(k)\,r_{2\,n}(l))$. Due to the $(2\,n)$-periodicity ({\it lemma 13}) this is $p_n(r_{2\,n}(k)\,r_{2\,n}(l)\, mod(2\,n))\sspeq p_n(r_{2\,n}(k\,l))\sspeq p(k\,l)$, due to {\it lemma 15} with $n\sspto 2\,n$, and again the $(2\,n)$-periodicity. In the first alternative, the $+1$ case, $a_n(a_n(k)\,a_n(l))\sspeq r_n(k)\,r_n(l)\,(mod\, n) \sspeq r_n(k\,l)$, again from {\sl lemma 15}. This is just $a(k\,l)$ if  $p_n(k\,l)\sspeq +1$, proving the assertion for this alternative. In the other case, $p_n(k\,l)\sspeq -1$,  $a_n(r_n(k)\,r_n(l))\sspeq -(r_n(k)\,r_n(l))(mod\, n)$ which is rewritten with $mod\,n$-arithmetic and {\it lemma 15} as $-r_n(k\,l) (mod\, n)\sspeq -(k\,l)(mod\, n)\sspeq r_n(-k\,l)$. This coincides with $a(k\,l)$ for this alternative, proving the assertion. \psn
{\bf (-,-)}:  Now we have from {\it lemma 16}
\Beq
{\bf (-,-)}:\hskip 2cm \  r_{2\,n}(k)\sspeq n\sspp r_n(k) \ {\text and}\  r_{2\,n}(l)\sspeq n\sspp r_n(l)\,.
\Eeq 
Here $p_n((r_n(-k)\,r_n(-l))$ is rewritten with {\it lemma 17}, eq. $(77)$, the  $2\,n$ periodicity, and {\it lemma 15}, as follows.  $p_n((n-r_n(k))\,(n-r_n(l)))\sspeq p_n((2\,n-r_{2\,n}(k))\,(2\,n -r_{2\,n}(l)))\sspeq p_n(r_{2\,n}(k)\,r_{2\,n}(l))\sspeq p_n(r_{2\,n}(r_{2\,n}(k)\,r_{2\,n}(l))) \sspeq p_n(r_{2\,n}(k)\,r_{2\,n}(l)\,(mod\,2\,n)).$ With {\it lemma 15} (with $n\sspto 2\,n$ this becomes $p_n(r_{2\,n}(k\,l))\sspeq p_n(k\,l)$, again from the $(2\,n)$-periodicity. In the first alternative $p_n(k\,l)\sspeq +1$ and  $a_n(r_n(-k)\,r_n(-l))\sspeq r_n(r_n(-k)\,r_n(-l))$.  With {\it lemma 17} and  $mod\, n$-arithmetic this is $r_n((-r_n(k))\,(-r_n(l)))\sspeq r_n(r_n(k)\,r_n(l))$, and with {\it lemma 15} this becomes $r_n(k\,l)$, coinciding with $a_n(k\,l)$ for this alternative. For the other alternative, $p_n(k\,l)\sspeq -1$, $a_n(r_n(-k)\,r_n(-l))\sspeq -(r_n(-k)\, r_n(-l))(mod\, n)$. This becomes, with {\it lemma 17}, $mod\, n$-arithmetic and {\it lemma 15} \pn
$-(r_n(k)\,r_n(l))\,(mod\, n) \sspeq  -r_n(k\,l)(mod\, n)$. This vanishes if $r_n(k\,l)\sspeq 0$ which means  $k\,l\sspequiv 0\,(mod\,n)$ (which  is possible, \eg $n=6, k=2, l=3$), and then this coincides with the claim which is for this alternative $a(k\,l)\sspeq r_n(-(k\,l))\sspeq 0$. If $r_n(k\,l)\sspneq 0$ then $ -r_n(k\,l)\,(mod\, n)\sspeq n\sspm r_n(k\,l)$ which also coincides with the claim  $a(k\,l)\sspeq n\sspm r_n(k\,l)$ if $k\,l\sspnotequiv 0\,(mod\, n)$. \psn
We skip the proofs of the other two cases, $(+,-)$ and $(-,+)$, which run along the same line. Here one arrives first at $p_n(-(k\,l))$, and in order to compare it with $p(k\,l)$  both alternatives in {\it lemma 17} have to be considered like in the just considered  second alternative. \hskip 8cm $\square$ \psn  
For the computation of the cycle structure of the {\sl Galois} group ${\cal G}_n\sspeq {\cal G}{\sl al}(\mathbb Q(\rho(n))/\mathbb Q)$ we are only interested, due to eq. $(20)$, in odd numbers relatively prime to $n$. Contrary to ordinary $mod\, n$-arithmetic where the set of odd numbers $\mathbb O\sspdef \{2l+1\,|\,l\sspin \mathbb Z\}$ is in general not closed under multiplication (\eg\ $5\sspcdot 5\sspeq 25 \sspequiv 4\, (mod\,7)$), it will be shown that  $\mathbb O$ is closed under \Modd\,$n$ multiplication. Of interest are the units, the elements which have inverses, in order to see the expected group structure. First consider the reduced set ${\mathbb O}^*_n$, given by the odd numbers relatively prime to $n$. The negative odd numbers in this set are just the negative of the positive odd numbers, therefore it will suffice to consider  ${\mathbb O}^*_{n,>}\sspdef \{2l+1\,|\,l\sspin \mathbb N_0,\ gcd(2l+1,n)\sspeq 1\}$.  If $n$ is a power of $2$ the set ${\mathbb O}^*_2$ will be $\mathbb O$, with the \ogf  \dstyle{ G(x)\sspeq {\frac{x}{(1-x)^2}}\, (1\sspp  x)} for the sequence of positive odd numbers \dstyle{\{o^*_{2,>}(k)\sspdef 2\,k+1\}_{k=0}^\infty}. For the other even numbers $n$ only the odd numbers relatively prime to the squarefree kernel of $n$, called $sqfk(n)$, (see \seqnum{A007947},  encountered already several times), will enter the discussion.  Therefore, besides the just considered (trivial) case of the even prime $2$, the set ${\mathbb O}^*_n$ is  only relevant for  squarefree odd moduli, either prime or composite. We consider first the case of odd primes $n\sspeq p$, and give the \ogf of the sequences of numbers from  ${\mathbb O}^*_{p,>}$, called $o^*_{p,>}(k)$, as well as an explicit formula in terms of floor functions. These are the positive odd integers without odd multiples of the odd $p$.\psn
{\bf Proposition 11:\  Odd prime moduli, o.g.f. and explicit formula for $\bf \mathbb O^*_{p,>}$ elements}\psn
With \dstyle{G_p(x)\sspdef \sum_{k=0}^\infty\, o^*_{p,>}(k)\,x^k}, for odd primes $p$, one has
\Beq
G_p(x)\sspeq\frac{x}{(1-x^{p-1})\,(1-x)}\,\left\{1\sspp 2\,\sum_{k=1}^{\frac{p-3}{2}}\, x^k\,(1\sspp x^{\frac{p-1}{2}}) \sspp 4^{\frac{p-1}{2}}\sspp x^{p-1}\right\}\ ,
\Eeq
and
\Beq
o^*_{p,>}(k)\sspeq {\Caseszwei{$0$}{{\text if} \ $k=0$\ ,}{2\,k\sspm 1\sspp 2\,\floor{\frac{k+\frac{p-3}{2}}{p-1}}}{{\text if} \ $k\sspin \mathbb N$\ .}}
\Eeq
We have used the even number $0$ for $k\sspeq 0$ such that the $G_p(x)$ sum can start with $k\sspeq 0$.\psn
{\bf Proof}:  This {\it proposition} will become a corollary to the later treated  general case of odd squarefree moduli $n$ in {\it proposition 13}.\psn
{\bf Example 15:\  O.g.f.s and formula for reduced odd numbers for modulus $\bf p\sspeq 7$.}\psn
 \Beqarray
G_7(x)&\sspeq&\frac{x}{(1-x^6)\,(1-x)}\,\left\{1 + 2\,(x\sspp x^2) \sspp 4\,x^3 \sspp 2\,(x^4\sspp x^5) \sspp x^6 \right\}\ , \\
o^*_{7,>}(n) &\sspeq& 2\,n\sspm 1\sspp 2\,\floor{\frac{n+2}{6}},\  n\sspgeq 1\, . 
\Eeqarray
The instances for $p\sspeq 3,5,7,11,13$ and $17$ are found under \seqnum{A007310},  \seqnum{A045572},  \seqnum{A162699}, \seqnum{A204454}, \seqnum {A204457}, and  \seqnum {A204458}, respectively.\psn
In order to prepare for the general case of odd squarefree modulus $n$, we state a {\it  proposition} on the structure of the reduced odd numbers set  $\mathbb O^*_{n,>}$.\psn
{\bf Proposition 12:\ Mirror symmetry and modular periodicity of  $\bf \mathbb{O}^*_{n,>}$ for odd $\bf n$} \psn
{\bf i)} Mirror symmetry.  For $k\sspin \{1,2,...,\delta(n)\}$ one has: 
\Beq
 o^*_{n,>}(2\,\delta(n)-(k-1))\sspeq 2\,n\sspm o^*_{n,>}(k),
\Eeq
where  $o^*_{n,>}(k) \sspeq m_{k}(n)$ from ${\cal M}(n)$, given in {\it definition 2}, eq. $(65)$.\psn
{\bf ii)} $mod\ 2\,n$  periodicity:  For $k\in \mathbb N$ one has:
\Beq
o^*_{n,>}(k)\sspeq  o^*_{n,>}(k\sspp 2\,\delta(n))\, (mod\,2\,n).
\Eeq
Written as a relation between neighboring fundamental units, numbered by $N\sspgeq 1$, this becomes the following statement. 
For $k\sspin \{2\,(N-1)\,\delta(n)+1,\,...,\,2\,N\,\delta(n)\}$, with $N\sspin \{2,\,3,\,4,\,...\}$, one has: 
\Beq
o^*_{n,>}(k)\sspeq (N-1)\,2\,n \sspp o^*_{n,>}(k\sspm (N-1)\,2\,\delta(n))\ .
\Eeq
$\delta(n)$ is the degree of the minimal polynomial $C(n,x)$ for the algebraic number $\rho(n)$ introduced in {\it section 2}. Note that if \dstyle{n\sspeq \prod_{j=1}^{\omega(n)}\, p_j} with distinct odd primes $p$, and $\omega(n)\sspeq $\seqnum{A001221}$(n)$, then \dstyle{2\,\delta(n)\sspeq \prod_{j=1}^{\omega(n)}\, (p_j-1)}. $L(n)\sspdef 2\,\delta(n)$ is the length of the fundamental $N$-units.\psn
Before we give the proof consider {\it figure 3} for the case $n\sspeq 3\cdot 5 \sspeq 15$ with $\delta(n)\sspeq 4$. The second statement {\bf ii)} concerns the relation of the numbers of the second fundamental unit ($N\sspeq 2$) to the one in the $N\sspeq 1$ unit. \Eg  the odd number $37$ for $k \sspeq 2\cdot 4+2\sspeq 10$ is equal to $ 30 \sspp o^*_{15,>}(10\sspm 2\cdot 4)\sspeq 30 \sspp o^*_{15,>}(2)\sspeq 30\sspp 7$, which checks.
The statement {\bf i)} shows the mirror symmetry within the first (and any other) unit of length $L(15)\sspeq 2\cdot 4\sspeq 8$. In {\it figure 3} this symmetry is indicated by the brackets below the first unit, and it is  a symmetry around the missing number $n\sspeq 15$. Missing numbers have been indicated by a dot. It is the pattern of missing odd numbers which is mirror symmetric, not the one of the actual values of the odd numbers. But the relation between the odd numbers in the second half of a unit and the first one then follows, and is given by the statement of the {\it proposition}. In {\it figure 3} $P(n)\sspdef 2\,n\sspeq 30$ is the shift for the $o^*_{>,n}$ values from the $N\sspeq 1$  to the $N\sspeq 2$ unit (or any of the neighboring units), and $p(n)\sspdef n\sspp 1\sspeq 16$ is the shift for these values from the first half of every unit to the second half.\psn
%%%%%%%%%%%%%%%%%%%%%%%%%%%%%%%%%%%%
%% Figure 3: Odd numbers mod n structure
\pbn 
\begin{center}
\parbox{16cm}{
{\includegraphics[height=8cm,width=.8\linewidth]{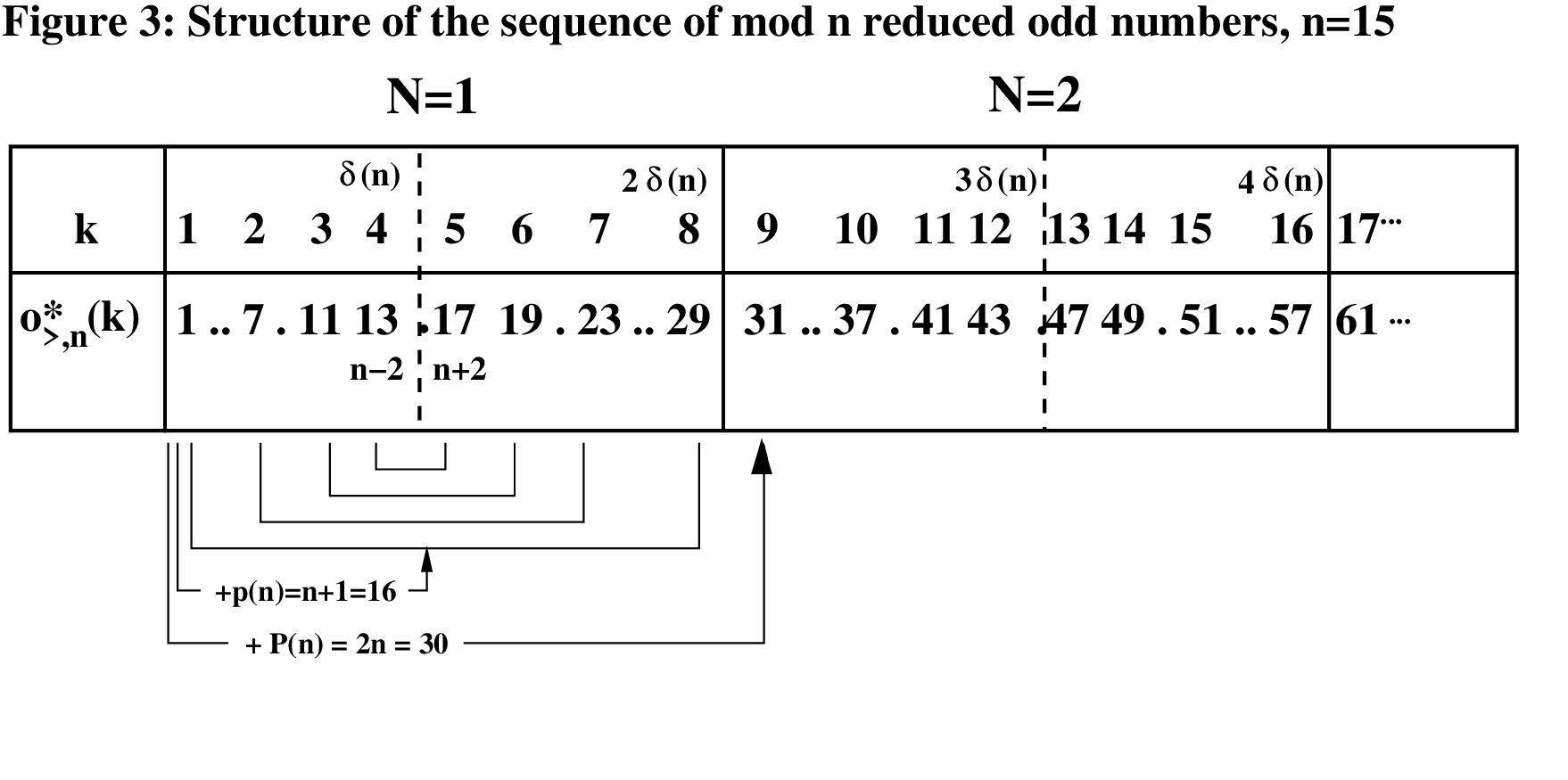}}
}
\end{center}
\pbn
%%%%%%%%%%%%%%%%%%%%%%%%%%%%%%%%%%%%
{\bf Proof:  i)} It is clear from the degree of the minimal polynomial and eq. $(20)$ that the number of $mod\, n$ reduced positive odd numbers smaller than $2\,n$ is $\delta(2\,n)\sspeq 2\,\delta(n)$ (from the definition of $\delta(n)$ in terms of {\sl Euler}'s $\varphi$ function). We now determine $\delta(n)$ reduced odd numbers $ mod \,n$ which lie between $n$ and $2\,n$, by mirroring the $\delta(n)$ elements of ${\cal M}(n)$ around the position where the missing number $n$ is situated, which is between the position $k\sspeq \delta(n)$ and the next one. The mirror symmetry refers to the gaps. The number  $m_k(n)$ from  ${\cal M}(n)$  will be mapped  to $o^*_{n,>}(2\,\delta(n) - (k-1))$  at the mirrored position in the second half of the later defined first fundamental unit ($N\sspeq 1$). To find this odd number  the first difference set $\triangle\, \widehat {\cal M}(n)$, introduced in connection with ${\widehat {\cal M}}(n)$ from the {\it definition 3}, becomes important in order to count the gaps in the sequence of odd numbers when they are reduced $mod\, n$. $o^*_{n,>}(2\,\delta(n) - (k-1))$  is found by  adding to $m_k(n)$ twice the value of the sum of the gaps from the position $k$ to the center, the mirror axis. This is $2\,(\triangle m_{k+1}(n)\sspp \triangle m_{k+2}(n)\sspp  ... \sspp \triangle\, m_{\delta(n)}(n) \sspp 2)$. The $2$ in this sum is half the gap-length from the value $n-2$ at the position $\delta(n)$ and $n+2$ at the next position, the mirror-position $\delta(n)+1$ of $\delta(n)$. This is a telescopic sum which becomes $2\,(-m_k(n)\sspp m_{\delta(n)}(n)\sspp 2)$. Therefore, the value of  $o^*_{n,>}(2\,\delta(n) - (k-1))$ is $m_k \sspp 2\,(-m_k(n)\sspp m_{\delta(n)}(n)\sspp 2)\sspeq 4\sspp 2\,m_{\delta(n)}(n)\sspm m_{k}(n) \sspeq 2\,n\sspm m_{k}(n)$, because always $m_{\delta(n)}(n)\sspeq n-2$ holds.  See {\it figure 3}, $n=15$ with the values $30\sspm 1\sspeq 29$, $30\sspm 7\sspeq 23$, $30\sspm 11\sspeq 19$ and  $30\sspm 13\sspeq 17$ for $k\sspeq 1,\,2,\,3$ and $4$, respectively.  Now it is clear that this mirroring leads to the correct number of reduced odd numbers for the second half of the fundamental $N\sspeq 1$ unit, because $gcd(m_k(n),n)\sspeq 1$, as member of ${\cal M}(n)$, implies  for the mirror image also $gcd(o^*_{n,>}(2\,\delta(n)-(k-1)),n)\sspeq 1$. All positive odd numbers relatively prime to $n$ and not exceeding $2\,n\sspm m_{1}(n)\sspeq 2\, n\sspm 1 $ have thus be found.\psn
Proof of {\bf ii)}:  From above we know that the number at position $k\sspeq 2\,\delta(n)+1$ is $2\,n+1$ because the one for $k\sspeq 2\, \delta(n)$ has been shown to be $2\,n-1$ and $gcd(2\,n+1 ,n)=1$ (indirect proof by assuming the contrary, using odd $n$; this is similar to the proof given in connection with {\it definition 3} of $\widehat {\cal M}(n)$). Now it is clear that a shift in the $o^*_{n,>}$ numbers with $P(n)\sspdef 2\,n$ leads from the fundamental unit No. $N\sspeq 1$ to the second one, $N\sspeq 1$, by putting  $o^*_{n,>}(2\,\delta(n)\sspp k)\sspeq  o^*_{n,>}(k)\sspp P(n)$, for $k\sspin \{1,2,...,\delta(n)\}$.  This is obvious because the $gcd$ value $1$ is not changed by adding $2\,n$. This process can be iterated to find the $mod\, 2\,n$ periodicity structure stated in {\bf ii)}. See {\it figure 3}, $n=15$ with  $k=14$, $N\sspeq 2$: $o^*_{15,>}(14)\sspeq 49 \sspeq 2\cdot 15 \sspp o^*_{15,>}(14\sspm 1\cdot 2\cdot 4)\sspeq 30 \sspp \  o^*_{15,>}(6)\sspeq 30\sspp 19$. \hskip 3.7cm $\square$ 
The sequences $\{o^*_{n,>}\}$ for $n\sspeq 15$ and $ n\sspeq 21$ are found in \seqnum{A007775} and  \seqnum{A206547}, respectively.\psn
The sequences of the \Modd\,$n$ residues of the numbers $o^*_{n,>}(k)$, for $k\sspeq 1,\,2,\,...$, for prime moduli $n\sspeq p\sspeq 3,\, 5,\, 7,\, 11,\, 13,\, 17$, and for the first odd composed ones for $n=15,\, 21$ are found in \seqnum{A000012}, \seqnum{A084101}, \seqnum{A110551},  \seqnum{A206543}, \seqnum{A206544}, \seqnum{A206545},  and \seqnum{A206546}, \seqnum{A206548}, respectively.\psn 
For the following formula for the nonnegative odd numbers reduced $mod\, n$ involving floor functions we need the following list (increasingly ordered set) ${\cal F}(n)$ of length $L(n) \sspeq 2\,\delta(n)$ derived from the list $\triangle \widehat{\cal M}$.\psn
{\bf Definition 5: \ List $\mathbfcal {\cal F}(\bf n)$}\psn
$ {\cal F}(n)\sspeq \{f(n,1),...,f(n,2\,\delta(n)\}$ with
\Beqarray
f(n,j) &\sspeq& {\frac{\triangle m_{j+1}(n)\sspm 2}{2}}, \ \text{for}\ j\sspin \{1,2,...,\delta(n)-1\} \nonumber \\
f(n,\delta(n)) &\sspeq&1,\ \ f(n,\delta(n)\sspp j)\sspeq f(n,\delta(n)\sspm j),\ \text{for} \  j\sspin \{1,2,...,\delta(n)-1\},\  \text {and}\nonumber \\ 
 f(n,2\,\delta(n)) &\sspeq& 0.
\Eeqarray
Here $\triangle m_{k}(n)\sspdef m_k(n)\sspm m_{k-1}(n)$, (the $k$th element  of  $\triangle\widehat {\cal M}(n)$). This list ${\cal F}(n)$ is obtained from first enlarging  $\triangle\widehat {\cal M}$  by mirroring  the first $\delta(n)$ entries at the last element $4$, to obtain a list of order $2\,\delta(n)+1$. Then the first and last element of this new list is put to zero and all the other elements are diminished by 2, thus obtaining a set of only even numbers. Then one divides by $2$, omits the final $0$, and reverses the remaining list.  \psn
{\bf Example 16:}\ $\bf {\mathbfcal {\cal F}}(15)$\psn
$n\sspeq 15\sspeq 3\cdot 5$, $2\,\delta(15)\sspeq 2\cdot 4\sspeq 8$, ${\cal M}(15)\sspeq  \{1, 7, 11, 13\}$, $\widehat{\cal M}(15)\sspeq  \{0, 1, 7, 11, 13, 17\}$, $\triangle\widehat{\cal M}(1
5)\sspeq \{1, 6, 4, 2, 4\}$, the mirror extension is $\{1, 6, 4, 2, 4,2,4,6,1\}$, the reduction step leads to $\{ 0, 4, 2, 0, 2,0,2,4,0\}$, and finally, dividing by $2$,  omitting the last $0$ and reverting, leads to ${\cal F}(15)\sspeq \{2,1,0,1,0,1,2,0\}$. Except for the last $0$ there is a mirror-symmetry around the fourth entry $1$. See the first row of {\it table 5}. \psn 
For the odd squarefree composite numbers $n\sspeq$\seqnum{A024556}$(m)$, $m\sspeq 1,\,2,\, ...\,,17$, see {\it table 5} for $\widehat {\cal M}(n)$, $\triangle\,\widehat{\cal M}(n)$ and ${\cal F}(n)$.\psn 
{\bf Proposition 13:\  O.g.f. and formula  for $\bf \mathbb O^*_{n,>}$ elements.} \psn
With \dstyle{G_n(x)\sspdef \sum_{k=1}^{\infty}\, o^*_{n,>}(k)\,x^k}, for odd $n$, one has
\Beq
G_n(x) \sspeq  \frac{x}{(1-x^{2\,\delta(n)})\,(1-x)}\,\left\{1\sspp \sum_{k=1}^{2\,\delta(n)-1}\,\triangle o^*_{n,>}(k+1)\,x^k  \sspp x^{2\,\delta(n)}\right\}\    
\Eeq
with the first differences $\triangle o^*_{n,>}(j)\sspdef o^*_{n,>}(j)\sspm o^*_{n,>}(j-1) $\, .\psn
This generates 
\Beq
o^*_{n,>}(k)\sspeq 2\,k-1 + 2\, \sum_{j=1}^{2\,\delta(n)}\,f(n,j)\,\floor{\frac{k+(j-1)}{2\,\delta(n)}} .
\Eeq
Note that the numerator polynomial in $G_n(x)$ needs the first differences of the sequence members of the fundamental $N\sspeq 1$ unit, which due to the mirror symmetry of {\it proposition 12} can be reduced to the first differences of $\{m_2(n),\,m_3(n),\,..,\,m_{\delta(n)},n+2\}$ with $m_{\delta(n)}\sspeq n\sspm 2$. \psn
One could use, like in {\it proposition 11}, $o_{n,>}(0)\sspeq 0$ and let the sum in $G_n(x)$ start with $k\sspeq 0$.\psn 
{\bf Proof:} 
 We start with the periodicity of the sequence due to {\it proposition 12 ii)}. Taking the difference $o^*_{n,>}(k)\sspm o^*_{n,>}(k-1)\sspeq o^*_{n,>}(k-2\,\delta(n))\sspm o^*_{n,>}(k-1-2\,\delta(n))$ leads to the recurrence (we omit all unnecessary indices) $o(k)\sspeq o(k-1)\sspp o(k-2\,\delta(n))\sspm o(k-1-2\,\delta(n))$, for $k\sspgeq 2\,\delta(n)+1$. Here one needs also $o^*_{n,>}(0)\sspdef -1$. For some of these recurrences the {\it o.g.f.}s were determined by {\sl R. J. Mathar} in \eg \, \seqnum{A045572} and \seqnum{A162699}. I general \dstyle{G_n(x)\sspeq \sum_{k=1}^{2\,\delta(n)}\, o(k)\, x^k \sspp x\,\sum_{k=2\,\delta(n)+1}^{\infty}\,o(k-1)\,x^{k-1} \sspp x^{2\,\delta(n)} \, \sum_{k= 2\,\delta(n)+1}^{\infty}\,o(n-2\,\delta(n)) \sspm x^{2\,\delta(n)+1}\,\sum_{k= 2\,\delta(n)+1}^{\infty}\,o(n-2\,\delta(n)-1)\,x^{n-2\,\delta(n)-1}}, where the recurrence has been used (and the infinite sum has been reordered, not bothering about absolute convergence, in the sense of formal power series). Shifting the indices, one arrives at \dstyle{G_n(x)\sspeq\sum_{k=1}^{2\,\delta(n)}\, o(k)\, x^k \sspp x\,\left\{ G_n(x)\sspm ( o(1)\, x\sspp ....\sspp o(2\,\delta(n)-1)\,x^{2\,\delta(n)-1)}\right\} \sspp x^{2\,\delta(n)}\,G_n(x) \sspm x^{2\,\delta(n)+1}\,\left(G_n(x)\sspp (-1)\right)}, where the $-1$ resulted from putting $o(0)=-1$ (see above). This rearranges into $(1\sspm x\sspm x^{2\,\delta(n)}\sspp x^{2\,\delta(n)+1 })\, G\,_n(x)\sspeq x\,\left(1\sspp \triangle o(2)\,x \sspp ...\sspp \triangle o(2\,\delta(n))\, x^{2\,\delta(n)-1}\sspp x^{2\,\delta(n)} \right)$. Factorizing the bracket on the {\it l.h.s.} and division leads to the claimed form for $G_n(x)$. Of course, the denominator factor $(1-x^{2\,\delta(n)})$ can be factorized into cyclotomic polynomials.
\psn
For the proof of the second part we use $L \sspdef 2\,\delta(n)$. Given $G_n(x)$ with the $L\sspp 1$ input coefficients of the numerator polynomial one derives the claimed explicit form for $o^*_{n,>}(k)$, by defining first the sequence with entries \dstyle{b_L(k)\sspdef \floor{\frac{k+L-1}{L}}}, generated by \dstyle{\frac{x}{(1-x^L)\,(1-x)}} (partial sums of the characteristic sequence for multiples of $L$, then shifted).  The numerator polynomial leads for $o^*_{n,>}(k)$ to a sum of $L+1$  floor-functions with decreasing arguments, starting with \dstyle{\floor{\frac{k-0+L-1}{L}}}, ending with \dstyle{\floor{\frac{k-L+L-1}{L}}\sspeq\floor{\frac{k-1}{L}}}, and corresponding coefficients.  In order to find a standard form  for this sum we use the following floor-identity 
\Beq
\sum_{j=0}^{L-1}\,\floor{\frac{k+j}{L}}\sspeq k,\  \text{for}\ L\sspin \mathbb N,\ \text{and}\   k\sspin \mathbb Z\  , 
\Eeq
which can be seen from  \dstyle{\sum_{j=0}^{L-1}\, (k+j)\, (mod\, L)\sspeq \frac{L\,(L-1)}{2}},  $k\sspin \mathbb Z$, which is trivial (just add the $L$ terms which are $0,\,1,\,...\,,L-1$ in a certain cyclic order), and the relation between $mod\, L$ and the floor function\pn
\dstyle{\floor{\frac{k}{L}} \sspeq \frac{1}{L}\,(k\sspm k\,(mod\,L))}. This identity allows us to lower $L$ consecutive coefficients of this sum of $L\sspp 1$ terms by $1$, producing a term $k$ if one uses the identity for the first  $L$ terms because the second to last floor-argument is then  $\floor{\frac{k+0}{L}}$, and the identity is read backwards. If one uses the identity for the last $L$ terms one produces a $k-1$, because the last term has $ \floor{\frac{k-1}{L}}$. When we apply this identity twice in the described way we pick up  $n\sspp (n-1)\sspeq 2\,n\sspm1$ and the first and last coefficient, which were originally $1$, become $0$, and all other coefficients are diminished by $2$ because they participate in both applications of the identity. It is guaranteed that all coefficients are now even and $\sspgeq 0$, because the coefficients except the first and last one were even and $\sspgeq 2$, because these numerator polynomial coefficients resulted from first differences of the sequence of odd numbers. Therefore one can extract a factor $2$ and if the floor-functions are written with increasing arguments, starting with \dstyle{\floor{\frac{k}{L}}} (the second to last term in the original order), one has exactly the coefficients given by the list ${\cal F}(n)$ of length $L$ of {\it definition 5}. This proves the explicit form for $o^*_{n,>}(k)$. \hskip 11cm $\square$\psn
{\bf Example 17:  O.g.f. $\bf G_8(x)$}\psn
\Beq 
G_8(x)\sspeq x\,\frac{(1+2\,(x \sspp x^2 \sspp x^3 \sspp x^4 \sspp x^5 \sspp x^6 \sspp x^7) \sspp x^8}{(1-x^8)\,(1-x)}\sspeq x\,\frac{1\sspp x}{(1\sspm x)^2}\, 
\Eeq
generating, with offset $0$, the odd numbers \seqnum{A004273}.\psn 
Before coming to the multiplicative group \Modd\,$n$ we define reduced residue systems 
\Modd\,$n$ as well as reduced odd residue systems \Modd\,$n$. See \eg \cite{Apostol}, p. 113 for the $mod\, n$ case. 
\psn
\vfill
\eject
\noindent
{\bf Definition 6a: Reduced residue system Modd\, $\bf n$}\psn
A reduced residue system \Modd\,$n$ ({\it RRSn}) is any set of $\varphi(n)$  pairwise incongruent \Modd\,$n$ numbers, each of which is relatively prime to $n$. \psn
\Eg  $n=15$, $\varphi(15)\sspeq 8$, $\{1,2,4,7,8,11,13,14\}$ or $\{29,32,26,23,22,41,17,16\}$, {\it etc}. The first one is the smallest positive one.  These systems will not play a r\^ole later on.\psn
 {\bf Definition 6b: Reduced odd residue system Modd\, $\bf n$}\psn
A reduced odd residue system \Modd\,$n$ ({\it RoddRSn}) is any set of $\delta(n)$ odd pairwise incongruent \Modd\,$n$ numbers, each of which is relatively prime to $n$ \psn
\Eg  $n=15$, $\delta(15)\sspeq 4$, $\{1,7,11,13\}$ or $\{29,27,19,17 \}$. The first one coincides with ${\cal M}(15)$ from eq. $(65)$ and is the smallest positive reduced residue system \Modd\,$n$. Remember that the parity of the members in each \Modd\,$n$ residue class, not in class $[0]$, is the same (see {\it lemma 12}). Later on we will restrict ourselves mostly to the system ${\cal M}(n)$.\psn
We next  study the multiplicative group \Modd\,$n$. The elements are the residue classes $[m_j]$ , for $j\sspeq 1,\,2,\, ...,\,\delta(n)$, corresponding to the members of the reduced residue system ${\cal M}(n)$. In fact, we will take these representatives, multiplying \Modd\,$n$, as we have done above.\psn
In order to see the group structure one first convinces oneself that this set ${\cal M}(n)$ is closed under \Modd\,$n$ multiplication. This follows from {\sl proposition 10} and {\it lemma 12} which showed that these classes have only odd numbers, and every odd numbers appear exactly once because of the definition of these classes. The associativity, commutativity and the identity element $1$ are also clear. We do not have a formula how to find the inverse element $m_j^{-1}$ of $m_j$ but coming to this \Modd\,$n$ multiplication from the study of automorphisms of the splitting field for the minimal polynomial $C(n,x)$ for the algebraic number $\rho(n)$, it is clear that these inverses have to exist (the invariance property defines a group). We are dealing here with the {\sl Galois} groups for these polynomials. The cycle structure for $n\sspeq 1,\, 2,\,...,\,40$ is given in {\it table 6}. Because each group is Abelian one has for every group of  prime order a cyclic group, which can be checked for the examples given in {\it table 6}.\psn
{\bf Remark 6:\ } The multiplicative group \Modd\,$n$ is cyclic if $\delta(n)$ is prime.\psn
This is a corollary on the fundamental theorem on finite Abelian groups (see \eg \cite{Speiser}, p. 49, see also \cite{Cox}, p. 511, {\sl Cauchy}'s Theorem A.1.5), or use {\sl Lagrange}'s theorem on the order of subgroups, and the fact that the powers of an element, not the identity, generate a cyclic subgroup, to prove that in fact every group of prime order is cyclic. The values $n$ for which the order $\delta(n)$ is prime are given in \seqnum{A215046}.
This is the sequence $[4, 5, 6, 7, 9, 11, 23, 47, 59, 83, 107, 167, 179, 227, 263, 347, 359, ...]$. Of course, there are other values $n$ with cyclic \Modd\,$n$ group, like $n\sspeq 2, 3, 8, 10, ...$. For $n\sspeq 1$ one has the trivial case with cycle structure $[[0]]$, also a cyclic group, \viz\  $Z_1$. One can give a more general sequence of $n$ numbers with cyclic  ${\cal G}_n \sspiso$ multiplicative \Modd\,$n$ group, namely \seqnum{A210845}.\psn
{\bf Remark 7:} The multiplicative group \Modd\,$n$ is cyclic if $\delta(n)$ is squarefree.\psn
 The squarefree numbers are given in \seqnum{A005117}, namely  $[1, 2, 3, 4, 5, 6, 7, 9, 11, 13, 14, 18, 21, 22, 23, 25, 29,$\pn 
$31,...]$. This {\it remark 7} results from the fact that \seqnum{A000688}, giving the number of Abelian groups of order $n$, is $1$ exactly for the squarefree numbers \seqnum{A005117}. See the formula, based on the {\sl H.-E. Richert} reference quoted there.  Because for each order there is at least the cyclic group these values lead necessarily to a cyclic group. The above given \seqnum{A215046} values are a proper subset of  those from \seqnum{A210845}.  There are. however,  still more values $n$ with a cyclic \Modd\,$n$ group.  Missing are \eg $8, 10, 15, 16, 17, ... $. All the $n$ values with cyclic \Modd\,$n$ group  are in \seqnum{A206551}. The complementary sequence is \seqnum{A206552}, giving the $n$ values with non-cyclic \Modd\,$n$ group. See {\it Table 8} for all values $n\sspleq 100$.  \psn
As an aside we remark that special squarefree numbers are the so called cyclic numbers \seqnum{A003277}. It is known that if the order of a (finite) group is a cyclic number then there is only one group, the cyclic one. See a comment on \seqnum{A003277}. Therefore \seqnum{A000001}$(n)$, the number of groups of order $n$ is $1$ if $n$ is a cyclic number, in fact, the reverse also holds: if  \seqnum{A000001}$(n)\sspeq 1$ then $n$ is a cyclic number. This can be taken as an alternative definition for cyclic numbers because then there is only one (non-isomorphic) group of this order which has to be the cyclic group. See also Yimin Ge's Math Blog [YiminGe], where the `only if' statement in the proposition may be misleading but in the proof the given statement is correct. \Eg for order $6$ (not a cyclic number) there are groups other than $Z_6$. In fact \seqnum{A000001}$(6)\sspeq  2$, and there is the (non-Abelian) group D(3) (dihedral group). (It is clear that cyclic numbers are not the numbers $n$ for which the multiplicative group $mod\, n$, which is Abelian of order $\varphi(n)$, is cyclic. These numbers are given in  \seqnum{A033948}).  Because we are dealing with Abelian groups the squarefree numbers are more interesting here. \psn
Recall the situation for the {\sl Galois} group for the cyclotomic polynomials (the minimal polynomials for $\zeta(n)$, an $n$-th root of unity) which is isomorphic to the multiplicative (Abelian) group of integers modulo $n$.  See {\it Table 7} were we have listed these non-cyclic groups for $n\sspin \{1,\,...\,,100\}$. See also \seqnum{A282624} for a W. Lang link with a table for $n\sspin \{1,\,...\,,130\}$ also with generators. This $n-$value sequence is known as \seqnum{A033949}. (The case $n\sspeq 15$ is an exercise in \cite{Escofier}, p. 159, Ex. 9.6. 3)a)). The values $n$ with a cyclic {\sl Galois} group coincide with the moduli $n$ which possess primitive roots $r\sspeq r(n)$, \ie the order of $r$ modulo $n$ is $\varphi$: $r^{\varphi(n)}\sspequiv 1\,(mod\, n)$, and no smaller positive exponent $k$ satisfies this congruence. These moduli $n$ are known to be exactly $p^e,\,2\, p^{e'},\, 1$ and $2$, with some positive powers $e$ and $e'$. See \eg  \cite{NZM}, {\it Theorem 2.41}, p. 104.  All other $n$ lead to non-cyclic multiplicative groups modulo $n$ . See \seqnum{A033949}. For the smallest primitive roots in this case see \cite{ASt}, pp. 864-869, the column called $g$. We do not know a similar characterization of the non-cyclic numbers $n$ for the \Modd\,$n$ group which are shown in {\sl Table 8}. However, these numbers have  also to appear in \seqnum{A033949}. It is clear that we should determine primitive roots $r\sspequiv r(n)$ for the \Modd\,$n$ multiplication, \ie find those $r$ from ${\cal M}(n)$ which have order $\delta(n)$: the smallest positive $k$ which satisfies $r^k\sspequiv +1\,($\Modd\,$n$)  is $k\sspeq \delta(n)$. See \seqnum{A206550}  for these smallest positive primitive roots \Modd\,$n$, where a $0$ entry, except for $n\sspeq 1$, indicates that there exists no such primitive root. This sequence starts, with offset $1$, as $[0, 1, 1, 3, 3, 5, 3, 3, 5, 3, 3, 0, 7, 5, 7, 3, 3, 5, 3, 0, 11, 3, 3, 0, 3, 7,...]$. In general one does not expect a formula for these primitive roots $r\sspeq r(n)$, because also in the multiplication modulo $n$ case there is no one available. \psn
{\bf Proposition 14: Number of Modd\,$\bf n$ primitive roots}\psn
If a primitive root $r$ for the multiplicative group \Modd\,$n$ exists there are $\varphi(\delta(n))$ of them.\psn
This is the sequence \seqnum{A216322}. For example, $n\sspeq 13$, $\delta(13)\sspeq 6$, $\varphi(6)\sspeq 2$ with  $r_1\sspeq 7$ and $r_2\sspeq 11$; $n\sspeq 14$, $\delta(14)\sspeq 6$, $\varphi(6)\sspeq 2$ with  $r_1\sspeq 5$ and $r_2\sspeq 11$. \psn
{\bf Proof:} If a \Modd\,$n$ primitive root exists then the multiplicative group \Modd\,$n$ is cyclic.  The order of this group is $\delta(n)$, Then the number of pairwise incongruent primitive roots is obtained exactly like in the case of primitive m-roots of unity on the unit circle by the number of relatively prime numbers less than $m$ which is {\sl Euler}'s $\varphi(m)$. Here $m\sspeq \delta(n)$.\hskip 10cm $\square$\psn
Now the question whether the multiplicative group \Modd\,$p$ with $p$ a prime is cyclic is answered. Up to now we know from {\it remark 7} the positive answer only for those primes with \dstyle{\delta(p)\sspeq \frac{p-1}{2}} squarefree. These are the primes given in  \seqnum{A066651}. The case $p\sspeq 2$ with $\delta(2)\sspeq 1$ is trivially cyclic. For general prime $p$ one wants to show that there exists a primitive root in the multiplicative \Modd\,$p$ group. This is analog to the $mod\, p$ case, where a proof can be found in \cite{NZM}, {Theorem  2.36}, p. 99, or in \cite{Apostol}, ch. 10.4, pp 206 ff. One has, however, to be careful to use only the multiplicative group structure in the proof. Already in connection with {\it lemma 10} we have given a warning to use theorems valid in the modular arithmetic $mod \,n$ in our case. Another failure occurs for  the theorem \cite{NZM}, {\it Theorem 2.3(c)}, p. 49,  \eg  $x\sspequiv y\,(mod\, m_1), x\sspequiv y\,(mod\, m_2)$ if and only if \dstyle{x\sspeq y\,\left (mod\, \frac{m_1\,m_2}{gcd(m_1,m_2)}\right)}. For the  \Modd\,$n$ case one has the counterexample $19\sspeq 5\,(Modd\, 3)$ and $19\sspeq 5\, (Modd\, 7)$ but $19$ is obviously not congruent $5$\,(\Modd\,$21$).  Several {\sl lemmata}, the analoga of $mod\ n$ facts, are collected before stating {\it proposition 15}.\psn
{\bf Lemma 19: Analogon of the Fermat-Euler Theorem}\psn
If $a$ is odd and  $gcd(a,n)\sspeq 1$ then $a^{\delta(n)}\sspequiv +1\,(Modd\, n)$, with $delta(n)\sspeq $\seqnum{A055034}$(n)$ (see the start of {\it section 3)}.\psn
{\bf Proof:} Analog to \eg \cite{Apostol}, {\it Theorem 5.17}, p. 113, with $m\sspto n,$ $\varphi\sspto \delta$, the order of the set (reduced odd residue system) ${\cal M}(n)$ of eq. $(65)$, and modulo $\sspto$ {\it Modd\,}$n$ (we avoid the term Moddulo). Instead of ${\cal M}(n)$ one can use any other reduced odd residue system \Modd\,$n$ (see {\it definition 6b}), with $m_j(n)$, $j\sspin \{1,\,2,\,...\,,\delta(n)\}$, replaced by any member of the residue class $[m_j(n)]$. The cancellation used at the end of the proof in \cite{Apostol} can be done here by multiplying with the existing inverses of the $b_i$s there. We are dealing with the multiplicative group Modd\,$n$. \hskip 17cm $\square$\psn
Next follows the definition of the order of an element from a reduced odd residue system  \Modd\,$n$.\psn
{\bf Definition 7: Modd\, $\bf n$ order of $a$ from a RoddRSn}\psn
The \Modd\,$n$ order of a positive (odd) integer $a$ from a reduced odd residue system \Modd\,$n$ is the smallest positive integer $h$ such that $a^h\sspequiv 1\,$(\Modd\,$n$).\psn
\Eg $n\sspeq 10, a\sspeq 9, h\sspeq 2$ because $9^1\sspeq 9\,$(\Modd\,$10$), and $9^2\sspeq 1\,$(\Modd\,$10$). See the table \seqnum{A216320} corresponding to the \Modd\,$n$ orders $h$ of the smallest positive {\it RoddRSn} members \seqnum{A216319}.\psn
{\it Lemma 19} guarantees the existence of a \Modd\,$n$ order $h\sspleq \delta(n)$ for each positive odd number $a$ with $gcd(a,n)\sspeq 1$.\psn
{\bf Example 18:} $n\sspeq 12$, $\delta(12)\sspeq 4$, ${\cal M}(12)\sspeq \{1,5,7,11\}$, $1^4\sspequiv 1\, \text{(\Modd}\, 12)$, $5^2\sspequiv 1\,\text{(\Modd}\, 12) \sspfollows 5^4\sspequiv \text{(\Modd}\, 12)$, similarly for $7$ and $11$. See the cycle structure for $n\sspeq 12$ in {\it table 6}.\psn
{\bf Lemma 20: The Modd\,$\bf n$ order divides  $\bf \boldsymbol{\delta}(n)$}\psn 
If $h$ is the \Modd\,$n$ order of $a$ then $h|\delta(n)$. Moreover, $a^j\sspequiv a^k$ (\Modd\,$n$), {\it w.l.o.g.} $j\sspgr k$, if and only if $h|(j-k)$.\psn
{\bf Proof:} The first part is analog to \cite{NZM}, {\sl Corollary 2.32}, p. 98, with $\varphi\sspto \delta$. It uses {lemma 19}, the {\sl Fermat-Euler} analogon.  Remember that always $gcd(a,n)\sspeq 1$ from the order definition.  The second part, for which one can use the group property, follows from the analog of \cite{NZM}, {\it Lemma 2.31}, p. 98.
(In the older German version of this book, vol. I,  the analog of this {\it lemma 20} appears as {\it Satz 22.3} on p. 63.)\hskip 18cm $\square$ \psn
{\bf Lemma 21: On the Modd\,$\bf p$ order of $\bf a$, $\bf p$ a prime}\psn
If $h$ is the {\it Modd}\,$p$ order of $a$, with $p$ a prime, one has $(a^k)^h \sspequiv 1$\,(Modd\,$p$)  for all $k$, and $1\sspeq a^0,\, a^1,\, a^2,\,...\,, a^{h-1} $ are pairwise incongruent \Modd\,$p$.\psn
{\bf Proof:} First part:  $(a^k)^h\sspeq (a^h)^k$, $a^h\sspequiv 1\,$(\Modd\,$p$) 
and \Modd\,$p$ respects multiplication (see {\it proposition 10}). Second part: assume the contrary, \ie $a^i\sspequiv a^j$\,(\Modd\,$p$), $1\sspleq j\sspkl i\sspleq h-1$. Apply {\it lemma 20} for $n=p$, showing that $h|(i-j)$, but  $i-j\sspleq i\sspleq h-1$ leading to a contradiction. \hskip 4cm $\square$\psn
{\bf Lemma 22: $\bf X^h\sspequiv 1\, (\text{\bf Modd}\, p)$ has at most $h$ incongruent solutions}\psn
The congruence $ X^h\sspequiv 1$\, (\Modd\,$p$), $h$ a positive integer, has at most $h$ pairwise incongruent solutions \Modd\,$p$.\psn
{\bf Proof:} If \dstyle{\floor{\frac{X^h}{p}}} is even then $X^h\sspequiv 1\,(mod\,p)$. and due to the $mod\, p$ theorem, \eg \cite{Nagell}, {\it Theorem 42}, p. 80, there are at most $h$ incongruent solutions. If \dstyle{\floor{\frac{X^h}{p}}} is odd then $-X^h\sspequiv 1\,(mod\,p)$, \ie $X^h\sspequiv (p-1)\,(mod\,p)$ which again has at most $h$ incongruent solutions.\hskip 5.5cm $\square$\psn
This is a weak statement, but sufficient for the following. One could try to prove that the number of  \Modd\, $n$ incongruent solutions of the congruence $ X^h\sspequiv 1$\, (\Modd\,$p$) is  \dstyle{gcd(h,\frac{p-1}{2})} if $p$ is odd. This number is trivially $1$ for $p\sspeq 2$.\psn 
{\bf Lemma 23: On the Modd\,$\bf n$ order of $\bf a^k$}\psn
If $h$ is the  \Modd\,$ n$ order of $a$,  then \dstyle{\frac{h}{gcd(h,k)}} is the  \Modd\, $ n$ order of $a^k$.\psn
{\bf Proof:} This is the analogon of \cite{NZM} {\it lemma 2.33}, p. 98, using the present {\it lemma 20}, part 2. \hskip 4cm $\square$\psn 
Now we are ready for the following {\it proposition}.\psn 
{\bf Proposition 15: Existence of Modd\,$\bf p$ primitive roots}\psn
For every prime $p$ there exists a primitive root for the multiplicative group \Modd\,$p$.\psn
{\bf Proof:} One shows that there are precisely $\varphi(\delta(p))$ primitive roots and because this number is always $\geq 1$ the claim will then follow.\pn
Following the proof of \cite{NZM}, {\it Theorem 2.36}, p. 99, we infer from the present {\it lemma 19} (with $n\sspto p$) that every positive odd integer $a$ with $gcd(a,p)\sspeq 1$ has a \Modd\,$p$ order $h\sspequiv h_p(a)$, \dstyle{1\sspleq h\sspleq \delta(p)\sspeq \frac{p-1}{2}}, and from {\it lemma 20} $h|\delta(p)$. Moreover, $(a^h)^k\sspequiv \,1$(\Modd\, $p$) for all $k$, and $1,\, a^1,\, a^2\, ...\, a^{h-1}$ are $h$ pairwise incongruent odd numbers \Modd\,$p$. From {\it lemma 22} these are all the solutions of the congruence $X^h\sspequiv 1\,$ (\Modd\,$p$). {\it Lemma 23} shows that there are $\varphi(h)$ numbers $a^k$, with $k\sspin \{1,\,...\,h-1\}$ which have \Modd\,$p$ order $h\sspequiv h_p(a)$ because $\varphi(h)$ is the number of such $k$ with $gcd(k,h)\sspeq 1$. The case $h\sspeq 1$ is covered with $\varphi(1)\sspdef 1$.\pn
With \cite{Apostol}, p. 207, we define for each $h|\delta(p)$ the set 
\dstyle{A_p(h)\sspdef {\set{a}{a\sspin {\cal M}(p)\ {\text {and\ \Modd}}\,p\  {\text{order of}}\, a\ {\text{is}}\ h}}}. There are $\tau(\delta(p))$ ($\tau$\sspeq \seqnum{A000005}, the number of divisors) such disjoint sets. This is the sequence $[1,\,1,\,2,\,2,\,2,\,4,\,4,\,3,\,2,\,4,\,4,\,6,\,  ...] \sspeq$\seqnum{A216326}.  \Eg $p\sspeq 7,\ {\cal M}(7)\sspeq \{1,\, 3,\, 5\},\ \tau(3)\sspeq 2$  (see also row No. n\sspeq p\sspeq 7 of the \Modd\,$n$ order table \seqnum{A21630}). Call $\psi_p(h)\sspdef |A_p(h)|$, the number of elements of this set. We have seen that  $\psi_p(h)$ is $0$ or $\varphi(h)$ because for each $h|\delta(p)$ the set $A_p(h)$ is either empty (no $a\sspin {\cal M}(p)$ has \Modd\,$p$ order h) or this set has $\varphi(h)$ elements. Thus $\psi_p(h)\sspleq \varphi(h)$. It is clear that \dstyle{\sum_{h|\delta(h)}\, \psi_p(h)\sspeq \delta(p)} because each $a\sspin {\cal M}(p)$ belongs to one of these sets $A_p(h)$. Now a standard result  is \dstyle{\sum_{h|\delta(h)}\, \varphi(h)\sspeq \delta(p)} (see \eg \cite{Apostol}, {\it Theorem 2.2}, p. 26, wit h$n\sspto \delta(p)$), therefore, \dstyle{\sum_{h|\delta(h)}\, ( \varphi(h) \sspm \psi_p(h))\sspeq 0}, but because $ \varphi(h) \sspm \psi_p(h)\sspgeq 0$ it follows that $ \psi_p(h) \sspeq \varphi(h)\sspgr  1$, not $0$, for each $h|\delta(p)$. This holds especially for $h \sspeq \delta(p)$. Hence the number of primitive \Modd\,$p$ roots is $ \psi_p(\delta(p)) \sspeq \varphi(\delta(p)) \sspgr 1$, and the proof is complete. \hskip 14cm $\square$\psn
For the sequence $\delta(p(n))\, n\sspgeq 1,$ see \seqnum{A130290}. For the smallest positive primitive \Modd\,$n$ roots see the sequence \seqnum{A206550}. Here for prime $n$.\psn
{\bf Corollary 8: Cyclic Modd\,$\bf p$ group}\psn
The multiplicative (Abelian) \Modd\,$p$ group, for $p$ a prime, is the cyclic group $Z_{\delta(p)}$.\pbn
We now concentrate on those groups \Modd\,$n$ which are cyclic and notice that whenever $\delta(n)$ is even there exists a unique smallest positive odd number $>1$ which solves the congruence $x^2 \sspequiv 1\,$ (\Modd\,$n$). This is the nontrivial solution of this congruence. The trivial one is $x\sspeq 1$ (standing for the class of solutions $1\,$(\Modd\,$n$) which also includes $-1$). We prove this first for even $n$.\psn
{\bf Proposition 16: \ Nontrivial square-root of $\bf 1$\,(Modd\,$\bf 2\,k$), $\bf k\sspgeq 2$}\psn
If $n\sspeq 2\,k$, for $k\sspin \{2,3,4,...\}$, and if the {\sl Galois} group ${\cal G}{\sl al}(\mathbb Q(\rho(n))/\mathbb Q)$ is cyclic,  a unique smallest positive solution $r\sspgr 1$ of the congruence  $x^2 \sspequiv 1\,$(\Modd\,$2\,k$) exists, and it is $r\sspeq n-1$.\psn
{\bf Proof:} For even $n\sspeq 2\,k$, $k\sspgeq 2$, the cyclic group is then $Z_{\delta(2\,k)}$ and $\delta(2\,k)$ is even because with the standard prime factorization \dstyle{\delta(2^e\,p_1^{e_1} \cdots p_N^{e_N})\sspeq 2^{e-1}\,\prod_{j=1..N}\,p_j^{e_j-1}\,(p_j-1)}, and if $N=0$ then $e\sspgeq 2$, and  if $N\sspgeq 1$ then $e\sspgeq 1$. In both cases $\delta$ is even, because in the second case there is at least one $j$ and $p_j-1$ is even. Even order $\delta(2\,k)$ of the cyclic group means that powers of the generator $c$ of the $\delta$-cycle, the smallest positive primitive root (of $1$)\, $\Moddn{2\,k}$, come in pairs which are mutually inverse, except for $1\sspeq c^{\delta(2\,k)}$ and \dstyle{s\sspdef c^{\frac{\delta(2\,k)} {2}}}. This  $s\sspgr 1$ is unique and solves of the congruence  $x^2 \sspequiv 1\,($\Modd\,$2\,k$). Because \dstyle{\floor{\frac{(n-1)^2}{n}}\sspeq n\sspm 2 \sspp 0} is, with even $n$, even, and in order to compute 
$(n-1)^2\,$(\Modd\,$n$)  one has to determine $(n-1)^2\, (mod\,n)$ , which is $1$. Therefore, the unique $s$ for even $n$ is $s\sspeq n-1$. \hskip 15.5cm $\square$\psn   
{\bf Example 19:\  Nontrivial square-root of $\bf 1 \,(Modd\, 8)$}\psn
$s(8)\sspeq 7$: $7^2\sspeq 49 \sspequiv  1\,$(\Modd\,$8$), because \dstyle{\floor{\frac{49}{8}}\sspeq 6}, hence $p_8(49)\sspeq +1$, and $49 \sspequiv +1\, (mod \, 8)$, therefore  $49 \sspequiv +1\,$(\Modd\,$8$).\psn
Similarly, for $n\sspeq p$, $p$ an odd prime, the cycle length is \dstyle{\delta(p)\sspeq \frac{p-1}{2}}. In the event that \dstyle{\frac{p-1}{2}} is even, \ie if $p$ is of the form $4\,k+1$, \ie $p\sspequiv 1\,(mod\, 4)$, (see \seqnum{A002144}) a smallest positive nontrivial solution of the quadratic congruence $x^2\sspequiv 1 \,$(\Modd\,$n$) exists. We will numerate the  primes congruent $1\, (mod\, 4)$ by defining $\hat p(n)\sspdef$\seqnum{A002144}$(n)$, $n\sspgeq 1$. \Eg for $n\sspeq \hat p(2)\sspeq 13$ one has, besides the trivial solution $x\sspeq 1\, \Moddn{13}$ (this residue class includes $-1$),  $x\sspequiv 5\, \Moddn{13}$ as a solution. This means that in the group table the row for the element $5$ has a $1$ in the diagonal, like in the row for the  identity element $1$.  These elements are self-inverse like the identity element $1$. Note that in the ordinary reduced $mod\, n$ case there is always a $1$ in the diagonal for row $p-1$, besides the $1$ in the row for  the identity element,  but in that case  $p-1\sspequiv -1\, (mod\, n)$, and hence is a trivial solution. See \seqnum{A206549} for these non-trivial solutions in our case. The other inverses in these cyclic groups can be read off the cycle structure. One pairs the first entry, the generator $c$ of the cycle, with the second to last (the one before the $1$, the second with the third to last, etc. \Eg for $n=13$ one has the inverses $7^{-1}\sspequiv 11\,$(\Modd\,$13),\  11^{-1}\sspequiv 7\,($\Modd\,$13$), and $3^{-1}\sspequiv 9\,$(\Modd\,$13$) , $9^{-1}\sspequiv 3\,$(\Modd\,$13$), and the left over $5$ is self-inverse, as mentioned above. This pairing of inverses is clear from the cyclic structure. \psn
The next objective is to find some method to compute this uniquely existing smallest positive non-trivial solution $s(p)$ of the congruence $x^2\sspequiv 1$\,(\Modd\,$p$) for given prime $p\sspeq \hat p(n)\sspeq$\seqnum{A002144}$ (n)$ for the $k(n)\sspdef \frac{\hat p(n) \sspm 1}{4}$ values given in \seqnum{A005098}. This $k$-sequence starts with $[1, 3, 4, 7, 9, 10, 13, 15, 18, ...]$.  The following  $l$-algorithm  will determine all these non-trivial solutions $s(\hat p(n))$ of the congruence $x^2\sspequiv 1\,$(\Modd\,$\hat p(n)$).  We start with the {\it Ansatz} $s(p)\sspeq \sqrt{(o+1)\,p\sspm 1}$, with some odd number $o$. This means that we are looking for primes of the form \dstyle{\frac{s^2\sspp 1}{o+1}}, and we are interested only in primes of the form $4\,k\sspp 1$, \ie $\hat p(n)$. Because \dstyle{s(\hat p)^2\sspeq (o+1)\,\hat p\sspm 1}, $s(\hat p)$ is odd, say $2\,K\sspp 1$. This implies \dstyle{2\,K\,(K+1)\sspp 1\sspeq {\frac{o+1}{2}}\,\hat p}, which shows that $ {\frac{o+1}{2}}$ has to be odd, \ie $o\sspeq 4\,\hat l\sspp 1$. Now we compute the even number \dstyle{2\, k\sspeq \frac{p-1}{2}\sspeq {\frac{(2\,K+1)^2\sspm (4\,\hat l+1)\,\hat p}{2}}\sspeq 4\,T(K)\sspm 2\,\hat l \sspm 2\,k\,(4\,\hat l+1)}, where we have introduced the triangular numbers $T(K)\sspdef \frac{K\,(K+1)}{2}$, given in \seqnum{A000217}. This implies that  \dstyle{2\,k\,(2\,\hat l\sspp 1) \sspeq 2\,T(K) \sspm \hat l} . Therefore, $\hat l$ is even, say $2\,l$, and we have, with inputs $k,\ l$ and $K$:\ $ p\sspeq 4\,k\sspp 1$,  $s^2(\hat p)\sspeq (o\sspp 1)\,p\sspm 1$,  $o\sspeq 8\,l \sspp 1$  and $s(\hat p)\sspeq 2\,K(\hat p)\sspp 1$\ . Because now $ 2\,T(K)\sspeq 2\,k \sspp 2l\,(4\,k\sspp 1)\sspeq 2\, l\sspp 2k\,(4\,l\sspp 1)$, we find the `nice equation' 
\Beq
4\,T(K)\sspp 1\sspeq (4\,k\sspp 1)\,(4\,l\sspp 1)\ .
\Eeq
For given prime $\hat p$ of the form $4\,k\sspp 1$, i.e. for $k\sspeq \frac{p\sspm 1}{4}$ (see \seqnum{A005098}), we determine the minimal $l\sspin \mathbb N_0$, such that $k\,(1\sspp 4\,l)\sspp l$ produces a triangular number, namely $T(K)$ . The index of this triangular number is then $K\sspeq \frac{\sqrt{8\,T(K)\sspp 1}\sspm 1 }{2}$. We have uses Maple \cite{Maple} to compute these minimal $l$ values. It will be seen from the later {\sl proposition 15} that there is always such a minimal $l$, and not all non-negative values appear. Precisely the entries of \seqnum{A094178} and $0$ occur as minimal $l$ values. We will from now on denote this minimal $l$ also by $l$ (hoping to cause no confusion).  If  $l\sspeq 0$  then  $k$ is a triangular number, and conversely. The non-vanishing $l$ values are characterized by having as elements of the squarefree kernel set of $4\,l+1$ only distinct primes congruent $1\, (mod\, 4)$. \Eg $l=2$ will not appear because the squarefree kernel of $9$ consists of the prime $3$ which is not  congruent $1\, (mod\, 4)$. In order to proof that only numbers from \seqnum{A094178} can appear as positive $l$ values, we use the `nice equation' and observe that $4\,T(K)\sspp 1\sspeq K^2\sspp (K+1)^2$, \ie it is the sum to two neighboring squares with $gcd(K,K+1)\sspeq 1$ (shown by indirect proof, see the remark following {\it definition 3} ). Now a theorem ensures that all prime factors of $4\,T(K)\sspp 1$ are congruent $1\,(mod\, 4)$. See the theorem 3.20 on p. 164 of the reference \cite{NZM}. From the `nice equation' we have either $4\,l\sspp 1\sspeq 1$, \ie $l=0$, or, if $l$ is non-vanishing, then the squarefree kernel set of $4\,l+1$ has only primes congruent $1\, (mod\ 4)$. Therefore $l$ can be $0$ or it is from \seqnum{A094178}, \ie from $[0, 1, 3, 4, 6, 7, 9, 10, 13, 15, 16, 18, 21, 22, 24, 25, 27, 28, ...]$. The $l\sspeq 0$ value appears for the primes $[5, 13, 41, 61, 113, 181, 313, 421, ...]$, which is the sequence \seqnum{A027862}. In this case $k\sspeq \frac{\hat p-1}{4}$ is a triangular number. $l\sspeq1$ for the primes $[17, 29, 53, 73, 109, 137, 281, 397, 449, 593, 757,...]$ This is the sequence \seqnum{A207337}$(n)$, $n\sspgeq 2$. These are the primes $\hat p$ such that $5\,\hat p\sspeq 4\,T(K) \sspp 1$ for some $K\sspeq K(\hat p)$. \Eg $5\,17\sspeq 85  \sspeq 4\,21 \sspeq 1$, thus $K(17)\sspeq 6$. These are also the primes of the form $(m^2 +1)/10$. Just use $m\sspeq m(K)\sspeq 2\ K\sspp 1$. In general, \dstyle{k\sspeq \frac{T(K)\sspm l}{4\,l \sspp 1}}, (see eq. $(90)$ and \dstyle{K\sspeq \frac{\sqrt{2\,(4\,l\sspp 1)\,\hat p\sspp 1}\sspm 1}{2}}. \Eg $\hat p=53$ belongs to $l\sspeq 1$ and $K\sspeq 11$, and $T(K)\sspeq 66$.  This checks with the `nice equation'. We do not give here the tables of the sequences  $l(n) \sspeq l(\hat p(n))$ and  $K(n)\sspeq K(l(n))$ (this incorrect renaming of arguments should not lead to confusion) for given $\hat p(n)$ because after the next {\sl proposition 17} we will find another way to obtain the $K(n)$ numbers directly, hence the corresponding $l(n)$ numbers, and from these the desired $s(\hat p(n))$ solutions from $s(\hat p(n))\sspeq 2\,K(n)\sspp 1$.
\psn
{\bf Proposition 17:\ Congruence $\bf  4\,T(X)\sspp 1\sspeq X^2\sspp (X+1)^2\sspequiv 0\, (mod\,p)$}\psn
There are precisely two incongruent solutions of the congruence  \dstyle{f(X)\sspdef 2\,X^2 \sspp 2\, X\sspp 1\sspequiv 0\, (mod\, p)}, provided $p=\hat p(n)\sspeq$\seqnum{A002144}$(n)$, $n\sspgeq 1$. The smallest positive representative will be called $K(n)$, and the next larger incongruent one is then $K2(n)\sspdef \hat p(n)\sspm 1\sspm K(n)$. \psn
{\bf Proof:} Because of the degree $2$ of $f$ this congruence has at most two incongruent solutions. We shall see that in fact there are two. This congruence is reduced to a problem of quadratic residues, following a standard prescription (see, \eg \cite{Nagell}, pp. 132-3). The discriminant of $f$ is $D\sspeq 2^2\sspm 4\,2\,1\sspeq -4$. Multiplying $f$ by $8$ yields $(4\,X+2)^2\sspp 4\sspequiv 0\,(mod\, 8\,p)$. With $Y\sspeq 4\,X+2$ this is $Y^2\sspeq D\, (mod\, n)$, with a composite modulus $n\sspeq 2^3\, p$. This is a quadratic residue problem,  but $gcd(D,n)$ is not $1$, but 4. Theorem 77 of \cite{Nagell} is applied with $d \sspeq 4$, $e\sspeq 2$, $f\sspeq 1$, $a+1\sspeq -1$ and $n_1\sspeq 2\,p$ Thus the problem is reduced to the quadratic residue problem $Z^2\sspequiv -1\, (mod\, 2\,p)$. This is solved by studying the congruences for the powers of primes, here just $2$ and $p$, separately (see \eg \cite{Nagell}, sect. 26, pp. 83-5). The congruence modulo $2$ has only the solution $+1$ (because $-1 \sspequiv +1\, (mod \,2))$, and modulo $p$ one can consult the {\sl Legendre} symbol \dstyle{\left( \frac{-1}{p}\right)} which is \dstyle{(-1)^{\frac{p-1}{2}}} (see \eg \cite{NZM}, {\it Theorem} $3.2\, (1)$, p. 132). Therefore $-1$ is a quadratic residue modulo $p$ if and only if this symbol is $+1$, demanding that $p\sspequiv 1\,(mod\, 4)$, \ie $p\sspeq \hat p$ from \seqnum{A002144}. Call the smallest positive solution $x_0$, then  $\hat p\sspm x_0$ is also an incongruent solution modulo $\hat p$, and two is the maximal number of solutions because of the degree $2$ of this congruence. This implies that there are $1\sspcdot 2\sspeq 2$ incongruent solutions of this congruence modulo $2\,p$ (see \eg \cite{Nagell}, Theorem 46, p. 84). Returning to the original problem this proves that there are also two incongruent solutions. If the smallest positive solution for $\hat p(n)\sspeq $\seqnum{A207337}$(n)$ is called $K(n)$, then the next larger incongruent solution of $ 4\,T(X)\sspp 1\sspeq X^2\sspp (X\sspp 1)^2)\sspequiv 0\, (mod\,\hat p(n)) $ is $K2(n)\sspdef \hat p(n)\sspm 1\sspm K(n)$, which is obvious. 
\hskip 16cm $\square$
\psn
The pair of sequences of all positive solutions modulo $5, 13,$ and $17$ are given in \seqnum{A047219}, \seqnum{A212160} and \seqnum{A212161}, respectively, where in each case  the even indexed members are the positive solutions congruent to $K(n)$ and the odd indexed ones are the positive solutions  congruent to $K2(n)$, \ie $a(2\,k)\sspeq k\,\hat p(n) \sspp K(n)$ and  $a(2\,k+1)\sspeq k\,\hat p(n) \sspp K2(n)$, $k\sspgeq 0$. For the three given examples $n\sspeq 1,\,2,\,3$\, respectively. The companions  $(K(n),K2(n))$ have been computed with Maple \cite{Maple} and they can be found as \seqnum{A212353} and \seqnum{A212354}. The first entries are for $K$: $[1, 2, 6, 8, 15, 4, 11, 5, 13, 27, 37, 45, 16, 7, 18, 52, 64,...]$ and for K2: $[3, 10, 10, 20, 21, 36, 41, 55, 59, 61, 59, 55, 92, 105, 118, 96, 92, 126,...]$ corresponding to the primes $\hat p$: $[5, 13, 17, 29, 37, 41, 53, 61, 73, 89, 97, 101, 109, 113, 137, 149, 157, 173,...]$. Note that {\it proposition 17} yields directly the searched for nontrivial solution of $s(\hat p(n))^2\sspequiv +1\,$ (\Modd\, $\hat p(n))$ {\it via}  $s(\hat p(n))$ $\sspeq 2\,K(n)\sspp 1$. Compare this with the $s$ values given in \seqnum{A206549} with first entries $[3, 5, 13, 17, 31, 9, 23, 11, 27, 55, 75, 91, 33, 15, 37, 105, 129, 93,... ]$. Because $\hat s(n)\sspdef s(\hat p(n))\sspeq 2\ K(n)\sspp 1 \sspleq \hat p(n)\sspm 2$ , \ie \dstyle{K(n) \sspleq \frac{\hat p\sspm 3}{2}} iff ${\widehat {s2}}(n)\sspdef K2(n)\sspp 1\sspgeq \hat p(n)$ and $\hat s(n)\sspleq {\widehat{s2}}(n)$, $n\sspgeq 1$. Thus only $\hat s(n)\sspin {\cal M}(\hat p)$, the restricted odd residue class \Modd\, $\hat p$, and the solution ${\widehat {s2}}(n)$ is discarded.\psn
This {\it proposition} and the `nice equation' eq. $(90)$ show that the $l-$algorithm from above will indeed produce a solution $l$, related to the existing $K(n)$ for given $\hat p(n)$, {\it via} \dstyle{l(n)\sspeq \left({\frac{4\,T(K(n))\sspp 1}{\hat p(n)}} \sspm 1\right)/4}. The smallest positive non-trivial solution $s(\hat p(n))$ of the congruence $x^2\sspequiv +1\, (Modd\,\hat p(n))$ is then $2\,K(n)\sspp 1$. See $K(n)\sspeq$\seqnum{A212353}$(n)$, $n\sspgeq 1$. \Eg $n=5,\ \hat p(5)\sspeq 37,\ K(5)\sspeq 15$, \dstyle{l\sspeq \floor{\frac{4\,120,\sspp 1}{37}} \sspm 1)/4\sspeq 3}, $\hat s(n)\sspeq  31$, with \dstyle{\floor{\frac{31^2}{37}}\sspeq 25\sspeq o } and $31^2\sspeq 961\sspeq -1\, (mod\, 37)\sspeq 36$.
\psn
 As a application involving $s(\hat p(n))$ we derive the analog of {\sl Wilson}'s theorem $\prod\, {\cal R}(p)\sspeq (p\sspm 1)! \sspequiv (p-1)\, (mod\, p)\sspequiv -1\,  (mod\, p)$, for each prime $p$ (see \eg \cite{Apostol}, Theorem 5.24, p. 116, or \cite{HardyWright}, Theorem 80, p. 68), where ${\cal R}(n)$is the set of the representatives of the smallest positive reduced residue system $mod\,n$ for $n\sspgeq 2$ which is of order $\varphi(n)$. In the \Modd\,$p$ case we have to replace the set ${\cal R}(n)$ by the set ${\cal M}(n)$ of eq. $(65)$ with order $\delta(n)$.
\psn
{\bf Proposition 18:  Analog of  Wilson's theorem for Modd\, ${\bf p}$}\psn
\Beq
\prod\,{\cal M}(p)\sspeq (p-2)!!\sspequiv  {\Caseszwei{$1$}{\text {if} $\frac{p-1}{2}$ \ \text{is odd}}{\hat s(n)}{\text {if} $\frac{p-1}{2}$ \ \text{is even},\  $p\sspeq \hat p(n),\  n\sspin \mathbb N$}}\, (\text{\Modd}\, p)\, ,
\Eeq
where $\hat s(n)\sspeq s(\hat p(n))$, $\hat p(n)\sspeq$\seqnum{A002144}$(n)$,  stands for the above treated nontrivial root \Modd\,$n$, \ie the solution of the congruence $ x^2\sspequiv 1\, \Moddn{$n$}$ which is not $1$, if it exists. Note that $-1\sspequiv +1\, \Moddn{$n$}$, $n\sspgeq 2$ (see {\it lemma 10}). For the double factorials $(p-2)!!$ see \seqnum{A207332}. \psn
{\bf Proof}: The multiplicative group group \Modd\,$p$, p a prime, is the cyclic group $Z_{\delta(p)}$ from {\it proposition 14}. The order $\delta(p)$ of this group is even if and only if $p\sspeq \hat p$, \ie a prime $1\,(mod\, 4)$. In this case we know from above that the non-trivial $\hat s(n)\sspeq s(\hat p(n))$ exists, and an algorithm for finding it has been given. Besides the unit element $1$ and this $\hat s(n)$ all other factors in $\prod\,{\cal M}(\hat p)$ can be paired such that their product is $1$\, (\Modd\,$\hat p$) (see the discussion after {\it proposition 16}), leaving only $s(\hat p)$. In the other case, when $p\sspequiv 3\,(mod\, p)$, the group order is odd. Then  all numbers besides $1$ can be paired in the product to produce $1\, \Moddn{$p$}$, and the result is therefore $1$.\hskip 9cm $\square$\psn
For the cyclic {\sl Galois} groups belonging to $C(n,x)$  the list of the smallest positive primitive roots $r(n)$ (\ie the smallest element from ${\cal{M}}(n)$ which generate these cyclic groups) are found under \seqnum{A206550}. As mentioned above, we do not have a formula for those $n>3$ values whose cycle does not  start with $3$ as smallest positive primitive root, like for $n\sspeq [6, 9, 13, 14, 15, 18, 21, 26, 27, 33, 37, 38, 39, ...]$. \pbn
{\bf Open problems:}\psn
$\bullet$  Proof of the {\it conjecture} on the $q-$sequence related to the discriminant of the minimal $C(n,x)$ polynomials; or find a counterexample.\psn
$\bullet$  Characterization of the values $n$ for which the {\sl Galois} group ${\cal G}_n$ is non-cyclic.\psn
$\bullet$ Characterization of the values $n$ for which the cycle of the cyclic multiplicative group \Modd\,$n\sspiso {\cal G}_n$ is generated by $3$. \psn  
$\bullet$ More theorems on multiplicative \Modd\,$n$ arithmetic.
\pbn \pbn
\hrule
\vfill 
\eject
%%%%%%%%%%%%%%%%%%%%%%%%%%%%%%%
%%%%%%%%%%
\section*{Appendix A}
The proof of eq. $(30)$, with ${\frac{\varphi(n)}{n}}$ replaced by the product-formula given there a bit later, is a (nice) application of {\it PIE} (the principle of inclusion and exclusion. See \eg \cite{Charalambides}, Theorem 4.2, pp. 134 ff). As mentioned above, this formula uses only the distinct prime factors of $n$, the elements of the set $sqfkset(n)$, the set of primes of the squarefree kernel of $n$. The product formula for ${\frac{\varphi(n)}{n}}$ can be read as the generating function for the {\it elementary symmetric functions} (here polynomials) in the variables ${\frac{1}{p_j}}, j\sspeq1,2,...,M(n)$. $M(n)$, also called $\omega(n)$, is given in \cite{Sloane} as \seqnum{A001221}$(n)$.   
\Beq
{\frac{\varphi(n)}{n}} \sspeq \prod_{j=1}^{M(n)}\,\left(1\sspm {\frac{1}{p_j}}\right)\sspeq \sum_{r=0}^{M(n)}\, (-1)^r\, \sigma_r\left({\frac{1}{p_1}},...,{\frac{1}{p_{M(n)}}} \right)\ ,
\Eeq
with $\sigma_0\sspeq 1$, and  symbolically \dstyle{\sigma_r\sspeq \sum_r \,{\frac{1}{p_{.}\, \cdots\, p_{.}}}}, where the sum extends over the \dstyle{{\binomial{M}{r}}} terms with $r$ factors  ${\frac{1}{p_{.}}}$  with increasing indices. \Eg $M=3$ with $\sigma_2\sspeq {\frac{1}{p_1\,p_2}} \sspp {\frac{1}{p_1\,p_3}}\sspp {\frac{1}{p_2\,p_3}}$.\pn
To calculate  \dstyle{\sum_{{k=1}\atop {gcd(k,n)=1}}^{n-1}k} one starts, at the zeroth step ($r=0$), with the unrestricted sum which is  ${\frac{n}{2}}\, (n-1)$, and subtracts, in the next step (r\sspeq 1), the sum of all multiples of each prime dividing $n$ which are $\sspleq (n-1)$. For the  $p_j$-multiples this sum  is \dstyle{p_j\, \sum_{k=1}^{{\frac{ n}{p_j}}\sspm 1}\,k\sspeq {\frac{n}{2}}\,\left ({\frac{n}{p_j}}\sspm 1\right )}, for $j\sspeq 1,2,...,M(n)$. This leads, in step $r\sspeq 1$ to the subtraction of \dstyle{{\frac{n}{2}}\,\sum_{j}\,\left ({\frac{n}{p_{j}}}-1\right )}. Now in this subtraction all multiples of the product of two different $p_j$s appeared twice, therefore one has, in the next step ($r=2$), to add them once. This is done by \dstyle{+1\,\sum_{i\sspkl j}\, p_i\, p_j\,\sum_{k=1}^{ {\frac{n}{p_i\,p_j}}\sspm 1}\ k\sspeq  {\frac{n}{2}}\,\sum_{i\sspkl j}\,\left( {\frac{n}{p_.\,p_.}}\sspm 1\right)}. In step $r=3$  one concentrates on products of three different primes $p_{i_1}\,p_{i_2}\,p_{i_3}$ with $i_1<i_2<i_3$. Now such a $3-$product appeared once in step $r=0$ ( in the unrestricted sum), $-3$ times in step $r=1$, originating from the multiples $(p_{i_1}\,p_{i_2})\cdot p_{i_3}$,  $(p_{i_1}\,p_{i_3})\cdot p_{i_2}$ and  $(p_{i_2}\,p_{i_3})\cdot p_{i_1}$, and $+3$  times in step $r\sspeq 2$, from the multiples  $p_{i_1}\cdot (\,p_{i_2}\,p_{i_3})$, $p_{i_2}\cdot (\,p_{i_1}\,p_{i_3})$, and  $p_{i_3}\cdot (\,p_{i_1}\,p_{i_2})$. Therefore, up to this stage, each such $3-$product appeared once too much, and it is subtracted in this step $r\sspeq 3$ when all these $3-$products terms are summed. Now the pattern starts to become clear. Up to, and including, step $r\sspeq 3$ one has for each $4-$product the counting $1\sspm 4\sspp 6\sspm 4\sspeq -1$. Hence in step $r\sspeq 4 $ one adds once the sum over all multiples of each such $4-$product. In general, in step $r$ this counting will produce $(-1)^{r+1}$ (from the alternating row $r\sspeq 4$ in the {\sl Pascal} triangle \seqnum{A007318}), and one will therefore add $(-1)^r$ times the sum over all multiples of each such $r-$product. This yields
\Beqarray
\sum_{{k=1}\atop {gcd(k,n)=1}}^{n-1}k  \sspeq &&\frac{n}{2}\,\left[ (n-1) - \sum_{i_1}\,\left ({\frac{n}{p_{i_1}}} \sspm 1\right) \sspp \sum_{i_1<i_2}\, \left ({\frac{n}{p_{i_1}\,p_{i_2}}} \sspm 1\right )\ ... -+ \right. \nonumber \\
&& \left. ...\  (-1)^{M(n)}\, \sum_{i_1<...<i_{M(n)}}\left ({\frac{n}{p_{1}\,p_{2}\, \cdots\, p_{M(n)}}}  \sspm 1\right)\right] \ . 
\Eeqarray
Now the $\pm1$ terms all cancel if one takes into account the number of terms of each $r-$ sum, which is ${\binomial{M}{r}}$.  This is due to the vanishing  alternating sum over {\sl Pascal} triangle's row $M$. Thus the result becomes the alternating sum over the elementary symmetric polynomials with the reciprocal distinct prime factors of $n$ as variables.
\Beq
\sum_{{k=1}\atop {gcd(k,n)=1}}^{n-1}k  \sspeq  {\frac{n^2}{2}}\, \sum_{r=0}^{M(n)}\, (-1)^r\, \sigma_r\left({\frac{1}{p_1}}, ...,{\frac{1}{p_{M(n)}}}\right)\ . 
\Eeq
As mentioned above this is written as the generating function of these polynomials and the final result is with the second part of eq. $(92)$
\Beq
\sum_{{k=1}\atop {gcd(k,n)=1}}^{n-1}k  \sspeq  {\frac{n^2}{2}}\, \prod_{j=0}^{M(n)}\, \left(1\sspm {\frac{1}{p_j}}\right)\ .
\Eeq
This shows that \dstyle{s(n)\sspdef {\frac{2}{n^2}}\,\sum_{{k=1}\atop {gcd(k,n)=1}}^{n-1}k} is indeed \dstyle{{\frac{\varphi(n)}{n}}}, which appears in \cite{Sloane} as \seqnum{A076512}$(n)$/\seqnum{A109395}$(n)$.  $s(n)$ is identical for all $n$ with the same distinct prime factors, independent of their multiplicity, \ie it depends only on the elements of $sqfset(n)$. From the multiplicativity of $\varphi$ one sees that this scaled sum $s(n)$ is also multiplicative. \Eg $ s(12) \sspeq s(6)\sspeq s(2)\,s(3)$.
%%%%%%%%%%%%%%%%%%%%%%%%%%%%%%%%%%%%%%%%
\section*{Appendix B}
For the proof of eq. $(56)$ we follow \cite{Remmert}, p.33,  and consider the rewritten version
\Beq
\prod_{k=1}^{n-1}\, \left( 1\sspm e^{2\,\pi\,i\,{\frac{k}{n}}}\right) \sspeq {\Caseszwei{1}{\text{if}\ $n$ \text{is odd\,,}}{0}{\text{if}\ $n$ \text{is even\ .}}}
\Eeq        
In order to see that this is identical to eq. $(56)$, just extract \dstyle{e^{\pi\,i\,{\frac{k}{n}}}}, leading to the ($2\, \cos$) factors under the product,  and in front summing in the exponent leads to the factor \dstyle{e^{{\frac{1}{2}}\,\pi\,i\,(n-1)}\sspeq i^{n-1}\sspeq (-1)^{\frac{n-1}{2}}}\ . \pn
The \lhs of eq. $(96)$ is then computed using (for the last step see \eg \cite{GKP}, Exercise 50, p. 149) 
\Beq
1\sspp z\sspp z^2\sspp\ ...\ \sspp z^{n-1}\sspeq {\frac{z^n\sspm 1}{z\sspm 1}}\sspeq \prod_{k=1}^{n-1}\,(z\sspm e^{2\,\pi\,i\,{\frac{k}{n}}})\ , 
\Eeq
where one specializes to  $z=-1$. (By putting $z\sspeq +1$ one finds the result \dstyle{\prod_{k=1}^{n-1}\,2\,\sin\left(\pi\,{\frac{k}{n}}\right)\sspeq n\, , \ n\sspgeq 2}, with $1$ for $n\sspeq 1$.)
\pbn
\hrule
\vfill 
\eject

\pbn\pbn\pbn
\hrule
\pbn\pbn\pbn 
Keywords: regular $n$-gons, algebraic number, minimal polynomial, zeros, factorization, congruences, {\sl Galois} theory, cycle graphs.\psn 
AMS MSC numbers: 11R32, 11R04, 08B10,  13F20, 12D10, 13P05 \psn
OEIS A-numbers: \seqnum{A000001}, \seqnum{A000005}, \seqnum{A000007}, \seqnum{A000010}, \seqnum{A000012}, \seqnum{A000035}, \seqnum{A000265}, \seqnum{A000668}, \seqnum{A000688}, \seqnum{A001221}, \seqnum{A002110}, \seqnum{A002144}, \seqnum{A003277}, \seqnum{A004124}, \seqnum{A005013}, \seqnum{A005098}, \seqnum{A005117}, \seqnum{A006053}, \seqnum{A006054},  \seqnum{A007310}, \seqnum{A007318}, \seqnum{A007775}, \seqnum{A007814}, \seqnum{A007947}, \seqnum{A007955}, \seqnum{A008683}, \seqnum{A023022},  \seqnum{A024556}, \seqnum{A027862}, \seqnum{A033949}, \seqnum{A045572}, \seqnum{A049310}, \seqnum{A0503384}, \seqnum{A052547}, \seqnum{A047219}, \seqnum{A053120}, \seqnum{A055034}, \seqnum{A066651}, \seqnum{A076512},  \seqnum{A077998}, \seqnum{A084101}, \seqnum{A085810}, \seqnum{A088520}, \seqnum{A090298}, \seqnum{A092260}, \seqnum{A094178}, \seqnum{A106803}, \seqnum{A109395}, \seqnum{A110551}, \seqnum{A113807}, \seqnum{A116423}, \seqnum{A120757}, \seqnum{A127672}, \seqnum{A128672}, \seqnum{A130290}, \seqnum{A130777}, \seqnum{A147600},  \seqnum{A162699}, \seqnum{A181875}, \seqnum{A181876}, \seqnum{A181878}, \seqnum{A181879}, \seqnum{A181880}, \seqnum{A193376}, \seqnum{A193377}, \seqnum{A193679}, \seqnum{A193680}, \seqnum{A193681},   \seqnum{A193682},  \seqnum{A203571}, \seqnum{A203572}, \seqnum{A203575}, \seqnum{A204453},  \seqnum{A204454},  \seqnum{A204457}, \seqnum{A204458}, \seqnum{A206543}, \seqnum{A206544}, \seqnum{A206545}, \seqnum{A206546}, \seqnum{A206547}, \seqnum{A206548}, \seqnum{A206549}, \seqnum{A206550}, \seqnum{A206551}, \seqnum{A206552}, \seqnum{A207333}, \seqnum{A207334}, \seqnum{A207337}, \seqnum{A210845}, \seqnum{A212160}, \seqnum{A212161}, \seqnum{A212353}, \seqnum{A212354}, \seqnum{A215041}, \seqnum{A215046}, \seqnum{A216319}, \seqnum{A216320}, \seqnum{A216322}, \seqnum{A216326}, \seqnum{A255237}. \seqnum{A282624}.
\vfill 
\eject
\noindent
%\begin{landscape}
\begin{center}
{\large {\bf Table 1: Reduced DSR-algebras (over $\boldsymbol {\mathbb Q}$), $\bf n=3,...,12$. }} 
\end {center}
\begin{center}  
%{\large  $\bf{ }$ }
%\end{center}
%\begin{center}
\begin{tabular}{|l|l|l|l|c|}\hline
&& &&\\
\hskip .2cm$\bf n$ & $\bf \boldsymbol{\rho}\sspequiv \boldsymbol{\rho(n)} $  &\bf{ reduced DSR-algebra}   & $\bf\boldsymbol{\delta}(n)$  & $\bf DSR-basis$ \\
&& && \\ \hline\hline
&& && \\
\hskip .2cm$\bf 3$  &   1    &  $ \rho^2\sspeq 1$ &  $1$      & $<1>$  \\
&& &&  \\
\hskip .2cm$\bf 4$  &   $\sqrt{2}$    & $\rho^2\sspeq 2$ &   $2$     &$ <1,\rho>$                \\
&& &&\\
\hskip .2cm$\bf 5$  &  \dstyle{\varphi\sspeq \frac{1}{2}\,(1\sspp \sqrt{5})}     &   $\rho^2\sspeq \rho\sspp 1$ &     $2$  &  $<1,\rho>$       \\
&& && \\
\hskip .1cm\fbox{{\color{Bittersweet} ${\bf 6}$}}  &  $\sqrt{3}$     &   $\rho^2\sspeq 3,\ \ [\sigma\sspeq \rho^2\sspm 1\sspeq 2]$ &    $2$    &    $<1,\rho>$         \\
&& && \\
\hskip .2cm$\bf 7$  &  $2\,\cos\left( {\frac{\pi}{7}}\right)$ &   $\rho^2\sspeq 1\sspp \sigma,\ \ \sigma^2\sspeq 1\sspp \rho\sspp \sigma ,\ \ \rho\,\sigma\sspeq \rho\sspp \sigma$  &   $3$  &  $<1,\rho,\sigma>$   \\
&& && \\
\hskip .2cm$\bf 8$  &   $\sqrt{2\sspp \sqrt{2}\,}$   &   $\rho^2\sspeq 1\sspp \sigma,\ \sigma^2\sspeq 2\,\sigma\sspp 1,\ \ \tau^2\sspeq 2\,(1\sspp \sigma),$   &   $4$  &  $<1,\rho,\sigma,\tau>$       \\ 
   &       & $\rho\,\sigma\sspeq \rho\sspp \tau,\ \ \rho\,\tau\sspeq 2\,\sigma,\ \ \sigma\, \tau\sspeq 2\,\rho\sspp \tau$,  &      &    \\
   &       &  $(\sigma\sspeq 1\sspp \sqrt{2},\ \ \tau\sspeq \sqrt{2}\,\rho)$ &      &    \\
&& &&  \\
\hskip .1cm\fbox{{\color{Bittersweet} ${\bf 9}$}}  &  $2\,\cos\left( {\frac{\pi}{9}}\right)$   &  $ \rho^2\sspeq 1\sspp \sigma,\ \ \sigma^2\sspeq 2\sspp \rho\sspp \sigma,\ \ \rho\,\sigma\sspeq 2\,\rho\sspp 1, $ &   $3$    & $<1,\rho,\sigma>$     \\
   &       & $ (\sigma\sspeq \rho^2\sspm 1), \ \  [\tau\sspeq \rho\,(\sigma\sspm 1)\sspeq 1\sspp \rho]$  &      &    \\
&& && \\
\fbox{{\color{Bittersweet} ${\bf 10}$}} &   $\varphi\,\sqrt{3\sspm \varphi}$   & $ \rho^2\sspeq 1\sspp \sigma,\ \ \sigma^2\sspeq -1\sspp 3\,\sigma,\ \ \tau^2\sspeq-1\sspp4\, \sigma,\ \ $ &   $4$   &  $<1,\rho,\sigma,\tau>$     \\
   &       & $ \rho\,\sigma\sspeq \rho\sspp \tau, \ \ \rho\,\tau\sspeq 3\,\sigma\sspm 2,\ \ \sigma\,\tau\sspeq \rho\sspp 2\,\tau$,  &      &    \\
   &   & $ (\sigma\sspeq \rho^2\sspm 1\sspeq 1\sspp \varphi, \ \ \tau \sspeq \rho\,(\sigma\sspm 1)$, &   &    \\
 &       & \hskip .5cm $\sspeq (1\sspp \varphi)\,\sqrt{3\sspm \varphi}) $  &      &    \\
 &       & $[\omega\sspeq 2\,(-1\sspp \sigma)\sspeq 2\,\varphi]$  &      &    \\
&& && \\
\hskip .1cm$\bf 11$  & $2\,\cos\left( {\frac{\pi}{11}}\right)$  &$ \rho^2\sspeq 1\sspp \sigma,\ \ \sigma^2\sspeq 1\sspp \sigma\sspp \omega, $  & $5$ & $<1,\rho,\sigma,\tau,\omega>$\\
 &       & $ \tau^2\sspeq 1\sspp \sigma\sspp\tau\sspp \omega, \ \ \omega^2\sspeq 1\sspp \rho\sspp\sigma\sspp \tau\sspp \omega$, &      &    \\
 &       & $\rho\,\sigma\sspeq \rho\sspp \tau,\ \ \rho\, \tau\sspeq \sigma\sspp \omega,\ \ \rho\,\omega\sspeq \tau\sspp \omega$, &      &    \\
&       & $\sigma\,\tau\sspeq \rho\sspp \tau\sspp \omega,\ \ \sigma\,\omega\sspeq \sigma\sspp\tau\sspp\omega$,  &      &    \\
&       & $\tau\,\omega\sspeq \rho\sspp \sigma\sspp \tau\sspp \omega$,  &      &    \\
&       & $(\sigma\sspeq \rho^2\sspm 1, \ \ \tau \sspeq \rho\,(\sigma\sspm 1),\ \ \omega\sspeq \sigma\,(\sigma\sspm 1)\sspm 1) $  &      &    \\
&& &&\\
\fbox{{\color{Bittersweet} ${\bf 12}$}}  &  $\sqrt{2\sspp \sqrt{3}\,}$  &  $\rho^2\sspeq 1\sspp \sigma,\ \ \sigma^2\sspeq 2\,(1\sspp \tau),\ \ \tau^2\sspeq 3\,(1\sspp \sigma),\ \ $  &  $4$   &  $<1,\rho,\sigma,\tau>$  \\
&       & $\rho\,\sigma\sspeq \rho\sspp \tau ,\ \ \rho\, \tau\sspeq 1\sspp 2\,\sigma, 
\ \  \sigma\,\tau\sspeq 3\,\rho\sspp \tau$, &      &    \\
&       &$[\omega \sspeq 1\sspp \sigma,\ \ \Chi\sspeq 2\,\rho] $,&      &    \\
&       & $(\sigma\sspeq \rho^2\sspm 1 \sspeq 1\sspp \sqrt{3}$, \ \ \dstyle{\tau \sspeq {\frac{\sqrt{2}}{2}}\,(3\sspp \sqrt{3}\,))} &      &    \\
$\hskip .2cm\bf \vdots$&& && \\
\hline
\end{tabular}
\end{center}
\psn
\psn
\begin{center}
\dstyle{\rho(n)\sspdef 2\,\cos\left( {\frac{\pi}{n}}\right),\ \ R^{(n)}_k\sspeq S(k-1,\rho(n)),\ \ k\sspeq 1, ..., \floor{\frac{n}{2}}, }\psn
$R_1\sspeq 1,\ R_2\sspequiv \rho,\   R_3\sspequiv \sigma,\ R_4\sspequiv \tau,\ R_5\sspequiv \omega,\ R_6\sspequiv \Chi$ (dependence on $n$ suppressed).
\end{center}\psn
In round brackets the values for the basis elements are given in terms of $\rho$. In square brackets the linear dependent DSRs are given. Boxed $n$-numbers indicate linear dependent DSRs. 
%\end{landscape}
\vfill
\eject
\noindent
\begin{center}
{\large {\bf Table 2: Minimal polynomials of \dstyle{\bf 2\,cos\left ({\frac{\boldsymbol{\pi} }{n}}\right)} }} {\bf for}\  $\bf n\sspeq 1,2,...,30$.
\end {center}
\begin{center}  
\begin{tabular}{|l|l|}\hline
\bf n & \hskip 3cm $\bf C(n,x)$ \\ 
\hline\hline
$\bf 1$    &  $ x+2 $   \\        
\hline
$\bf 2$    &  $ x$      \\ 
\hline
$\bf 3$    &  $ x-1 $    \\    
\hline
$\bf 4$    &  $x^2 -2$  \\
\hline
$\bf 5$    & $x^2-x-1  $ \\    
\hline
$\bf 6$    & $ x^2-3 $      \\
\hline
$\bf 7$    & $x^3-x^2-2\,x+1  $  \\     
\hline
$\bf 8$    &  $x^4-4\,x^2+2  $   \\
\hline
$\bf 9$    &$x^3-3x-1  $ \\    
\hline
$\bf 10$  & $ x^4-5\,x^2+5 $   \\
\hline
$\bf 11$  &$x^5-x^4-4\,x^3+3\,x^2+3\,x-1   $ \\  
\hline
$\bf 12$  &  $ x^4-4\,x^2+1  $  \\            
\hline
$\bf 13$  & $ x^6-x^5-5\,x^4+4\,x^3+6\,x^2-3\,x-1 $ \\       
\hline
$\bf 14$  & $x^6-7\,x^4+14\,x^2-7   $        \\
\hline
$\bf 15$  &  $x^4+x^3-4\,x^2-4\,x+1     $\\  
\hline
$\bf 16$  & $ x^8-8\,x^6+20\,x^4-16\,x^2+2  $ \\
\hline
$\bf 17$  & $x^8-x^7-7\,x^6+6\,x^5+15\,x^4-10\,x^3-10\,x^2+4\,x+1   $  \\  
\hline
$\bf 18$  &  $ x^6-6\,x^4+9\,x^2-3  $ \\    
\hline
$\bf 19$  & $ x^9-x^8-8\,x^7+7\,x^6+21\,x^5-15\,x^4-20\,x^3+10\,x^2+5\,x-1 $ \\     
\hline
$\bf 20$  &  $ x^8-8\,x^6+19\,x^4-12\,x^2+1  $ \\   
\hline
$\bf 21$  &  $  x^6+x^5-6\,x^4-6\,x^3+8\,x^2+8\,x+1 $   \\
\hline
$\bf 22$  & $ x^{10}-11\,x^8+44\,x^6-77\,x^4+55\,x^2-11 $ \\
\hline
$\bf 23$  &$ x^{11}-x^{10}-10\,x^9+9\,x^8+36\,x^7-28\,x^6-56\,x^5+35\,x^4+ $ \\
 & $35\,x^3-15\,x^2-6\,x+1,    $ \\ 
\hline
$\bf 24$  &  $x^8-8\,x^6+20\,x^4-16\,x^2+1 $    \\
\hline
$\bf 25$  &  $  x^{10}-10\,x^8+35\,x^6-x^5-50\,x^4+5\,x^3+25\,x^2-5\,x-1  $  \\
\hline
$\bf 26$  & $ x^{12}-13\,x^{10}+65\,x^8-156\,x^6+182\,x^4-91\,x^2+13  $ \\   
\hline
$\bf 27$  & $x^9-9\,x^7+27\,x^5-30\,x^3+9\,x-1   $   \\
\hline
$\bf 28$  & $ x^{12}-12\,x^{10}+53\,x^8-104\,x^6+86\,x^4-24\,x^2+1 $       \\
\hline
$\bf 29$  &  $ x^{14}-x^{13}-13\,x^{12}+12\,x^{11}+66\,x^{10}-55\,x^9-165\,x^8 + $\\
&$ 120\,x^7+210\,x^6-126\,x^5-126\,x^4+56\,x^3+28\,x^2-7\,x-1  $   \\
\hline
$\bf 30$  &  $  x^8-7\,x^6+14\,x^4-8\,x^2+1$  \\ 
\hline 
$\hskip .2cm\bf \vdots$&\\
\hline
\end{tabular}
\end{center}
\vfill
\eject
\noindent
%\begin{landscape}
\begin{center}
{\large {\bf Table 3: \seqnum{A187360}$\bf (n,m)$  coefficient array of \\ minimal polynomials of \dstyle{\bf 2\, cos\left ({\frac{\boldsymbol{\pi} }{n}}\right)}, rising powers}}
\end {center}
\begin{center}  
%{\large  $\bf{ }$ }
%\end{center}
%\begin{center}
\begin{tabular}{|l|r|r|r|r|r|r|r|c|}\hline
&& && &&  &&\\
\bf n/m & $\bf 0$  & $\bf 1$  &$\bf 2$  &$\bf 3$  &$\bf 4$  & $\bf 5$ & $\bf 6$ &...\\
&& && && &&\\ \hline\hline
&& && && &&\\
$\bf 1$  &   2     &  1     &        &        &          &       &   &          \\
&& && &&  && \\
$\bf 2$  &   0   &     1 &        &        &          &       &&        \\
&& && &&  && \\
$\bf 3$  &   -1    &     1 &        &        &          &       &&          \\
&& && &&  && \\
$\bf 4$  &   -2    &      0 &   1     &      &          &       &&            \\
&& && &&  &&\\
$\bf 5$  &  -1     &   -1 &     1  &        &          &       & &        \\
&& && &&  &&\\
$\bf 6$  &   -3    &    0&    1    &       &             &         &    &    \\
&& && &&  &&\\
$\bf 7$  &   1     &   -2  &   -1  &    1   &      &        &     &  \\ 
&& && &&  && \\
$\bf 8$  &  2     &  0  &   -4    &  0     & 1  &  &  &      \\
&& && &&  &&\\
$\bf 9$  &   -1     &    -3 &   0    &    1    &   &  &  &      \\
&& && &&  &&\\
$\bf 10$  &   5     &    0     &   -5   &   0   & 1   &      &  &    \\
&& && &&  &&\\
$\bf 11$  &  -1   &  3  &  3   & -4  &  -1  & 1 &  &   \\
&& && &&  &&\\
$\bf 12$  &  1   &  0  &  -4   & 0 &  1  &  &  &   \\
&& && &&  &&\\
$\bf 13$  &  -1   &  -3  &  6   &  4  &  -5  & -1 &   1 & \\
&& && &&  &&\\
$\bf 14$  &  -7   &  0  &  14   &  0  &  -7  & 0  &   1 & \\
&& && &&  &&\\
$\bf 15$  &  1   &  -4  &  -4   &  1  &  1  &   & &    \\
%% && && &&  &&\\ 
$\hskip .2cm\bf\vdots$&& && && &&\\
\hline
\end{tabular}
\end{center}
%\end{landscape}
\vfill
\eject
\noindent
\begin{center}
{\large {\bf Table 4: Zeros of  $\bf C(n,x)\bf$ in power basis (rising powers of $\bf {\boldsymbol \rho}(n)$)\ {\bf for}\  $\bf n\sspeq 1,2,...,30$.}}
\end {center}
\begin{center}  
\begin{tabular}{|l|l|}\hline
\bf n & \hskip 2cm {\bf  coefficients of $\bf C$-zeros in power basis}\  $\bf <{\boldsymbol \rho}^0, ... ,{\boldsymbol \rho}^{{\boldsymbol\delta}(n)-1}>$ \\ 
\hline\hline
$\bf 1$    &  $[[-2]] $   \\        
\hline
$\bf 2$    &  $[[0]]  $      \\ 
\hline
$\bf 3$    &  $ [[1]] $    \\    
\hline
$\bf 4$    &  $ [[0, 1], [0, -1]]$  \\
\hline
$\bf 5$    & $ [[0, 1], [1, -1]] $ \\    
\hline
$\bf 6$    & $ [[0, 1], [0, -1]] $      \\
\hline
$\bf 7$    & $  [[0, 1], [-1, -1, 1], [2, 0, -1]]  $  \\     
\hline
$\bf 8$    &  $  [[0, 1], [0, -3, 0, 1], [0, 3, 0, -1], [0, -1]]     $   \\
\hline
$\bf 9$    &$  [[0, 1], [-2, -1, 1], [2, 0, -1]]       $ \\    
\hline
$\bf 10$  & $  [[0, 1], [0, -3, 0, 1], [0, 3, 0, -1], [0, -1]]      $   \\
\hline
 $\bf 11$  &$ [[0, 1], [0, -3, 0, 1], [1, 2, -3, -1, 1], [-2, 0, 4, 0, -1], [2, 0, -1]]   $ \\  
\hline
$\bf 12$  &  $[[0, 1], [0, 4, 0, -1], [0, -4, 0, 1], [0, -1]]        $  \\            
\hline
$\bf 13$  & $[[0, 1], [0, -3, 0, 1], [0, 5, 0, -5, 0, 1], [1, -3, -3, 4, 1, -1], [-2, 0, 4, 0, -1], [2, 0, -1]]           $ \\       
\hline
$\bf 14$  & $   [[0, 1], [0, -3, 0, 1], [0, 5, 0, -5, 0, 1], [0, -5, 0, 5, 0, -1], [0, 3, 0, -1], [0, -1]]   $        \\
\hline
$\bf 15$  &  $  [[0, 1], [-2, 3, 1, -1], [-1, -4, 0, 1], [2, 0, -1]]    $\\  
\hline
$\bf 16$  & $  [[0, 1], [0, -3, 0, 1], [0, 5, 0, -5, 0, 1], [0, -7, 0, 14, 0, -7, 0, 1], [0, 7, 0, -14, 0, 7, 0, -1], [0, -5, 0, 5, 0, -1],  $ \\
               &$  [0, 3, 0, -1], [0, -1]]  $   \\
\hline
$\bf 17$  & $  [[0, 1], [0, -3, 0, 1], [0, 5, 0, -5, 0, 1], [0, -7, 0, 14, 0, -7, 0, 1], [-1, 4, 6, -10, -5, 6, 1, -1],  $  \\  
    &$[2, 0, -9, 0, 6, 0, -1], [-2, 0, 4, 0, -1], [2, 0, -1]]    $   \\
\hline
$\bf 18$  &  $ [[0, 1], [0, 5, 0, -5, 0, 1], [0, -4, 0, 5, 0, -1], [0, 4, 0, -5, 0, 1], [0, -5, 0, 5, 0, -1], [0, -1]] $ \\    
\hline
$\bf 19$  & $ [[0, 1], [0, -3, 0, 1], [0, 5, 0, -5, 0, 1], [0, -7, 0, 14, 0, -7, 0, 1], [1, 4, -10, -10, 15, 6, -7, -1, 1],   $ \\  
    &$ [-2, 0, 16, 0, -20, 0, 8, 0, -1], [2, 0, -9, 0, 6, 0, -1], [-2, 0, 4, 0, -1], [2, 0, -1]]   $   \\
\hline
$\bf 20$  & $ [[0, 1], [0, -3, 0, 1], [0, -7, 0, 14, 0, -7, 0, 1], [0, 8, 0, -18, 0, 8, 0, -1], [0, -8, 0, 18, 0, -8, 0, 1], $ \\   
               &$[0, 7, 0, -14, 0, 7, 0, -1], [0, 3, 0, -1], [0, -1]]   $   \\
\hline
 $\bf 21$  &  $  [[0, 1], [0, 5, 0, -5, 0, 1], [2, 3, -4, -1, 1], [-3, -9, 1, 6, 0, -1], [-2, 0, 4, 0, -1], [2, 0, -1]]        $   \\
\hline
$\bf 22$  & $ [[0, 1], [0, -3, 0, 1], [0, 5, 0, -5, 0, 1], [0, -7, 0, 14, 0, -7, 0, 1], [0, 9, 0, -30, 0, 27, 0, -9, 0, 1],  $ \\
               &$   [0, -9, 0, 30, 0, -27, 0, 9, 0, -1], [0, 7, 0, -14, 0, 7, 0, -1], [0, -5, 0, 5, 0, -1], [0, 3, 0, -1], [0, -1]] $   \\
\hline
$\bf 23$  &$  [[0, 1], [0, -3, 0, 1], [0, 5, 0, -5, 0, 1], [0, -7, 0, 14, 0, -7, 0, 1], [0, 9, 0, -30, 0, 27, 0, -9, 0, 1],  $ \\ 
  &$ [-1, -5, 15, 20, -35, -21, 28, 8, -9, -1, 1], [2, 0, -25, 0, 50, 0, -35, 0, 10, 0, -1], $ \\
  &$ [-2, 0, 16, 0, -20, 0, 8, 0, -1], [2, 0, -9, 0, 6, 0, -1], [-2, 0, 4, 0, -1], [2, 0, -1]]  $ \\
\hline
$\bf 24$  &  $   [[0, 1], [0, 5, 0, -5, 0, 1], [0, -7, 0, 14, 0, -7, 0, 1], [0, -8, 0, 6, 0, -1], [0, 8, 0, -6, 0, 1],    $    \\
   &$ [0, 7, 0, -14, 0, 7, 0, -1], [0, -5, 0, 5, 0, -1], [0, -1]]     $ \\
\hline
$\bf 25$  &  $ [[0, 1], [0, -3, 0, 1], [0, -7, 0, 14, 0, -7, 0, 1], [0, 9, 0, -30, 0, 27, 0, -9, 0, 1],     $  \\
  &$ [0, -10, 5, 30, -5, -27, 1, 9, 0, -1], [0, 10, -15, -15, 20, 7, -8, -1, 1], [-2, 0, 16, 0, -20, 0, 8, 0, -1],    $\\
  &$ [2, 0, -9, 0, 6, 0, -1], [-2, 0, 4, 0, -1], [2, 0, -1]]    $  \\
\hline
$\bf 26$  & $ [[0, 1], [0, -3, 0, 1], [0, 5, 0, -5, 0, 1], [0, -7, 0, 14, 0, -7, 0, 1], [0, 9, 0, -30, 0, 27, 0, -9, 0, 1],     $ \\  
  &$ [0, -11, 0, 55, 0, -77, 0, 44, 0, -11, 0, 1], [0, 11, 0, -55, 0, 77, 0, -44, 0, 11, 0, -1],     $  \\ 
   &$ [0, -9, 0, 30, 0, -27, 0, 9, 0, -1], [0, 7, 0, -14, 0, 7, 0, -1], [0, -5, 0, 5, 0, -1], [0, 3, 0, -1], [0, -1]]   $\\
\hline
$\bf 27$  & $  [[0, 1], [0, 5, 0, -5, 0, 1], [0, -7, 0, 14, 0, -7, 0, 1], [-2, 7, 1, -14, 0, 7, 0, -1], [2, -5, -4, 5, 1, -1],  $   \\
&$  [2, -1, -16, 0, 20, 0, -8, 0, 1], [-2, 0, 16, 0, -20, 0, 8, 0, -1], [-2, 0, 4, 0, -1], [2, 0, -1]]        $\\
\hline
$\hskip .2cm\bf\vdots$&  {\bf continued on next page}   \\
\hline
\end{tabular}
\end{center}
\pbn
{\bf Note: the higher power coefficients not shown are all zero}.\psn
{\bf Example: $\bf n=27$:  the first zero of  $\bf C(27,x)$ is  ${\boldsymbol \rho}(27)$, the second zero is $\bf 5\, {\boldsymbol \rho}(27) \sspm 5\,{\boldsymbol \rho}(27)^3\sspm 1\,{\boldsymbol \rho}(27)^5$, {\bf etc}.} 
\vfill
\eject
\begin{center}
{\large {\bf Table 4 continued:  Zeros of  $\bf C(n,x)\bf$ in power basis (rising powers of $\bf {\boldsymbol \rho}(n)$)\ {\bf for}\  $\bf n\sspeq 28, ..., 30$.}}
\end {center}
\begin{center}  
\begin{tabular}{|l|l|}\hline
\bf n & \hskip 2cm {\bf power basis coefficients for}\  $\bf <{\boldsymbol \rho}^0, ... ,{\boldsymbol \rho}^{{\boldsymbol\delta}(n)-1}>$ \\ 
\hline\hline
$\bf 28$  & $ [[0, 1], [0, -3, 0, 1], [0, 5, 0, -5, 0, 1], [0, 9, 0, -30, 0, 27, 0, -9, 0, 1], [0, -11, 0, 55, 0, -77, 0, 44, 0, -11, 0, 1],  $       \\
 &$ [0, 12, 0, -67, 0, 96, 0, -52, 0, 12, 0, -1], [0, -12, 0, 67, 0, -96, 0, 52, 0, -12, 0, 1],  $\\
 &$ [0, 11, 0, -55, 0, 77, 0, -44, 0, 11, 0, -1], [0, -9, 0, 30, 0, -27, 0, 9, 0, -1], [0, -5, 0, 5, 0, -1], [0, 3, 0, -1],  $\\
&$  [0, -1]]    $\\
\hline
$\bf 29$  &  $ [[0, 1], [0, -3, 0, 1], [0, 5, 0, -5, 0, 1], [0, -7, 0, 14, 0, -7, 0, 1], [0, 9, 0, -30, 0, 27, 0, -9, 0, 1],  $\\
              &$ [0, -11, 0, 55, 0, -77, 0, 44, 0, -11, 0, 1], [0, 13, 0, -91, 0, 182, 0, -156, 0, 65, 0, -13, 0, 1], $   \\
              &$[1, -7, -21, 56, 70, -126, -84, 120, 45, -55, -11, 12, 1, -1], [-2, 0, 36, 0, -105, 0, 112, 0, -54, 0, 12, 0, -1],    $   \\
 &$ [2, 0, -25, 0, 50, 0, -35, 0, 10, 0, -1], [-2, 0, 16, 0, -20, 0, 8, 0, -1], [2, 0, -9, 0, 6, 0, -1], [-2, 0, 4, 0, -1],   $   \\ 
  &$ [2, 0, -1]]     $   \\ 
\hline
$\bf 30$  &  $ [[0, 1], [0, -7, 0, 14, 0, -7, 0, 1], [0, -7, 0, 22, 0, -13, 0, 2], [0, 4, 0, -13, 0, 7, 0, -1], [0, -4, 0, 13, 0, -7, 0, 1],  $  \\ 
              &$ [0, 7, 0, -22, 0, 13, 0, -2], [0, 7, 0, -14, 0, 7, 0, -1], [0, -1]]   $   \\
\hline 
$\hskip .3cm\bf\vdots$&\\
\hline
\end{tabular}
\end{center}
\vfill
\eject
\begin{landscape}
\begin{center}
{\large {\bf Table 5: Extended set $\bf \widehat{\mathbfcal M}(n)$, first differences $\bf {\boldsymbol \triangle}\,\widehat{\mathbfcal  M}(n) $, and floor-pattern $\bf{\mathbfcal F}(n)$ for composed odd squarefree modulus $\bf n$.}} 
\end {center}
\begin{center}  
%{\large  $\bf{a }$ }
%\end{center}
%\begin{center}
\begin{tabular}{|c|r|c|c|c|c|}\hline
&& && &\\
$\bf m$ &  $\bf n(m) $ & $\bf 2\,{\boldsymbol \delta}($n\bf $)$   &$\bf \widehat{\mathbfcal M}(n)$ & $\bf \triangle\,\widehat{\mathbfcal  M}(n) $ & $\bf{\mathbfcal F}(n)$ \\
&& && & \\ \hline\hline
$\bf 1$ &   $\bf 15$  & $8$ &  $ [0, 1, 7, 11, 13, 17]  $  &  $[1, 6, 4, 2, 4] $    & $[2, 1, 0, 1, 0, 1, 2, 0]  $  \\
\hline
$\bf 2$  &   $\bf 21$ &$12$  & $ [0, 1, 5, 11, 13, 17, 19, 23]$ &   $[1, 4, 6, 2, 4, 2, 4] $   & $[1, 2, 0, 1, 0, 1, 0, 1, 0, 2, 1, 0]  $  \\
\hline
$\bf 3$  &  $\bf 33 $  &$20  $ &  $ [0, 1, 5, 7, 13, 17, 19, 23, 25, 29, 31, 35]$ &     $[1, 4, 2, 6, 4, 2, 4, 2, 4, 2, 4] $  &  $ [1, 0, 2, 1, 0, 1, 0, 1, 0, 1,    $       \\
&& &  &  & $ 0, 1, 0, 1, 0, 1, 2, 0, 1, 0]  $  \\
\hline
$\bf 4$  &  $\bf 35$ & $ 24$ &  $[0, 1, 3, 9, 11, 13, 17, 19, 23, 27, 29, 31, 33, 37]$ & $[1, 2, 6, 2, 2, 4, 2, 4, 4, 2, 2, 2, 4]$  &  $[0, 2, 0, 0, 1, 0, 1, 1, 0, 0, 0, 1,  $  \\
&$5\sspcdot 7$  & & & & $  0, 0, 0, 1, 1, 0, 1, 0, 0, 2, 0, 0]     $  \\
\hline
$\bf 5$  &   $\bf 39$  & $24$ &   $ [0, 1, 5, 7, 11, 17, 19, 23, 25, 29, 31, 35, 37, 41]$   &   $[1, 4, 2, 4, 6, 2, 4, 2, 4, 2, 4, 2, 4]  $  &  $[1, 0, 1, 2, 0, 1, 0, 1, 0, 1, 0, 1,      $       \\
&$3\sspcdot 13$&  & & & $  0, 1, 0, 1, 0, 1, 0, 2, 1, 0, 1, 0]    $  \\
 \hline
$ \bf 6$  & $\bf 51$ & $32$ & $[0, 1, 5, 7, 11, 13, 19, 23, 25, 29, $  & $ [1, 4, 2, 4, 2, 6, 4, 2, 4, 2, 4, 2, 4, 2, 4, 2, 4] $ & $ [1, 0, 1, 0, 2, 1, 0, 1, 0, 1, 0, 1, 0, 1, 0, 1, $ \\
&$3\sspcdot 17$& &$ 31, 35, 37, 41, 43, 47, 49, 53]$ & & $ 0, 1, 0, 1, 0, 1, 0, 1, 0, 1, 2, 0, 1, 0, 1, 0]  $  \\
 \hline
$\bf 7 $  & $\bf 55 $ &$40$ &$ [0, 1, 3, 7, 9, 13, 17, 19, 21, 23, 27, 29, $  & $  [1, 2, 4, 2, 4, 4, 2, 2, 2, 4, 2, $ & $ [0, 1, 0, 1, 1, 0, 0, 0, 1, 0, 0, 2, 0, 0,  $\\
&$5\sspcdot 11$& &$29, 31, 37, 39, 41, 43, 47, 49, 51, 53, 57]   $ & $  2, 6, 2, 2, 2, 4, 2, 2, 2, 4]$& $0, 1,0, 0, 0, 1, 0, 0, 0, 1, 0, 0, 0, 2, 0,   $  \\
&& &  &  & $ 0, 1, 0, 0, 0, 1, 1, 0, 1, 0, 0]   $  \\
 \hline
$\bf 8$  & $\bf 57$ &$36$ &$[0, 1, 5, 7, 11, 13, 17, 23, 25, 29, 31, $  & $[1, 4, 2, 4, 2, 4, 6, 2, 4,$ & $ [1, 0, 1, 0, 1, 2, 0, 1, 0, 1, 0, 1, 0, 1, 0, $\\
&$3\sspcdot 19$& &$31, 35, 37, 41, 43, 47, 49, 53, 55, 59] $ & $2, 4, 2, 4, 2, 4, 2, 4, 2, 4]$& $ 1, 0, 1, 0, 1, 0, 1, 0, 1, 0, 1, 0, 1, 0, $  \\
&& &  &  & $  2, 1, 0, 1, 0, 1, 0]  $  \\
\hline
$\bf 9$  & $\bf 65$ & $48$ &$[0, 1, 3, 7, 9, 11, 17, 19, 21, 23, 27, 29, 31, 33, $  & $[1, 2, 4, 2, 2, 6, 2, 2, 2, 4, 2, 2,    $ & $[0, 1, 0, 0, 2, 0, 0, 0, 1, 0, 0, 0, 1, 1, 0, 1,   $\\
&$5\sspcdot 13 $& &$ 33, 37, 41, 43, 47, 49, 51, 53, 57, 59, 61, 63, 67]  $ & $  2, 4, 4, 2, 4, 2, 2, 2, 4, 2, 2, 2, 4] $& $ 0, 0, 0, 1, 0, 0, 0, 1, 0, 0, 0, 1, 0, 0, 0, 1,   $  \\
&& &  &  & $ 0, 1, 1, 0, 0, 0, 1, 0, 0, 0, 2, 0, 0, 1, 0, 0]   $  \\
 \hline
$\bf 10$  & $\bf 69$ &$44$  &$ [0, 1, 5, 7, 11, 13, 17, 19, 25, 29, 31, 35, 37,  $  & $ [1, 4, 2, 4, 2, 4, 2, 6, 4, 2, 4, 2, $ & $[1, 0, 1, 0, 1, 0, 2, 1, 0, 1, 0, 1, 0, 1, 0,  $\\
&$3\sspcdot 23$  &  &$41, 43, 47, 49, 53, 55, 59, 61, 65, 67, 71]  $ & $4, 2, 4, 2, 4, 2, 4, 2, 4, 2, 4]   $& $1, 0, 1, 0, 1, 0, 1, 0, 1, 0, 1, 0, 1, 0, 1,  $  \\
&& &$  $ & $   $& $ 0, 1, 0, 1, 0, 1, 2, 0, 1, 0, 1, 0, 1, 0]   $  \\
 \hline
$\bf 11$   & $\bf 77$ & $60$  &$[0, 1, 3, 5, 9, 13, 15, 17, 19, 23, 25, 27,   $  & $[1, 2, 2, 4, 4, 2, 2, 2, 4, 2, 2, 2, 2, 6, 2, 2,  $ & $ [0, 0, 1, 1, 0, 0, 0, 1, 0, 0, 0, 0, 2, 0, 0,       $\\
& $7\sspcdot 11$& &$29, 31, 37, 39, 41, 43, 45, 47, 51, 53,    $ & $ 2, 2, 2, 4, 2, 4, 2, 2, 4, 2, 2, 2, 2, 2, 4]   $& $ 0, 0, 0, 1, 0, 1, 0, 0, 1, 0, 0, 0, 0, 0, 1,   $  \\
&& & $57, 59, 61, 65, 67, 69, 71, 73, 75, 79] $ & $   $& $0, 0, 0, 0, 0, 1, 0, 0, 1, 0, 1, 0, 0, 0, 0,   $  \\
&& &  &  & $ 0, 2, 0, 0, 0, 0, 1, 0, 0, 0, 1, 1, 0, 0, 0]  $  \\
\hline \pbn
{\bf continued}
\end{tabular}
\end{center}
\end{landscape}
\vfill
\eject
\begin{landscape}
\begin{center}
{\large {\bf Table 5 ctnd.: Extended set $\bf \widehat{\mathbfcal M}(n)$, first differences $\bf {\boldsymbol \triangle}\,\widehat{\mathbfcal  M}(n) $, and floor-pattern $\bf{\mathbfcal F}(n)$ for composed odd squarefree modulus $\bf n$.}} 
\end {center}
\begin{center}  
%{\large  $\bf{a }$ }
%\end{center}
%\begin{center}
\begin{tabular}{|c|r|c|c|c|c|}\hline
&& && &\\
$\bf m$ &  $\bf n(m) $ & $\bf 2\,{\boldsymbol \delta}($n\bf $)$   &$\bf \widehat{\mathbfcal M}(n)$ & $\bf \triangle\,\widehat{\mathbfcal  M}(n) $ & $\bf{\mathbfcal F}(n)$ \\
&& && & \\ \hline\hline
$\bf 12$  & $\bf 85$  & $64$ &$[0, 1, 3, 7, 9, 11, 13, 19, 21, 23, 27, 29, 31,  $  & $ [1, 2, 4, 2, 2, 2, 6, 2, 2, 4, 2,       $ & $ [0, 1, 0, 0, 0, 2, 0, 0, 1, 0, 0, 0, 1, 0, 0, 0,     $\\
&$5\sspcdot 17$ & & $ 33, 37, 39, 41, 43, 47, 49, 53, 57, 59, 61,       $ & $2, 2, 4, 2, 2, 2, 4, 2, 4, 4, 2,        $& $ 1, 0, 1, 1, 0, 0, 0, 1, 0, 0, 0, 1, 0, 0, 0, 1,       $  \\
&& & $ 63, 67, 69, 71, 73, 77, 79, 81, 83, 87]    $ & $ 2, 2, 4, 2, 2, 2, 4, 2, 2, 2, 4]  $& $ 0, 0, 0, 1, 0, 0, 0, 1, 0, 0, 0, 1, 1, 0, 1, 0, $  \\
&& & $     $ & $    $& $0, 0, 1, 0, 0, 0, 1, 0, 0, 2, 0, 0, 0, 1, 0, 0] $  \\
$\bf 13$   & $\bf 87$ & $ 56$ &$ [0, 1, 5, 7, 11, 13, 17, 19, 23, 25, 31$  & $[1, 4, 2, 4, 2, 4, 2, 4, 2, 6, 4, 2, 4, 2, 4, $ &$[1, 0, 1, 0, 1, 0, 1, 0, 2, 1, 0, 1, 0, 1,   $\\
&$3\sspcdot 29$&&$35, 37, 41, 43, 47, 49, 53, 55, 59, 61,   $ & $ 2, 4, 2, 4, 2, 4, 2, 4, 2, 4, 2, 4, 2, 4]   $& $ 0, 1, 0, 1, 0, 1, 0, 1, 0, 1, 0, 1, 0, 1      $  \\ 
&& &$ 65, 67, 71, 73, 77, 79, 83, 85, 89]  $ & $    $& $ 0, 1, 0, 1, 0, 1, 0, 1, 0, 1, 0, 1, 0, 1,       $  \\ 
&& &$   $ & $    $& $ 0, 1, 0, 1, 2, 0, 1, 0, 1, 0, 1, 0, 1, 0]  $  \\ 
\hline
$\bf 14$  & $\bf 91$ & $72$ &$[0, 1, 3, 5, 9, 11, 15, 17, 19, 23, 25,$  & $ [1, 2, 2, 4, 2, 4, 2, 2, 4, 2, 2, 2, 2,  $ & $ [0, 0, 1, 0, 1, 0, 0, 1, 0, 0, 0, 0, 0, 1, 1,   $\\
&$7\sspcdot 13$& & $27, 29, 31, 33, 37, 41, 43, 45, 47,    $ & $ 2, 4, 4, 2, 2, 2, 4, 2, 2, 2, 2, 2, 6,    $& $0, 0, 0, 1, 0, 0, 0, 0, 0, 2, 0, 0, 0, 0, 1,  $  \\ 
&& &$ 51, 53, 55, 57, 59, 61, 67, 69, 71,   $ & $ 2, 2, 2, 2, 4, 2, 2, 2, 2, 2, 4]   $& $ 0, 0, 0, 0, 0, 1, 0, 0, 0, 0, 0, 1, 0, 0, 0,     $  \\ 
&& &$ 73, 75, 79, 81, 83, 85, 87, 89, 93]  $ & $   $& $0, 0, 0, 1, 0, 0, 1, 0, 1, 0, 0, 0]   $  \\  
\hline
$\bf 15$  & $\bf 93$ & $60$ &$ [0, 1, 5, 7, 11, 13, 17, 19, 23, 25, 29, 35,  $ & $ [1, 4, 2, 4, 2, 4, 2, 4, 2, 4, 6,  $ & $ [1, 0, 1, 0, 1, 0, 1, 0, 1, 2, 0, 1, 0, 1, 0,        $ \\
&$3\sspcdot 31$&  &$37, 41, 43, 47, 49, 53, 55, 59, 61, 65, 67,    $ & $ 2, 4, 2, 4, 2, 4, 2, 4, 2, 4, 2,   $& $ 1, 0, 1, 0, 1, 0, 1, 0, 1, 0, 1, 0, 1, 0, 1,   $  \\ 
&& &$ 71, 73, 77, 79, 83, 85, 89, 91, 95]  $ & $ 4, 2, 4, 2, 4, 2, 4, 2, 4]   $& $0, 1, 0, 1, 0, 1, 0, 1, 0, 1, 0, 1, 0, 1, 0,    $  \\ 
&& &$  $ & $   $& $ 1, 0, 1, 0, 2, 1, 0, 1, 0, 1, 0, 1, 0, 1, 0]  $  \\ 
\hline
$\bf 16$  & $\bf 95$ & $72$  &$[0, 1, 3, 7, 9, 11, 13, 17, 21, 23, 27,  $  & $ [1, 2, 4, 2, 2, 2, 4, 4, 2, 4, 2, 2, 2,   $   & $  [0, 1, 0, 0, 0, 1, 1, 0, 1, 0, 0, 0, 1, 0, 0,        $\\
&$5\sspcdot 19 $& &$29, 31, 33, 37, 39, 41, 43, 47, 49,   $ & $ 4, 2, 2, 2, 4, 2, 2, 2, 6, 2, 2, 4, 2,   $& $ 0, 1, 0, 0, 0, 2, 0, 0, 1, 0, 0, 0, 1, 0, 0,   $  \\ 
&& &$51, 53, 59, 61, 63, 67, 69, 71, 73,   $ & $ 2, 2, 4, 2, 2, 2, 4, 2, 2, 2, 4]  $& $ 0, 1, 0, 0, 0, 1, 0, 0, 0, 1, 0, 0, 0, 1, 0,   $  \\ 
&& &$ 77, 79, 81, 83, 87, 89, 91, 93, 97] $ & $   $& $ 0, 0, 1, 0, 0, 2, 0, 0, 0, 1, 0, 0, 0, 1, 0,    $  \\ 
&& &$  $ & $   $& $ 0, 0, 1, 0, 1, 1, 0, 0, 0, 1, 0, 0 ]   $  \\ 
\hline
$\bf 17$  & $\bf 105$ & $48$  &$[0, 1, 11, 13, 17, 19, 23, 29, 31, 37,      $  & $   [1, 10, 2, 4, 2, 4, 6, 2, 6,      $ & $  [4, 0, 1, 0, 1, 2, 0, 2, 1, 0, 1, 2,    $\\
&$3\sspcdot 5\sspcdot 7$ &  &$41, 43, 47, 53, 59, 61, 67, 71, 73,   $ & $4, 2, 4, 6, 6, 2, 6, 4, 2,    $& $ 2, 0, 2, 1, 0, 2, 1, 2, 3, 1, 0, 1,  $  \\ 
&& &$79, 83, 89, 97, 101, 103, 107 ]  $ & $6, 4, 6, 8, 4, 2, 4]   $& $  0, 1, 3, 2, 1, 2, 0, 1, 2, 0, 2, 2,   $  \\ 
&& &$  $ & $   $& $  1, 0, 1, 2, 0, 2, 1, 0, 1, 0, 4, 0]   $  \\
\hline
$\hskip .2cm\bf \vdots$&& &&& \\
\hline
\end{tabular}
\end{center}
\end{landscape}
\vfill
\eject
\noindent
\begin{center}
{\large {\bf Table 6: Cycle structure of \dstyle{\bf {\cal G}{\it al}(\boldsymbol{\mathbb Q}(\boldsymbol{\rho}(n))/\boldsymbol{ \mathbb Q})} {\bf for}\  $\bf n\sspeq 1,2,...,40$.}}
\end {center}
\begin{center}  
\begin{tabular}{|l|l|}\hline
$\hskip .2cm\bf n$ & \hskip 6cm  \bf cycles \\ 
\hline\hline
$\hskip .2cm\bf 1$    &  $[[0]]\sspeq [[1]]$   \\        
\hline
$\hskip .2cm\bf 2$    &  $ [[1]]$      \\ 
\hline
$\hskip .2cm\bf 3$    &  $ [[1]] $    \\    
\hline
 $\hskip .2cm\bf 4$    &  $ [[{\bf 3}, 1]] $  \\
\hline
$\hskip .2cm\bf 5$    & $ [[{\bf 3}, 1]] $ \\    
\hline
$\hskip .2cm\bf 6$    & $  [[{\bf 5}, 1]] $      \\
\hline
$\hskip .2cm\bf 7$    & $ [[3, 5, 1]] $  \\     
\hline
$\hskip .2cm\bf 8$    &  $ [[3, {\bf 7}, 5, 1]]  $   \\
\hline
$\hskip .2cm\bf 9$    &$ [[5, 7, 1]] $ \\    
\hline
$\hskip .1cm\bf 10$  & $  [[3, {\bf 9}, 7, 1]]  $   \\
\hline
$\hskip .1cm\bf 11$  &$  [[3, 9, 5, 7, 1]]  $ \\  
\hline
 \fbox{{\color{Bittersweet} ${\bf 12}$}}  &  $   [[{\bf 5}, 1], [{\bf 7}, 1], [{\bf 11}, 1]]  $  \\            
\hline
$\hskip .1cm\bf 13$  & $ [[7, 3, {\bf 5}, 9, 11, 1]] $ \\       
\hline
$\hskip .1cm\bf 14$  & $ [[5, 3, {\bf 13}, 9, 11, 1]]  $        \\
\hline
$\hskip .1cm\bf 15$  & $ [[7, {\bf 11}, 13, 1]]    $\\  
\hline
$\hskip .1cm\bf 16$  & $ [[3, 9, 5, {\bf 15}, 13, 7, 11, 1]]  $ \\
\hline
$\hskip .1cm\bf 17$  & $ [[3, 9, 7, {\bf 13}, 5, 15, 11, 1]]  $  \\  
\hline
$\hskip .1cm\bf 18$  &  $  [[5, 11, {\bf 17}, 13, 7, 1]]   $ \\    
\hline
$\hskip .1cm\bf 19$  & $ [[3, 9, 11, 5, 15, 7, 17, 13, 1]]     $\\     
\hline
\fbox{{\color{Bittersweet} ${\bf 20}$}}  & $ [[3, {\bf 9}, 13, 1], [7, {\bf 9}, 17, 1], [{\bf 11}, 1], [{\bf 19}, 1]]   $\\     
\hline
$\hskip .1cm\bf 21$  &  $  [[11, 5, {\bf 13}, 17, 19, 1]]  $   \\
\hline
$\hskip .1cm\bf 22$  & $[[3, 9, 17, 7, {\bf 21}, 19, 13, 5, 15, 1]]   $ \\
\hline
$\bf 23$  &$   [[3, 9, 19, 11, 13, 7, 21, 17, 5, 15, 1]]  $ \\ 
\hline
\fbox{{\color{Bittersweet} ${\bf 24}$}}  &  $[[5, {\bf 23}, 19, 1], [{\bf 7}, 1], [11, {\bf 23}, 13, 1], [{\bf 17}, 1]] $    \\
\hline
$\hskip .1cm\bf 25$  &  $ [[3, 9, 23, 19, {\bf 7}, 21, 13, 11, 17, 1]]  $  \\
\hline
$\hskip .1cm\bf 26$  & $ [[7, 3, 21, 9, 11, {\bf 25}, 19, 23, 5, 17, 15, 1]]   $ \\   
\hline
$\bf 27$  &  $ [[5,25,17,23,7,19,13,11,1]] $   \\
\hline
\fbox{{\color{Bittersweet} ${\bf 28}$}}  &  $[[3, 9, {\bf 27}, 25, 19, 1], [5, 25, {\bf 13}, 9, 11, 1],[17, 9, {\bf 15}, 25, 23, 1]] $       \\
\hline
$\hskip .1cm\bf 29$  &  $  [[3, 9, 27, 23, 11, 25, {\bf 17}, 7, 21, 5, 15, 13, 19, 1]]         $\\
\hline
 \fbox{{\color{Bittersweet} ${\bf 30}$}}  &  $[[7, {\bf 11}, 17, 1], [13, {\bf 11}, 23, 1], [{\bf 19}, 1], [{\bf 29}, 1]] $\\ 
\hline 
$\hskip .1cm\bf 31$  &  $   [[3,9,27,19,5,15,17,11,29,25,13,23,7,21,1]]      $\\
\hline
$\hskip .1cm\bf 32$  &  $  [[3, 9, 27, 17, 13, 25, 11, {\bf 31}, 29, 23, 5, 15, 19, 7, 21, 1]] $  \\ 
\hline 
$\hskip .1cm\bf 33$  &  $  [[5, 25, 7, 31,{\bf 23}, 17, 19, 29, 13, 1]]      $\\
\hline
$\hskip .1cm\bf 34$  &  $  [[3, 9, 27, 13, 29, 19, 11, {\bf 33}, 31, 25, 7, 21, 5, 15, 23, 1]]   $  \\ 
\hline 
$\hskip .1cm\bf 35$  &  $ [[3, 9, 27, 11, 33, {\bf 29}, 17, 19, 13, 31, 23, 1]]      $\\
\hline
\fbox{{\color{Bittersweet} ${\bf 36}$}}   &  $  [[5, 25, {\bf 19}, 23, 29, 1], [7, 23, {\bf 17}, 25, 31, 1], [11, 23, {\bf 35}, 25, 13, 1]]     $  \\ 
\hline 
$\hskip .1cm\bf 37$  &  $  [[5, 25, 23, 33, 17, 11, 19, 21, {\bf 31}, 7, 35, 27, 13, 9, 29, 3, 15, 1]]   $\\
\hline
$\hskip .1cm\bf 38$  &  $ [[13, 17, 7, 15, 33, 27, 29, 3, {\bf 37}, 25, 21, 31, 23, 5, 11, 9, 35, 1]]   $  \\ 
\hline 
$\hskip .1cm\bf 39$  &  $  [[7, 29, 31, 17, 37, {\bf 25}, 19, 23, 5, 35, 11, 1]]   $\\
\hline
\fbox{{\color{Bittersweet} ${\bf 40}$}}   &  $  [[3, {\bf 9}, 27, 1], [7, {\bf 31}, 23, 1], [11, {\bf 39}, 29, 1], [13, {\bf 9}, 37, 1], [17, {\bf 31}, 33, 1], [19, {\bf 39}, 21, 1]]   $  \\ 
\hline 
$\hskip .2cm\bf \vdots$&\\
\hline
\end{tabular}
\end{center}
\psn
\begin{center}
\noindent
Boxed and colored $\bf n$-numbers indicate non-cyclic Galois groups. See Table 8.\pn
Boldface numbers are nontrivial square roots Modd\, $\bf n$, denoted by $\bf s$ in the text.
\end{center}
\vfill
\eject
\noindent
%\begin{landscape}
\begin{center}
{\large {\bf Table 7: Non-cyclic Galois groups \dstyle{\bf {\cal G}{\it al}(\boldsymbol{\mathbb Q}(\boldsymbol{\zeta}(n))/\boldsymbol{ \mathbb Q})}\, ,\,  $\bf n\sspleq 100$ }} 
\end {center}
\begin{center}  
%{\large  $\bf{ }$ }
%\end{center}
%\begin{center}
\begin{tabular}{|l|c|c|c|c|}\hline
&& &&\\
\hskip .2cm$\bf n$ & $\bf \boldsymbol{\varphi}(n)$  &\bf{cycle structure}& \bf{no. of cycles} & $\bf Galois\, group$ \\
&& && \\ \hline\hline
\fbox{{\color{Bittersweet}${\bf  8}$}}  &   $4$    &  $2_3$ &  $3$      & $ Z_2\ssptimes Z_2 $  \\
\hline
$\hskip .1cm\bf 12$  &   $4 $    & $2_3$ &   $3$     & $ Z_2\ssptimes Z_2 $  \\
\hline
\fbox{{\color{Bittersweet}${\bf 15}$}}  &  $8$    &   $4_2\,2_2$ &     $4$  &  $ Z_4\ssptimes Z_2 $       \\
\hline
$\hskip .1cm\bf 16$  &  $8$    &   $4_2\,2_2$ &     $4$  &  $ Z_4\ssptimes Z_2 $       \\
\hline
$\hskip .1cm\bf 20$  &  $8$    &   $4_2\,2_2$ &     $4$  &  $ Z_4\ssptimes Z_2 $       \\
\hline
\fbox{{\color{Bittersweet}${\bf 21}$}}  &  $12$ &   $ 6_3 $ &     $3$  &  $ Z_3\ssptimes Z_2^2 $  \\
\hline
\fbox{{\color{Bittersweet}${\bf 24}$}}  & $8$  &$2_7 $  & $7$ & $  Z_2^3 $\\
\hline
$\hskip .1cm{\bf 28}$  &  $12$ &   $ 6_3 $ &     $3$  &  $ Z_3\ssptimes Z_2^2 $  \\
\hline
$\hskip .1cm\bf 30$  &   $8$   &   $4_2\,2_2$   &   $4$  &  $Z_4\ssptimes Z_2$       \\
 \hline
\fbox{{\color{Bittersweet}${\bf 32}$}}  & $16$  &$8_2\,4_1\, 2_2$  & $5$ & $ Z_8\ssptimes Z_2 $\\
 \hline
\fbox{{\color{Bittersweet}${\bf 33}$}}   & $20$  &$10_3$  & $3$ & $ Z_5\ssptimes Z_2^2 $\\
\hline
\fbox{{\color{Bittersweet}${\bf 35}$}}   & $24$  &$12_2\,6_2$  & $4$ & $ Z_4\ssptimes Z_3 \ssptimes Z_2 $\\
\hline
$\hskip .1cm\bf 36$  &   $12$   &   $6_3$   &   $3$  &  $Z_3\ssptimes Z_2^2$       \\
 \hline
$ \hskip .1cm\bf 39$   & $24$  &$12_2\,6_2$  & $4$ & $ Z_4\ssptimes Z_3 \ssptimes Z_2 $\\
\hline
\fbox{{\color{Bittersweet}${\bf 40}$}}   & $16$  &$4_4\,2_6$  & $10$ & $ Z_4\ssptimes Z_2^2 $\\
 \hline
$\hskip .1cm\bf 42$  &   $12$   &   $6_3$   &   $3$  &  $Z_3\ssptimes Z_2^2$       \\
 \hline
$\hskip .1cm\bf 44$  &   $20$   &   $10_3$   &   $3$  &  $Z_5\ssptimes Z_2^2$       \\
 \hline
$ \hskip .1cm\bf 45$   & $24$  &$12_2\,6_2$  & $4$ & $ Z_4\ssptimes Z_3 \ssptimes Z_2 $\\
\hline
$\hskip .1cm \bf 48$   & $16$  &$4_4\,2_6$  & $10$ & $ Z_4\ssptimes Z_2^2 $\\
 \hline
\fbox{{\color{Bittersweet}${\bf 51}$}}  & $32$  &$16_2\,8_1\,4_1\, 2_2$  & $6$ & $ Z_{16}\ssptimes Z_2 $\\
 \hline
$ \hskip .1cm\bf 52$   & $24$  &$12_2\,6_2$  & $4$ & $ Z_4\ssptimes Z_3 \ssptimes Z_2 $\\
\hline
\fbox{{\color{Bittersweet}${\bf 55}$}}  & $40$  &$20_2\,10_2$  & $4$ & $ Z_5\ssptimes Z_4 \ssptimes Z_2 $\\
\hline
\fbox{{\color{Bittersweet}${\bf 56}$}}   & $24$  &$6_7$  & $7$ & $ Z_3\ssptimes Z_2^3 $\\
 \hline
\fbox{{\color{Bittersweet}${\bf 57}$}}   & $36$  &$18_3$  & $3$ & $ Z_9\ssptimes Z_2^2 $\\
 \hline
$\hskip .1cm \bf 60$   & $16$  &$4_4\,2_6$  & $10$ & $ Z_4\ssptimes Z_2^2 $\\
 \hline
\fbox{{\color{Bittersweet}${\bf 63}$}}   & $36$  &$6_{12}$  & $12$ & $ Z_3^2\ssptimes Z_2^2 $\\
 \hline
$\hskip .1cm \bf 64$ & $32$  &$16_2\,8_1\,4_1\, 2_2$  & $6$ & $ Z_{16}\ssptimes Z_2 $\\
 \hline
\fbox{{\color{Bittersweet}${\bf 65}$}}   & $48$  &$12_6$  & $6$ & $ Z_4^2\ssptimes Z_3 $\\
 \hline
$\hskip .1cm\bf 66$   & $20$  &$10_3$  & $3$ & $ Z_5\ssptimes Z_2^2 $\\
\hline
$\hskip .1cm \bf 68$ & $32$  &$16_2\,8_1\,4_1\, 2_2$  & $6$ & $ Z_{16}\ssptimes Z_2 $\\
 \hline
\fbox{{\color{Bittersweet}${\bf 69}$}}   & $44$  &$22_3$  & $3$ & $ Z_{11}\ssptimes Z_2^2 $\\
 \hline
$ \hskip .1cm\bf 70$   & $24$  &$12_2\,6_2$  & $4$ & $ Z_4\ssptimes Z_3 \ssptimes Z_2 $\\
 \hline
$\hskip .1cm\bf 72$   & $24$  &$6_7$  & $7$ & $ Z_3\ssptimes Z_2^3 $\\
 \hline
$\hskip .1cm\bf 75$ & $40$  &$20_2\,10_2$  & $4$ & $ Z_5\ssptimes Z_4 \ssptimes Z_2 $\\
\hline
$\hskip .1cm\bf  76$  & $36$  &$18_3$  & $3$ & $ Z_9\ssptimes Z_2^2 $\\
 \hline
\fbox{{\color{Bittersweet}${\bf 77}$}}   & $60$  &$30_3$  & $3$ & $ Z_5\ssptimes Z_3\ssptimes Z_2^2 $\\
 \hline
$ \hskip .1cm\bf 78$   & $24$  &$12_2\,6_2$  & $4$ & $ Z_4\ssptimes Z_3 \ssptimes Z_2 $\\
\hline
\fbox{{\color{Bittersweet}${\bf 80}$}}   & $32$  &$4_{12}\,2_4$  & $16$ & $ Z_4^2\ssptimes Z_2 $\\
\hline
\end{tabular}
\end{center}
\psn
\psn
\hskip 2cm Continued on the next page.
\vfill
\eject
\noindent
\begin{center}
{\large {\bf Table 7 continued: Non-cyclic Galois groups \dstyle{\bf {\cal G}{\it al}(\boldsymbol{\mathbb Q}(\boldsymbol{\zeta}(n))/\boldsymbol{ \mathbb Q})}\, ,\,  $\bf n\sspleq 100$ }} 
\end {center}
\begin{center}  
%{\large  $\bf{ }$ }
%\end{center}
%\begin{center}
\begin{tabular}{|l|c|c|c|c|}\hline
&& &&\\
\hskip .2cm$\bf n$ & $\bf \boldsymbol{\varphi}(n)$  &\bf{cycle structure}& \bf{no. of cycles} & $\bf Galois group$ \\
&& && \\ \hline\hline
$\hskip .1cm\bf 84$   & $24$  &$6_7$  & $7$ & $ Z_3\ssptimes Z_2^3 $\\
 \hline
\fbox{{\color{Bittersweet}${\bf 85}$}}   & $64$  &$16_4\,8_2\,4_4$  & $10$ & $ Z_{16}\ssptimes Z_4 $\\
\hline
\fbox{{\color{Bittersweet}${\bf 87}$}}   & $56$  &$28_2\,14_2$  & $4$ & $ Z_7\ssptimes Z_4 \ssptimes Z_2$\\
\hline
\fbox{{\color{Bittersweet}${\bf 88}$}}   & $40$  &$10_7$  & $7$ & $ Z_5\ssptimes Z_2^3 $\\
\hline
$ \hskip .1cm\bf 90$   & $24$  &$12_2\,6_2$  & $4$ & $ Z_4\ssptimes Z_3 \ssptimes Z_2 $\\
\hline
\fbox{{\color{Bittersweet}${\bf 91}$}}   & $72$  &$12_8\,6_8$  & $16$ & $ Z_4\ssptimes Z_3^2\ssptimes Z_2 $\\
\hline
\fbox{{\color{Bittersweet}${\bf 92}$}}   & $44$  &$22_3$  & $3$ & $ Z_{11}\ssptimes Z_2^2 $\\
\hline
$\hskip .1cm\bf 93$  & $60$  &$30_3$  & $3$ & $ Z_5\ssptimes Z_3\ssptimes Z_2^2 $\\
 \hline
\fbox{{\color{Bittersweet}${\bf 95}$}}   & $72$  &$36_2\,18_2$  & $4$ & $ Z_9\ssptimes Z_4\ssptimes Z_2 $\\
\hline
\fbox{{\color{Bittersweet}${\bf 96}$}}   & $32$  &$8_4\,4_3\,2_6$  & $13$ & $ Z_8\ssptimes Z_2^2 $\\
\hline
$\hskip .1cm\bf 99$  & $60$  &$30_3$  & $3$ & $ Z_5\ssptimes Z_3\ssptimes Z_2^2 $\\
 \hline
$\hskip .1cm\bf 100$ & $40$  &$20_2\,10_2$  & $4$ & $ Z_5\ssptimes Z_4 \ssptimes Z_2 $\\
\hline
$\hskip .2cm\bf \vdots$&& && \\
\hline
\end{tabular}
\end{center}
\psn
\psn
The cyclic group of order $\bf m$ is denoted by $\bf Z_m$. For all other values  $\bf n\sspleq 100$ the Galois group is the cyclic group $\bf Z_{{\boldsymbol \varphi}(n)}$.\psn
Only independent cycles are counted, i.e., cycles which appear as sub-cycles of the given ones have been omitted. \psn
The notation, e.g., $\bf 16_2\,8_1\,4_1\, 2_2$, means that there are $\bf 2$ cycles of order (length) $\bf 16$, one cycle of order $\bf 8$, one cycle of order $\bf 4$ and two cycles of order $\bf2$.\psn 
Direct products of identical cyclic groups are sometimes written in exponent form, e.g., $\bf Z_2^2$ stands for $\bf Z_2\ssptimes Z_2$. \psn
Boxed and colored $\bf n$-numbers indicate where some non-cyclic Galois group appears for the first time. Some of  the cycle graphs are shown in  Fig. 4.
%\end{landscape}
\vfill
\eject
\noindent
%\begin{landscape}
\begin{center}
{\large {\bf Table 8: Non-cyclic Galois groups \dstyle{\bf {\cal G}{\it al}(\boldsymbol{\mathbb Q}(\boldsymbol{\rho}(n))/\boldsymbol{ \mathbb Q})}\, ,\,  $\bf n\sspleq 100$ }} 
\end {center}
\begin{center}  
%{\large  $\bf{ }$ }
%\end{center}
%\begin{center}
\begin{tabular}{|l|c|c|c|c|}\hline
&& &&\\
\hskip .2cm$\bf n$ & $\bf \boldsymbol{\delta}(n)$  &\bf{cycle structure}& \bf{no. of cycles} & $\bf Galois group$ \\
&& && \\ \hline\hline
\fbox{{\color{Bittersweet}${\bf 12}$}}  &   $4$    &  $2_3$ &  $3 $      & $ Z_2\ssptimes Z_2 $  \\
\hline
\fbox{{\color{Bittersweet}${\bf 20}$}}  &   $8 $    & $4_2\,2_2$ &   $4$     & $ Z_4\ssptimes Z_2 $  \\
\hline
$\hskip .1cm\bf 24$  &  $8$    &   $4_2\,2_2$ &     $4$  &  $ Z_4\ssptimes Z_2 $       \\
\hline
\fbox{{\color{Bittersweet}${\bf 28}$}}  &  $12$ &   $ 6_3 $ &     $3$  &  $ Z_3\ssptimes Z_2^2 $  \\
\hline
$\hskip .1cm\bf 30$  &   $8$   &   $4_2\,2_2$   &   $4$  &  $Z_4\ssptimes Z_2$       \\
 \hline
$\hskip .1cm\bf 36$  & $12$  &$6_3$  & $3$ & $Z_3\ssptimes Z_2^2 $\\
 \hline
\fbox{{\color{Bittersweet}${\bf 40}$}}  & $16$  &$4_6$  & $6$ & $ Z_4\ssptimes Z_4 $\\
 \hline
$\hskip .1cm\bf 42$  & $12$  &$6_3$  & $3$ & $Z_3\ssptimes Z_2^2 $\\
 \hline
\fbox{{\color{Bittersweet}${\bf 44}$}}   & $20$  &$10_3$  & $3$ & $ Z_5\ssptimes Z_2^2 $\\
 \hline
\fbox{{\color{Bittersweet}${\bf 48}$}}  & $16$  &$8_2\,4_1\, 2_2$  & $5$ & $ Z_8\ssptimes Z_2 $\\
 \hline
\fbox{{\color{Bittersweet}${\bf 52}$}}   & $24$  &$12_2\,6_2$  & $4$ & $ Z_4\ssptimes Z_3 \ssptimes Z_2 $\\
\hline
$\hskip .1cm\bf 56$  & $24$  &$12_2\,6_2$  & $4$ & $ Z_4\ssptimes Z_3 \ssptimes Z_2 $\\
\hline 
\fbox{{\color{Bittersweet}${\bf 60}$}}   & $16$  &$4_4\,2_6$  & $10$ & $ Z_4\ssptimes Z_2^2 $\\
 \hline
\fbox{{\color{Bittersweet}${\bf 63}$}}  & $18$  &$6_4$  & $4$ & $ Z_3^2\ssptimes Z_2$\\
 \hline
$\hskip .1cm\bf 65$  & $24$  &$12_2\,6_2$  & $4$ & $ Z_4\ssptimes Z_3 \ssptimes Z_2 $\\
\hline
$\hskip .1cm\bf 66$  & $20$  &$10_3$  & $3$ & $ Z_5\ssptimes Z_2^2 $\\
 \hline
\fbox{{\color{Bittersweet}${\bf 68}$}}  & $32$  &$16_2\,8_1\,4_1\, 2_2$  & $6$ & $ Z_{16}\ssptimes Z_2 $\\
 \hline
$\hskip .1cm\bf 70$  & $24$  &$12_2\,6_2$  & $4$ & $ Z_4\ssptimes Z_3 \ssptimes Z_2 $\\
\hline
$\hskip .1cm\bf 72$  & $24$  &$12_2\,6_2$  & $4$ & $ Z_4\ssptimes Z_3 \ssptimes Z_2 $\\
\hline
\fbox{{\color{Bittersweet}${\bf 76}$}}   & $36$  &$18_3$  & $3$ & $ Z_9\ssptimes Z_2^2 $\\
 \hline
$\hskip .1cm\bf 78$  & $24$  &$12_2\,6_2$  & $4$ & $ Z_4\ssptimes Z_3 \ssptimes Z_2 $\\
\hline
\fbox{{\color{Bittersweet}${\bf 80}$}}  & $32$  &$8_4\,4_4$  & $8$ & $ Z_8\ssptimes Z_4 $\\
 \hline
\fbox{{\color{Bittersweet}${\bf 84}$}}   & $24$  &$6_7$  & $7$ & $ Z_3\ssptimes Z_2^3 $\\
 \hline
$\hskip .1cm\bf 85$  & $32$  &$16_2\,8_1\,4_1\, 2_2$  & $6$ & $ Z_{16}\ssptimes Z_2 $\\
 \hline
\fbox{{\color{Bittersweet}${\bf 88}$}}  & $40$  &$20_2\,10_2$  & $4$ & $ Z_5\ssptimes Z_4 \ssptimes Z_2 $\\
\hline
$\hskip .1cm\bf 90$  & $24$  &$12_2\,6_2$  & $4$ & $ Z_4\ssptimes Z_3 \ssptimes Z_2 $\\
\hline
\fbox{{\color{Bittersweet}${\bf 91}$}}   & $36$  &$12_4$  & $4$ & $ Z_4\ssptimes Z_3^2 $\\
 \hline
$\hskip .1cm\bf 92$  & $44$  &$22_3$  & $3$ & $ Z_{11}\ssptimes Z_2^2 $\\
 \hline
$\hskip .1cm\bf 96$  & $32$  &$16_2\,8_1\,4_1\, 2_2$  & $6$ & $ Z_{16}\ssptimes Z_2 $\\
 \hline
$\bf 100$  & $40$  &$20_2\,10_2$  & $4$ & $ Z_5\ssptimes Z_4 \ssptimes Z_2 $\\
\hline
$\hskip .2cm\bf \vdots$&& && \\
\hline
\end{tabular}
\end{center}
\psn
\psn
The cyclic group of order $\bf m$ is denoted by $\bf Z_m$. For all other $\bf n\sspleq 100$ cases the Galois group is the cyclic group $\bf Z_{{\boldsymbol \delta}(n)}$. The $\bf n$ values are given in \seqnum{A206552}.\psn
Only independent cycles are counted, i.e., cycles which appear as sub-cycles of the given ones have been omitted.\psn
 The notation, e.g., $\bf 16_2\,8_1\,4_1\, 2_2$, means that there are $\bf 2$ cycles of order (length) $\bf 16$, one cycle of order $\bf 8$, one cycle of order $\bf 4$ and two cycles of order $\bf2$.\psn 
Direct products of identical cyclic groups are sometimes written in exponent form, e.g., $\bf Z_2^2$ stands for $\bf Z_2\ssptimes Z_2$. \psn
Boxed and colored $\bf n$-numbers indicate where some non-cyclic Galois group appears for the first time. For the cycle graphs see Fig. 4.
%\end{landscape}
\vfill
\eject
\vfill
\eject
%%%  Figure 4 
\parbox{16cm}{
{\includegraphics[height=22cm,width=1\linewidth]{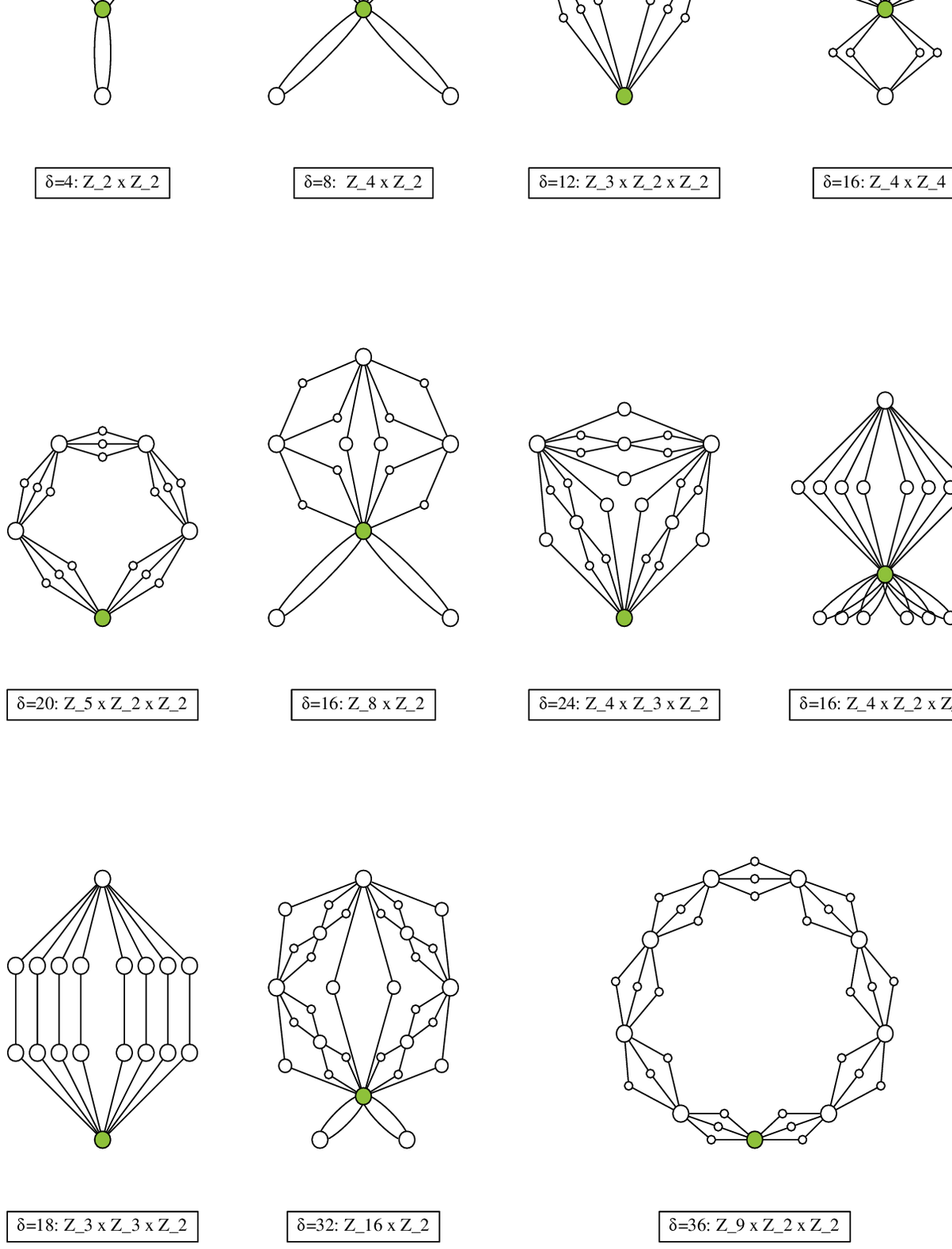}}
}
\psn
Figure 4:  Cycle graphs for non-cyclic Galois groups \dstyle{\bf {\cal G}{\it al}(\boldsymbol{\mathbb Q}(\boldsymbol{\rho}(n))/\boldsymbol{ \mathbb Q})} appearing for  $\bf n\sspeq 1..100$. $\delta$ is the degree of $\bf C(n,x)$, the maximal polynomial of  $\boldsymbol{\rho}(n)$, hence the order of the Galois group. See  Table 8. Continued on next page.
\psn
\vfill
\eject
%%%  Figure 4 cntd.
\parbox{16cm}{
{\includegraphics[height=15cm,width=1\linewidth]{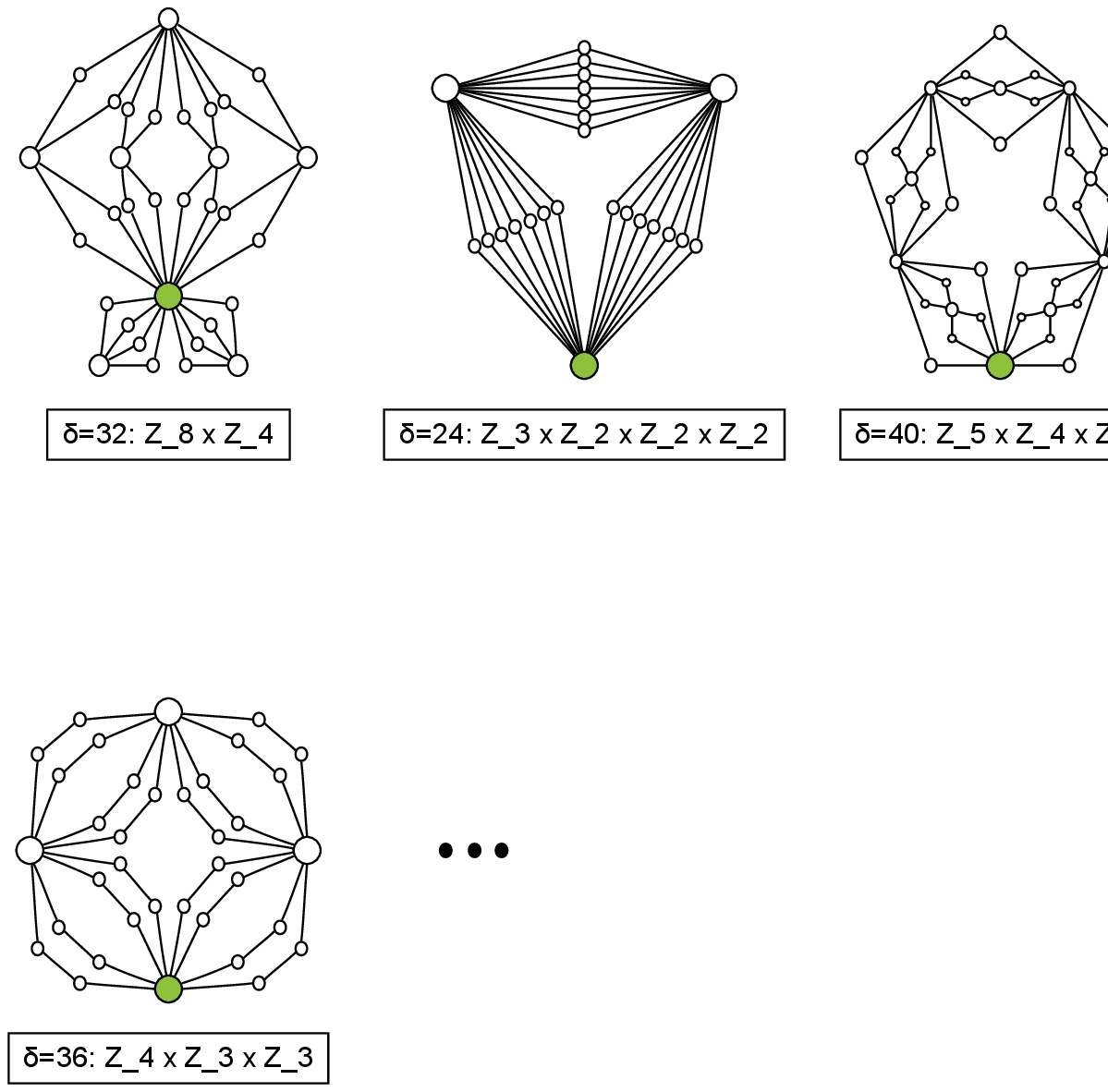}}
}
\psn
Figure 4 continued:  Cycle graphs for non-cyclic Galois groups \dstyle{\bf {\cal G}{\it al}(\boldsymbol{\mathbb Q}(\boldsymbol{\rho}(n))/\boldsymbol{ \mathbb Q})} appearing for  $\bf n\sspeq 1..100$. See  Table 8. 
%%%%%%%%%%%%%%%%%%%%%%
\vfill
\eject

\end{document}